\documentclass[onefignum,onetabnum]{siamart190516}

\usepackage{lipsum}
\usepackage{amsfonts}
\usepackage{amsmath}
\usepackage{amssymb}
\usepackage{graphicx}
\usepackage{epstopdf}
\usepackage{algorithmic}
\ifpdf
  \DeclareGraphicsExtensions{.eps,.pdf,.png,.jpg}
\else
  \DeclareGraphicsExtensions{.eps}
\fi

\newsiamremark{example}{Example}
\newsiamremark{remark}{Remark}
\newsiamremark{hypothesis}{Hypothesis}
\crefname{hypothesis}{Hypothesis}{Hypotheses}
\newsiamthm{claim}{Claim}

\headers{Predicting bifurcations of almost-invariant patterns}{M. Ndour, K. Padberg-Gehle}

\title{Predicting bifurcations of almost-invariant patterns: a set-oriented approach\thanks{Submitted December 24, 2019.
\funding{This research has received funding from the European Union's Horizon 2020 Research and Innovation Programme under the Marie Sklodowska-Curie grant agreement no. 643073, and from the Deutsche Forschungsgemeinschaft within Priority Programme SPP 1881 Turbulent Superstructures.}}}

\author{Moussa Ndour\thanks{Technische Universit\"at Dresden,
Institute of Scientific Computing, 
Faculty of Mathematics, 
D-01062 Dresden, Germany  (\email{moudesa.ndour@mailbox.tu-dresden.de})}
\and Kathrin Padberg-Gehle\thanks{Leuphana Universit\"at L\"uneburg,
Institute of Mathematics and its Didactics, 
Faculty of Education,
Universit\"atsallee 1,
D-21335 L\"uneburg, Germany 
  (\email{padberg@leuphana.de}).}}

\usepackage{amsopn}

\begin{document}

\maketitle

\begin{abstract}
 The understanding and prediction of sudden changes in flow patterns is of paramount importance in the analysis of geophysical flows as these rare events relate to critical phenomena such as atmospheric blocking, the weakening of the Gulf stream, or the splitting of the polar vortex. In this work our aim is to develop first steps towards a theoretical understanding of vortex splitting phenomena.
To this end, we study bifurcations of global flow patterns in parameter-dependent two-dimensional incompressible flows, with the flow patterns of interest corresponding to specific invariant sets. Under small random perturbations these sets become almost-invariant and can
be computed and studied by means of a set-oriented approach, where the underlying dynamics is described in terms of a reversible finite-state Markov chain. Almost-invariant sets are obtained from the sign structure of leading eigenvectors of the corresponding transition matrix. By a flow pattern bifurcation we mean a qualitative change in the form of a break-up of an almost-invariant set, when a critical external parameter of the underlying dynamical system is reached. For different examples and settings we follow the spectrum and the corresponding eigenvectors under continuous changes of the underlying system and yield indicators for different bifurcation scenarios for almost-invariant sets. In particular, we study a Duffing-type oscillator, which is known to undergo a classic pitchfork bifurcation. We find that the set-oriented analogue of this classical bifurcation includes a splitting of a rotating pattern, which has generic precursor signal that can be deduced from the behavior of the spectrum.
\end{abstract}

\begin{keywords}
 almost-invariant sets, transfer operator, bifurcation, Hamiltonian system, set-oriented approach 
 \end{keywords}

\begin{AMS}
	37J20,37M20,37M25
  \end{AMS}

\section{Introduction}
Understanding critical transitions in the macroscopic dynamics of a given complex system is, nowadays, of high interest due to the 
emergence of a new scientific challenge towards developing mathematical formulations of bifurcations in complex system models \cite{Mscheffer,MschefferandCo}. A concrete case study may consist of finding indicators or early warning signals of sudden changes in flow patterns emerging from real world systems. 
Famous examples include the Antarctic polar vortex break up scenario in late September 2002, where the rotating atmospheric pattern 
suddenly underwent a radical split (see e.g.\ \cite{breakup,split,LekRoss}).
One may classify this splitting event in the range of critical transitions  in real world complex flows. 
Thus, a legitimate question arises: How could such an event be predicted before its occurrence?

A possible answer to this question can be eventually made by combining a model-based approach and its set-oriented investigation. 
That is, first, one needs to find a simple but relevant mathematical model whose dynamics is sufficiently representative of the underlying real world complex system phenomena. 
Second, from the newly built model, one needs to computationally extract  patterns that, a priori, exhibit polar vortex-like dynamics. 
Finally, from the dynamics of the chosen model, it follows that the task of finding 
early warning signals of any radical split of the resulting pattern will progressively lead to predicting the sudden change.

A set-oriented dynamical systems approach aims at finding particular measurable sets and studying the probability of transport between them under the evolution of the system. 
These sets should be robust under small random external perturbations. In this case, mixing is equivalent to transport through the boundaries of the sets. 
That is, the set-oriented approach is a probabilistic method of computing slowly mixing sets.
A
 mathematical representation of slowly mixing patterns in dynamical systems was broadly studied recently; see e.g.\ \cite{DellOJ,GFMDell,GFnonauto1,GFnonauto2,GFKPG}. 
These studies were built around the idea of finding measurable partitions of the phase space of a given dynamical system in terms of phase space regions 
with minimal transport through their boundaries. For an autonomous dynamical system, those regions are referred to as \textit{almost-invariant sets} \cite{DellOJ, GFMDell}, 
since they mitigate transport between their interior 
and the rest of the phase space. They are called \textit{coherent sets} in the context of nonautonomous systems, as they move over finite time intervals with minimal 
dispersion \cite{GFnonauto1, GFnonauto2}. 

In this work, we aim to provide first steps towards predicting bifurcation of patterns that can mathematically be represented as optimal almost-invariant sets. 
Our approach is probabilistic and will be mainly based on analysing the spectral behavior of discrete Markov chains subject to external perturbations. The corresponding stochastic transition matrices are finite rank approximations of the Perron-Frobenius operator and its diffused version. 
Almost-invariant sets are numerically approximated by means of the dominant eigenvector basis of the transition matrix. Moreover, we use the sign structures of these vectors to systematically 
design the meaningful patterns that emerge from the dynamical system under study. Meanwhile, the behavior of the corresponding dominant eigenvalues under parameter-variation indicates when radical changes of patterns occur.
Indeed, as the bifurcation parameter is varied, eigenvalues change continuously with respect to the parameter. 

To the best of our knowledge, a set-oriented bifurcation analysis is still a broadly open topic that may require a new theoretical approach beyond classical bifurcation theory. 
Known previous works in this direction include using the discrete spectrum of the deterministic Perron-Frobenius operator (referred to as transfer operator)
generated by dissipative and non-dissipative systems.
In \cite{GasNic} a transfer operator based framework was successfully developed for studying the one-dimensional pitchfork normal form. Indeed, particular changes 
in the discrete spectrum of the transfer operator, including a clustering process towards $1$ as the critical parameter is reached, yield indicators of the pitchfork bifurcation.
In \cite{alexis} early warning indicators for transitions between atmospheric flow regimes were defined based on the transfer operator of a dissipative atmospheric model. 
In that work, the discrete spectrum of the transfer operator was initially used to approximate two isolated regimes as almost-invariant sets.
Closer to the setting of the present work, bifurcations of almost-invariant  and almost-cyclic sets in two-dimensional conservative systems and corresponding changes in the spectrum of the transition matrices were observed in \cite{Oliver,Grover} but not systematically studied. The aim of the present paper is to provide further methodological steps towards a better understanding of such global bifurcations. 

This work is organized as follows: In section \ref{sec1} we briefly review the concept of almost-invariant sets within a set-oriented numerical framework,  which yields stochastic transition matrices for reversible finite-state Markov chains.  In section \ref{pertubedMotiv}, we address results from the perturbation theory of stochastic matrices and show how their dominant
spectrum is suitable for estimating almost-invariant sets that originate from invariant structures of the unperturbed dynamics. In section \ref{expModels}, in order to illustrate the numerical framework, we discuss the practical computation of the dominant 
almost-invariant sets for a two-dimensional non-dissipative flow.
In section \ref{numExp}, we start the foremost step by
systematically experimenting Markov chain toy models undergoing bifurcations of specifically constructed patterns. In section \ref{bfmodels}, we rigorously study two explicit dynamical models, a Duffing-type oscillator and a single gyre flow, and identify early warning 
signals for splittings of patterns through the trends of eigenvalues with respect to a bifurcation parameter. The paper concludes with a discussion and outlook in section \ref{conclusion}.

\section{Set-oriented approach and almost-invariant patterns}\label{sec1}
Here, we will review the probabilistic approach of analyzing the global evolution of a given dynamical system.  The goal is to identify patterns that remain invariant or almost-invariant under the time evolution of the resulting transition matrix.

Let us consider a $p$-parametrized ordinary differential equation in the domain $M\subset \mathbb{R}^d$ and let us suppose that $p \in \mathbb{R}$ is a 
bifurcation parameter. 
\begin{equation}\label{Deteq1}
 \dot{x} = F(x,p)=:F_p(x).
\end{equation}
We fix $p$ and assume that the vector field $F_p:M \to \mathbb{R}^d$  is sufficiently smooth to guarantee the existence and uniqueness of 
solutions of \eqref{Deteq1}. Thus, there exists a flow map $S^t: M \to M$ such that for any given initial solution $x(0)=x_0$ and flow time $t \in \mathbb{R}$
\begin{equation}\label{flm}
 x_0\mapsto S^t(x_0)\in M,\,\,  x_0\in  M
\end{equation}
yields the solution of the system at time $t$ for the initial value $x_0=x(0)$. In this paper, \eqref{Deteq1} is supposed to model the evolution of a rotating incompressible steady fluid 
flow such as a vortex. 
\subsection{Almost-invariant sets}
The ultimate goal is to find a partition of the phase space $M$ into $k$ sets $\{A_1,A_2,\ldots,A_k\}$ such that, under the evolution the dynamics, the transport 
 between these sets is very unlikely. 
In other words, $S^{-t} A_i \approx A_i$, $i=1,2,\ldots,k$. Thus, to make this approximation more precise, we consider the measure space $(M,\Sigma, \mu)$ such that the probability measure $\mu$ is $S^t$-invariant, i.e. $\mu(A) = \mu(S^{-t}(A))\,\forall \,\,A\in \Sigma$, 
and is absolutely continuous with respect to the Lebesgue measure.
Note that in this work $\mu$ is simply (normalized) Lebesgue measure itself since the dynamical system \eqref{Deteq1} under study is assumed to be incompressible. We refer to the 
measurable sets $\{A_1,A_2,\ldots,A_k\}$ as patterns because often they are labeled as such in 
real world applications, such as ocean eddies or atmospheric vortices. Finally, this macroscopic approach is set-oriented in the sense that the trajectory of a single point 
matters less than the ensemble evolution of a swarm of points 
or a measurable set of points. Following \cite{GaryF}, the \textit{invariance ratio} of a set $A_i,\,i=1,2,\ldots,k$ is defined as
 \begin{equation}\label{invariance1}
 \rho_\mu(A_i) = \frac{\mu(A_i \cap S^{-t}(A_i))}{\mu(A_i)}.
 \end{equation}
This is interpreted as the probability of a point in $A_i$ to stay in $A_i$ under the mapping $S^t$ . Hence, any measurable invariant set $A$ satisfies
$ \rho_\mu(A) = 1$. $\{A_1,A_2,\ldots,A_k\}$ is a family of \textit{almost-invariant sets} that partitions the phase space $M$ if 
$M=\cup_{i=1}^k A_i$ and 
\begin{equation}\label{invariance}
 \rho_\mu(A_i) \approx 1\,\,\,\forall \,\,i=1,2,\ldots,k.
\end{equation}
Finding such a family of almost-invariant sets is intractable in practice. Instead one seeks optimal solutions of a relaxed problem based on the description of the dynamics in terms of a finite-state Markov chain  and its spectral properties. 
\subsection{Discretization and stochastic matrices}
Without loss of generality, we discretize the phase space $M$ to obtain a finite state space $\mathcal{S}=\big\{ B_1, B_2,\ldots, B_N \big\}$ such that $m(B_i) = m(B_j)$, $i,j\leq N$ and 
$M= \cup_{j=1}^N B_j$ with $m(B_i \cap B_j) = 0$; where $m$ denotes the phase space volume measure - a normalized Lebesgue 
measure on $M$. Besides, let us define the lumped finite state \cite{GaryF}
$$\mathcal{C}_N = \Big\{\ A\subset M\,:\, A= \bigcup_{j\in \mathcal{I}}  B_j,\,\, \mathcal{I} \subset \{1,2,\ldots,N\}      \Big\}.$$
The time evolution of the dynamical system on the discretized phase space yields the transition matrix
\begin{equation}\label{matauto}
 (P^{t}_N)_{ij} =\frac{m(B_i \cap S^{-t}(B_j))}{m(B_i)}.
 \end{equation}
Each $(i,j)$-th entry is the probability that a randomly selected point $x\in B_i$ has its image in $B_j$.
$P^{t}_N$ is a row stochastic matrix and is interpreted as the transition matrix associated with an $N$-state Markov chain over the 
finite states $\big\{ B_i\big\}_{i=1}^N$.
Note that $P^{t}_N$ is actually a finite rank approximation of the Perron-Frobenius operator \cite{SUlam}.\\
The resulting dynamics of the Markov chain may be thought of as the dynamics of $S^t$ with a small amount of bounded noise added.\\
The Markov chain \eqref{matauto} from \eqref{Deteq1} is not in general reversible. However, reversible transition matrices yield important spectral properties which 
are dynamically efficient in terms of checking how mass is transported in both forward and backward time. Moreover, as we are ultimately interested in the macroscopic dynamics of patterns such as the transport and critical 
transition of optimal almost-invariant patterns, it is more relevant to use a reversibilised Markov chain. The latter comes as straightforward transformation of \eqref{matauto} as
\begin{equation}\label{introReversible}
 Q= \frac{(L+P)}{2},
\end{equation} 
 where $L = \Big(\frac{\pi_j P_{ji}}{\pi_i}\Big)_{i,j=1}^{N}$ is the transition matrix of the reversed Markov chain and $P:=P^{t}_N$ is assumed to have a unique positive stationary 
 distribution $\pi=[\pi_1,\pi_2,\ldots,\pi_N]$ with $\pi P = \pi$; in our case it holds that $\pi_i = m(B_i)$ with $m$ being normalized Lebesgue measure. It follows the approximation of the invariance ratio as follows 
 \begin{equation}\label{invarianceApprox}
\begin{split}
 \rho^N_\mu(A) &= \frac{\sum_{i,j\in \mathcal{I}} \pi_i(Q^{t}_N)_{ij}}{\sum_{i\in \mathcal{I}} \pi_i},\,\, \mathcal{I} \subset \{1,2,\ldots,N\},\,\,\,\text{invariance ratio},\\
               &=1\,\,\,\, \text{if}\,\,\, A\,\,\, \text{is invariant},\\
               &\approx 1\,\,\,\, \text{if}\,\,\, A\,\,\, \text{is almost-invariant}.  
\end{split}
\end{equation}
$Q$ is a transition matrix as the weighted average of two transition matrices $P$ and $L.$ Moreover, $Q$ is reversible since it satisfies the so-called
 detailed balance condition, $\pi_jQ_{j i} = \pi_i Q_{ij}$. Further important properties of $Q$ include:
 \begin{enumerate}
\item $Q$ is diagonalized by a basis of $\pi$-orthogonal right eigenvectors. 
\item $Q$ has only real eigenvalues contained in $[-1,\,1]$. Moreover, for any given eigenvalue with a corresponding right eigenvector $x$, there is an associated left eigenvector $y$ such that
$y = D_N x$, where $D_N= \text{diag}([\pi_1,\pi_2,\ldots,\pi_N]).$
\item 
$Q$ is symmetric or self-adjoint with respect to the weighted Euclidean space $\langle \cdot , \cdot\rangle_\pi$ defined in $\mathbb{R}^N$ 
such that $\langle x, y\rangle_\pi = \sum_{i=1}^N x_i y_i \pi_i$, and two vectors $x,y$ are orthogonal if $\langle x, y\rangle_\pi =0$. 
\end{enumerate} Moreover, it is easy to 
verify that the adjoint of $L$ with respect to $\langle \cdot , \cdot\rangle_\pi$ is the transition matrix $P$. Therefore, $Q$ is just the average of two adjoint matrices. Besides, 
in terms of transport, $Q$ is checking how mass is transported in forward and backward at stationarity. For more details on Markov chains we refer to e.g.\ \cite{Seneta}.
In this work, we will use the reversibilized transition matrix to almost-invariant patterns and their bifurcations.

\section{Perturbed invariant patterns and spectral configurations}\label{pertubedMotiv}
In this section, we consider a $k$-state, $k>2$, reducible Markov chain which becomes irreducible when it is subjected to small perturbations.
Then we assume the existence of disjoint strongly connected lumped 
states $\{ A_i\}_{i=1}^k$ and their perturbed versions $\{ A_i(\epsilon)\}_{i=1}^k$, $\epsilon\in \mathbb{R}$. Thus, under some convenient reordering within the states, 
the unperturbed and perturbed Markov chains are respectively given by 
\begin{equation}\label{pmatrix}	
Q = 
 \begin{pmatrix}
  Q_1 & 0 & \cdots & 0 \\
  0 & Q_2 & \cdots & 0 \\
  \vdots  & \vdots  & \ddots & \vdots  \\
  0 & 0 & \cdots & Q_k 
 \end{pmatrix}, 
 Q(\epsilon) = 
 \begin{pmatrix}
  Q_1(\epsilon) & E_{12} & \cdots & E_{1k} \\
  E_{21} & Q_2(\epsilon) & \cdots & E_{2k} \\
  \vdots  & \vdots  & \ddots & \vdots  \\
  E_{k1} & E_{k2} & \cdots & Q_k(\epsilon) 
 \end{pmatrix}, k>2.
 \end{equation}
Every $Q_i ,\,i=1,\ldots,k$,   
is a primitive $n_i\times n_i$ reversible transition matrix over the ``cloud'' of states $A_i$. Moreover, 
due to the reducibility of $Q$, a system described by $Q$ will always  
stay in state $A_i$ once it is initialized in $A_i$. This means that the conditional transition probability to map to $A_j$ when in $A_i$, 
$w(A_j, A_i) = \frac{\sum_{i\in I, j\in J }\pi_i Q_{ij} }{\sum_{i\in I} \pi_i}$, is the Kronecker symbol $\delta_{ij},\,\,\, j=1,\ldots,k,$ . Besides, the matrix $Q$ has 
an eigenvalue $1$ of multiplicity $k$. One may think of the matrix $Q$ as the reversibilized of version \eqref{matauto} generated from \eqref{Deteq1}.

The transition matrices $Q(\epsilon)$, $\epsilon\in\mathbb{R}$ are, however, irreducible and the magnitude of the off-diagonal blocks $E_{ij}$ is very small relative to $1$ 
with respect to any chosen matrix norm. 
This implies, $w(A_j(\epsilon),A_i(\epsilon))\approx \delta_{ij}$, $i,j=1,\ldots,k $, and means that when the dynamical system enters $A_i$, 
it will stay in $A_i$ for a long time  with high probability before 
it leaves. The patterns $A_i(\epsilon)$ are referred to as almost-invariant patterns. Following the settings in \cite{Kato}, $Q(\epsilon)$ 
is considered as an operator-valued function of $\epsilon$, which is analytic in $E\subset\mathbb{R}$, $0\in E$. Thus it can 
be, in general, expressed as $Q(\epsilon) = Q(0) + \epsilon Q^{(1)}$, which is the first order Taylor expansion of $Q(\epsilon)$.  
As a consequence of this regularity condition, the eigenvalues of $Q(\epsilon)$ are 
continuous in $\epsilon$. From this continuity and the fact that the $Q_i(\epsilon)$ are nearly stochastic matrices \cite{Meyer}, we have that the 
spectrum of $Q(\epsilon)$ includes three parts:
\begin{enumerate}
\item[\textbf{(a)}] the Perron root $\lambda_1(\epsilon) =1$, 
\item[\textbf{(b)}] the set of $k-1$ non-unit eigenvalues, $\{\lambda_2(\epsilon),\ldots,\lambda_k(\epsilon)\}$ that are clustered near $1$.  
\item[\textbf{(c)}] the remaining part of the spectrum which is bounded away from $1$, for small $\epsilon$. 
\end{enumerate}
Throughout this work, we set the ordering $1= \lambda_1(\epsilon) > \lambda_2(\epsilon)\geq \ldots \geq \lambda_k(\epsilon)$. Note that this section is motivated by the fact that the class 
of models we consider in this study yield invariant sets 
in their dynamical evolutions. This means that the almost-invariant patterns will be just considered as perturbed invariant sets.
\begin{example}\label{exsect3}
To illustrate this setting, we consider a $60$-state Markov chain with $\mathcal{S} = \{1,2,\ldots,60\}$. This is chosen to be reducible with three disjoint invariant 
patterns $A_1=\{1,2,\ldots,10\}$, $A_2=\{11,12,\ldots,40\}$ and $A_3=\{41,42,\ldots,60\}$. The corresponding transition matrix is shown in figure \ref{fig:tmatrices} (left) with the blue dots highlighting the nonzero entries. 
An example of a perturbed Markov chain, allowing for small amounts of transport between the three patterns, is shown in figure \ref{fig:tmatrices} (right), as the corresponding irreducible transition matrix $Q(\epsilon)$.
 \begin{figure}[!htb]
 \begin{center}
 \begin{tabular}{cc}
    \includegraphics[width=0.3\textwidth]{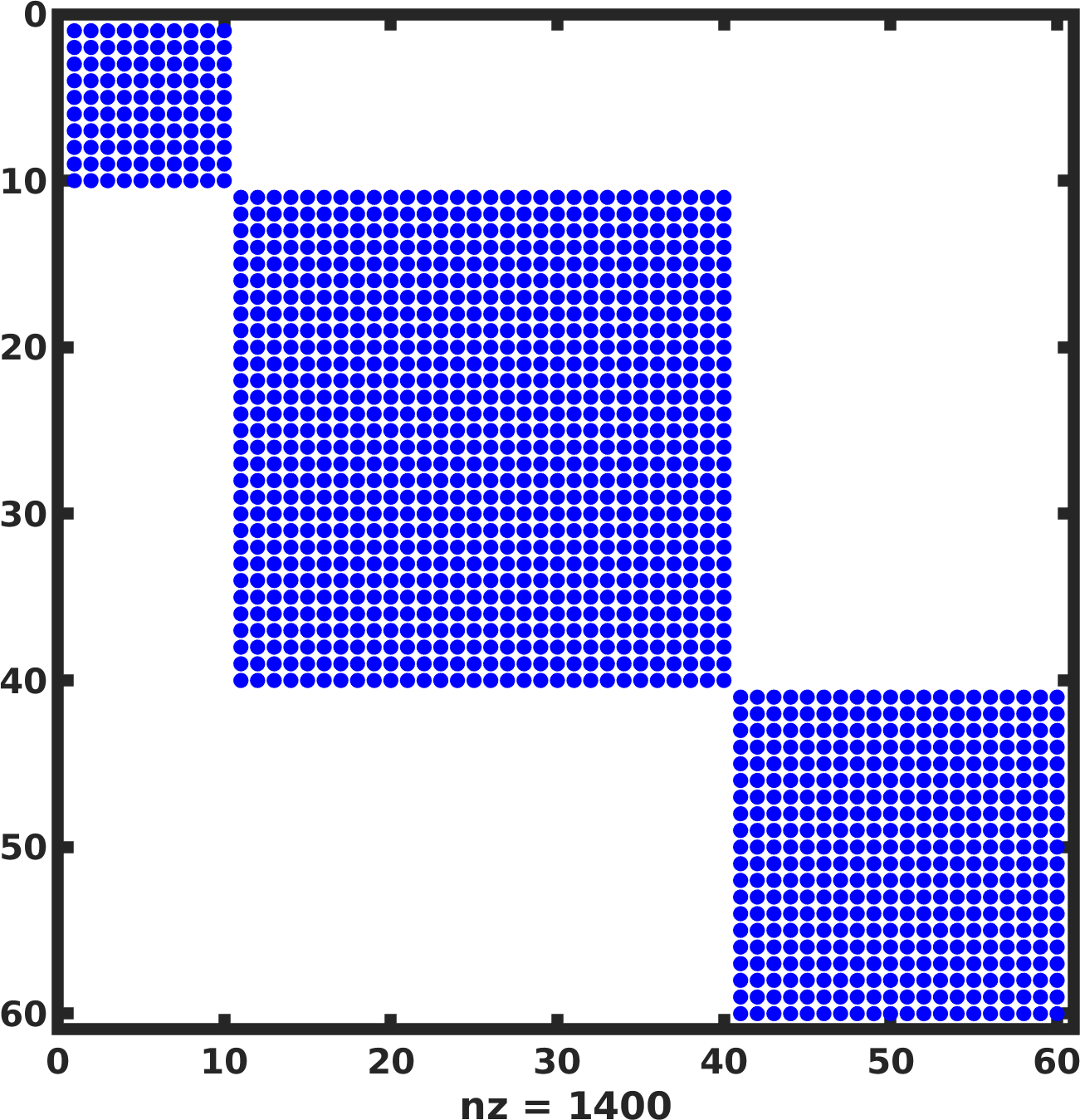} &
    \includegraphics[width=0.3\textwidth]{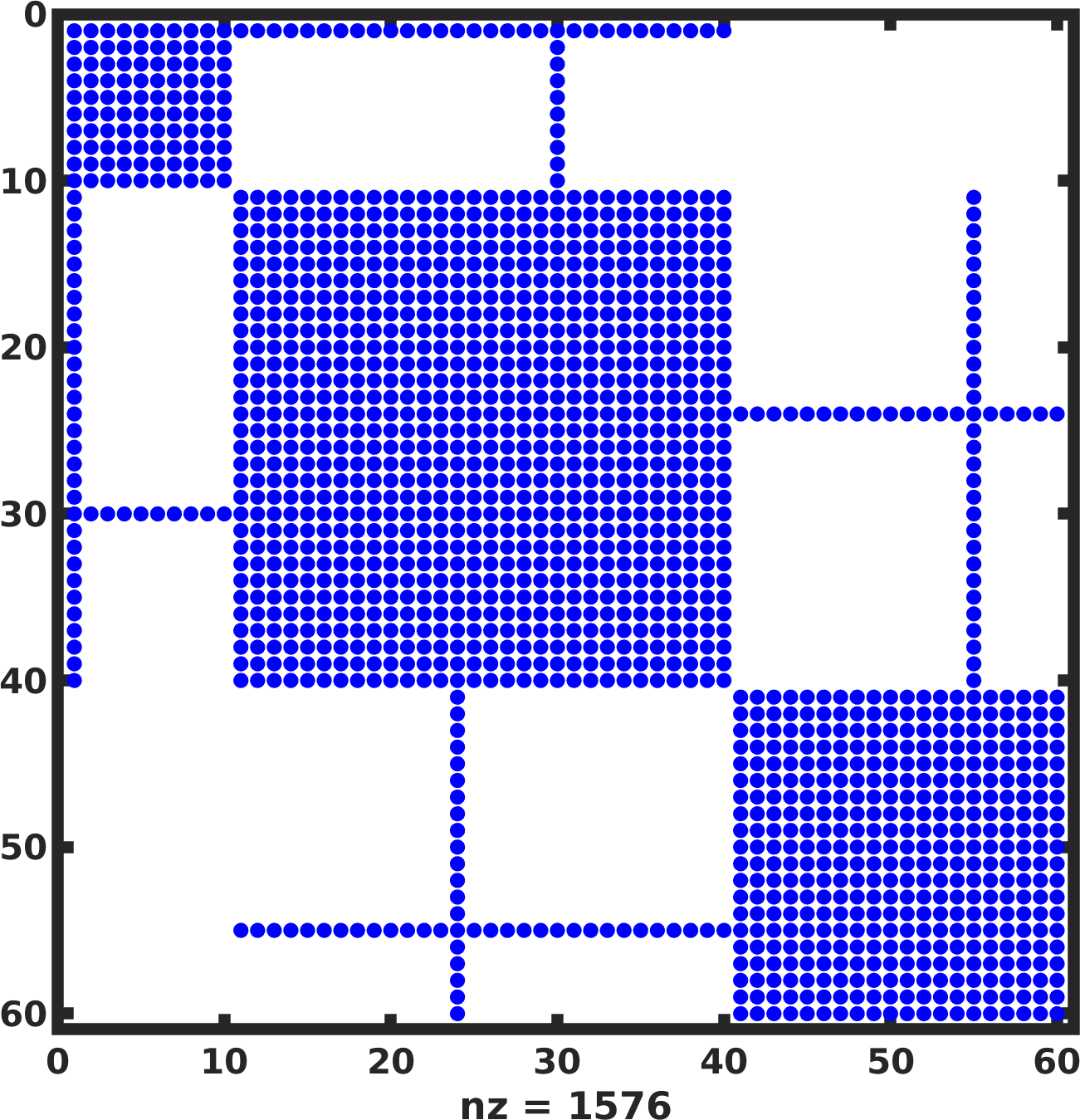} \\
        \end{tabular}
   \end{center}
\caption{\label{fig:tmatrices} Reducible and irreducible transition matrices $Q$ (left) and $Q(\epsilon)$ (right) of a $60$-state Markov chain (example \ref{exsect3}) exhibiting three invariant or three almost-invariant patterns, respectively.}
\end{figure}  
The corresponding eigenvalues of both matrices are shown in figure \ref{fig:upeigs}. As expected, the unperturbed matrix has an eigenvalue $1$ of multiplicity $3$ (figure \ref{fig:upeigs} (left)), while the perturbed matrix has two eigenvalues near the Perron root (figure \ref{fig:upeigs} (right)).
 \begin{figure}[!htb]
  \begin{center}
\begin{tabular}{cc}
    \includegraphics[width=0.45\textwidth]{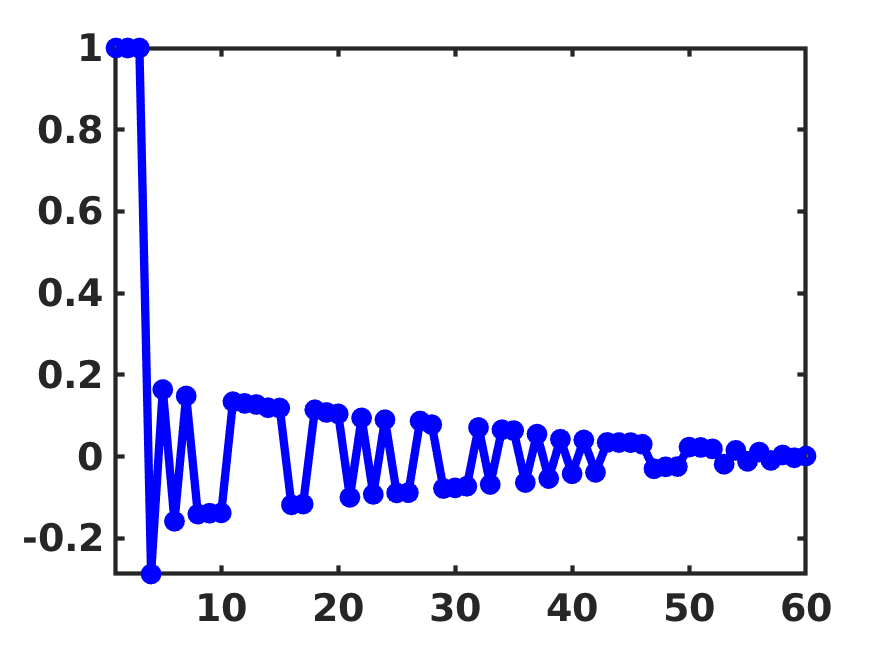} &
    \includegraphics[width=0.45\textwidth]{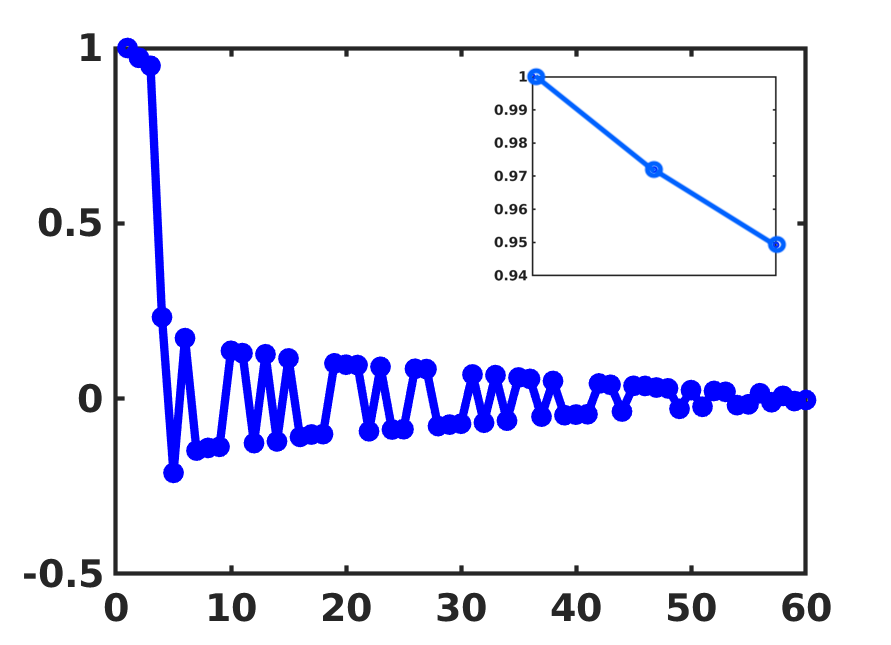}\\
  \end{tabular}
  \end{center}
\caption{\label{fig:upeigs} Eigenvalues (ordered by magnitude) of the unperturbed matrix $Q$ (left) and the perturbed matrix $Q(\epsilon)$ (right) for the $60$-state Markov chain model in example \ref{exsect3}. The perturbation results in two eigenvalues very close to one (right, see also inlet) which originate from three-fold eigenvalue $1$ (left) in the unperturbed situation. }
\end{figure}  
\end{example}

Due to reducibility, the global stationary distribution of the unperturbed transition matrix $Q$ in \eqref{pmatrix} is not unique. Indeed, each vector $V_i$, where
\begin{equation}\label{eigvecbasis_equ}
 V_i = (0,\ldots,0,\pi^{(i)},0,\ldots,0),\,\,\, i=1,\ldots,k,\,\,\text{with}\,\, {\pi^{(i)}} Q_i  = {\pi^{(i)}},\
\end{equation}
is a left eigenvector of $Q$ corresponding to the $k$-fold eigenvalue $\lambda_1 = 1$ of \eqref{pmatrix}. The eigenspace $E_{\lambda_1}$ is, thus, 
spanned by $\{V_i,\,\, i=1,\ldots,k\}$. The eigenvectors $V_i$ are only supported on $A_i$ where they have a constant sign. 
However, there exists other eigenvector bases $\{U_i,\,\, i=1,\ldots,k\}$ of $E_{\lambda_1}$ given by
\begin{equation}\label{ueigvec}
 U_i = \sum_{j=1}^k \alpha_{ij} V_j,\,\,i=1,\ldots,k,\,\, \alpha_{ij}\in \mathbb{R}.
\end{equation}
Thus, depending on the choice of $\alpha_{ij}$, each $U_i$ may partition the ``clouds'' $\{ A_i\}_{i=1}^k$ into configurations via its sign structure.  

\begin{figure}[!htb]
 \begin{center}
\begin{tabular}{ccc}    
\includegraphics[width=0.3\textwidth]{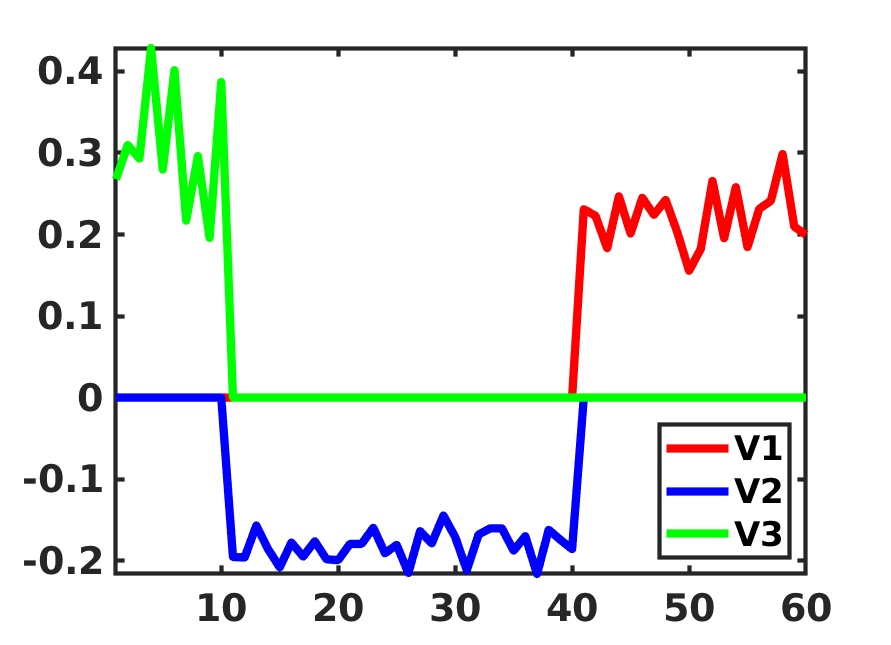} &
    \includegraphics[width=0.3\textwidth]{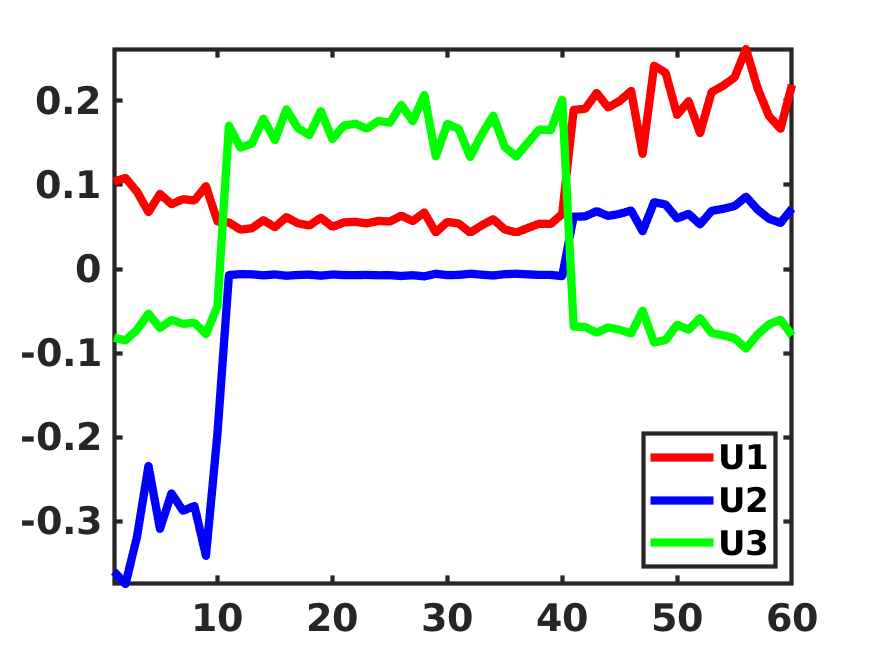} &
    \includegraphics[width=0.3\textwidth]{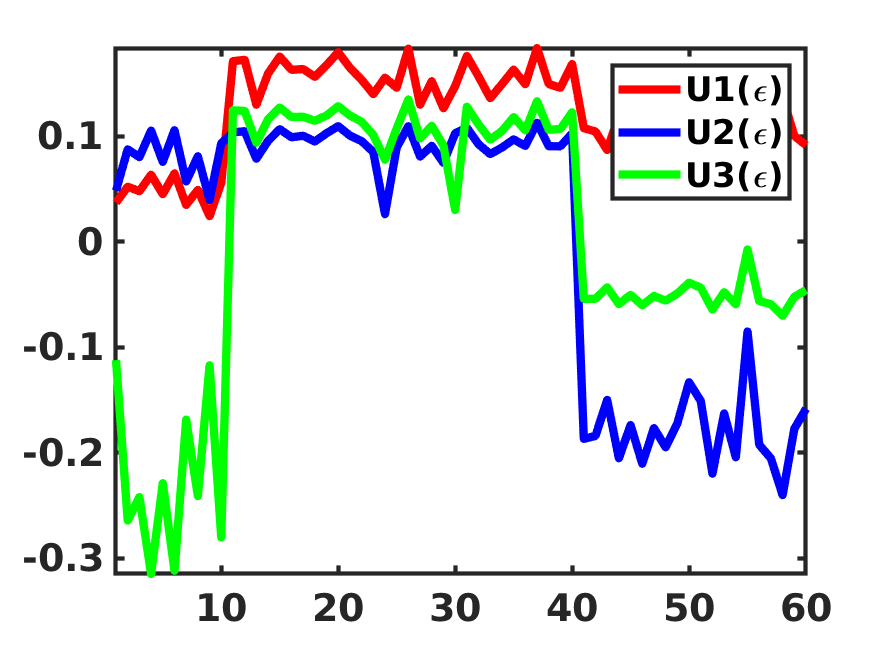}\\
       \end{tabular}
   \end{center}
   \caption{\label{fig:upeigvec} Eigenvectors $V_i$ (left) and $U_i$ (middle), $i=1,2,3$, to the $3$-fold eigenvalue $1$ of the unperturbed matrix $Q$ of the $60$-state Markov chain (example \ref{exsect3}), and eigenvectors to leading eigenvalues $U_i(\epsilon)$, $i=1,2,3$, for the perturbed matrix $Q(\epsilon)$ (right). }
\end{figure}  

\setcounter{theorem}{0}
\begin{example}[cont.]
The eigenvectors $V_i$, $i=1,2,3$ corresponding to the three-fold eigenvalue $1$ of the unperturbed matrix are only supported on the respective invariant patterns  (figure \ref{fig:upeigvec} (left)). 
In figure \ref{fig:upeigvec} (middle), the eigenvector $U_1$ is constructed to have a constant positive sign on the whole state space $\mathcal{S}$; it can be interpreted as  
a stationary distribution. $U_2$ yields a $2$-partition of $\mathcal{S}$ by lumping together $A_1$ and $A_2$. Finally $U_3$ yields a $3$-partition of $\mathcal{S}$, 
which corresponds exactly to the three invariant patterns that exist in the state space. 
In figure  \ref{fig:upeigvec} (right) the leading eigenvectors for the perturbed matrix $Q(\epsilon)$ are shown. From their sign structures a $3$-partition of $\mathcal{S}$ into almost-invariant patterns is obtained.
\end{example}

In the presence of perturbations, an explicit formula of the $k$ analytic eigenvectors corresponding to the dominant eigenvalues \,\,\textendash\,\,the eigenvalues 
clustered near $1$\,\,\textendash\,\,can be found as
\begin{equation}\label{peigvec}
\begin{split}
U_1(\epsilon) &= \pi(\epsilon) = [\pi_1(\epsilon),\pi_2(\epsilon),\ldots,\pi_N(\epsilon)],\,\,\pi_i(\epsilon)>0,\\
U_i(\epsilon) &= \sum_{j=1}^k (\alpha_{ij}+ \epsilon\beta_{ij}) V_j\\
              & + \epsilon\sum_{j=k+1}^N \frac{1}{1-\lambda_j(\epsilon)}\langle U_j, Q^{(1)}U_i\rangle_{\pi(\epsilon)} \  + O(\epsilon^2),\,\,i=2,\ldots,k,\,\,
               \alpha_{ij},\beta_{ij} \in \mathbb{R}.
\end{split}
\end{equation}
Formula \eqref{peigvec} was stated and proven in \cite{Deuflard} for the right eigenvectors of $Q(\epsilon)$. 
The proof is mainly based on (\cite{Kato}, Chp.\ $2$) but with a particular focus on reversible stochastic matrices. 
Here, we only use left eigenvectors of $Q(\epsilon)$ since left and right dominant eigenvectors are both analytic for $\epsilon\in \mathbb{R}$ and are related by 
$U_i(\epsilon) = D_N X_i(\epsilon)$, where $D_N= \text{diag}([\pi_1,\pi_2,\ldots,\pi_N])$ and $\{X_i(\epsilon),\,\, i=1,\ldots,N\}$ are the $\pi(\epsilon)-$orthogonal right 
eigenvectors of $Q(\epsilon)$. Note that $\text{sign} (U_i(\epsilon))  = \text{sign}(X_i(\epsilon))$. 

The first term in the second equation in \eqref{peigvec} suggests that the $U_i(\epsilon)'$s are actually $\epsilon$-up-or-down-shifts of the basis $V_j$ in 
equation \eqref{eigvecbasis_equ}, which were each supported on invariant patterns $A_j$. Thus, this shifting does not affect the sign structure of 
the unperturbed eigenvectors; see equation \eqref{ueigvec}. However, the second term depends on the spectral gap $\frac{1}{1-\lambda_j(\epsilon)}$ between 
the Perron root $1$ and the $N-k$ small magnitude eigenvalues of $Q(\epsilon)$. Therefore, this second term may have an influence on the sign structure of the unperturbed 
eigenvector, but only when a relatively small $\epsilon$ is chosen \cite{Deuflard}. 

With this setting of the Markov chain and the
lumped almost-invariant states, the sign structure of each dominant eigenvector in equation \eqref{eigvecbasis_equ} yields a partition of the state space. Indeed, 
each $U_i(\epsilon)$ defines a partition into $i$ nearly disjoint aggregates for $i=2,\ldots,k$, via its sign structure. In addition, $k$ dominant eigenvalues are a consequence of the occurrence of $k$ almost-invariant patterns given by the supports of $U_k(\epsilon)$.  
Finally, note that the remaining $N-k$ eigenvectors, corresponding to the spectrum $\{\lambda_j(\epsilon),j=k+1,\ldots,N \}$ bounded away from $1$, 
cannot be interpreted as \eqref{peigvec}. Indeed, the supports of these eigenvectors do not correspond to invariant patterns. However, they may play an important role 
when it comes to studying the changes of the dominant almost-invariant patterns with respect to an external bifurcation parameter.  

\section{Incompressible $2$D flows and almost-invariant sets}\label{expModels}
Since our study is motivated by geophysical applications including the splitting pattern of the Antarctic polar vortex in September 2002, we focus on models 
exhibiting vortices in their incompressible dynamics. As a first illustrative example, let us consider the following two-dimensional system of ordinary differential equations:
\begin{equation}\label{sys2}
\begin{aligned}
   \dot{x}(t)&=-\pi \sin(\pi x)\cos(\pi x) \\
    \dot{y}(t)&=\pi \cos( \pi x)\sin(\pi y)
\end{aligned}
\end{equation}
\begin{figure}[!htb]
\centering
\includegraphics[width=.4\textwidth]{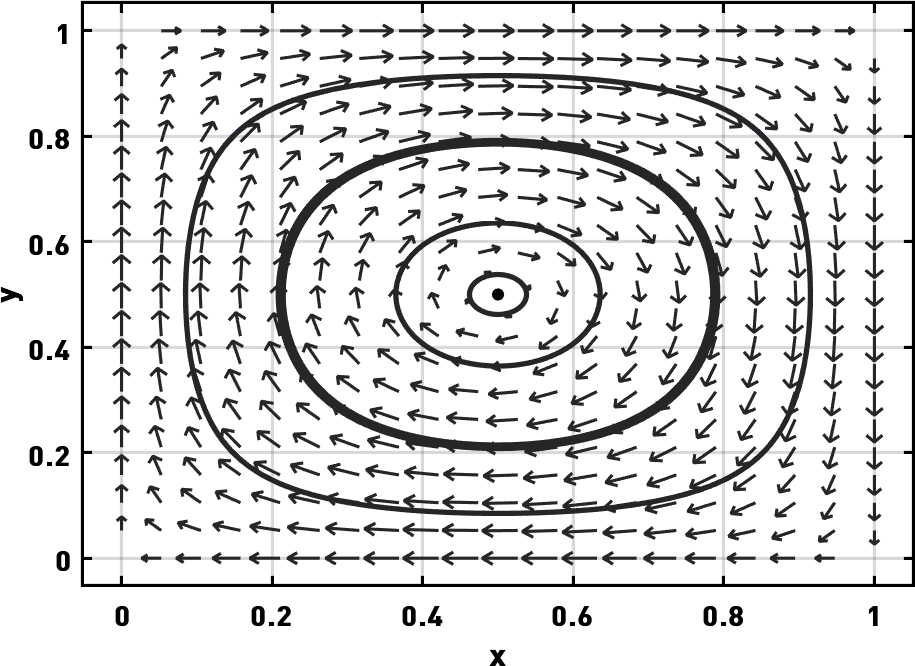}
\caption{\label{fig:DGphaseplot} Phase plane of system \eqref {sys2} consisting of periodic orbits.}
\end{figure}
From the stationary behavior of \eqref{sys2}, it is clear that every single orbit is periodic, see figure \ref{fig:DGphaseplot}. Hence, the ensemble evolution of a set of initial points under the flow 
map $S^t$ yields a bundle of closed curves for sufficiently large $t$. 
Under this rotational dynamics, one can always extract a finite number of disjoint ring-like sets $\{A_1,A_2,\ldots,A_k\}$ that partition the phase space $M$ so that the invariance 
equation $\rho_\mu(A_i) =1$ holds, for every $A_i$, $i=1, \ldots, k$. In this context, one may think of a set $A_i$ as a bundle of invariant orbits.
Note that this partition is not unique, given the particular behavior of \eqref{sys2}. We will, nevertheless, choose to work with a fixed partition of $k$ invariant 
sets. Therefore, as in section \ref{sec1}, let us suppose that the stationary dynamics within the discretized phase space yields a reducible diagonal block transition matrix $P_N$ with $k$ blocks. 
That is, the reversibilized transformation $Q_N$ in \eqref{introReversible} has the form \eqref{pmatrix}. In other words, the $k$ block matrices of $Q_N$ consist 
of clustered states such that each lumped state yields an invariant set $A_i,\,\,i=1,2,\ldots,k$.

In what follows, we will add an external perturbation to the reducible macroscopic dynamics so that the invariant sets persist but become almost-invariant sets 
 $\{A_1(\epsilon),A_2(\epsilon),\ldots,A_k(\epsilon)\}\subset \mathcal{C}_N$. That is, $\rho_\mu(A_i(\epsilon)) \approx1,\, i=1,2,\ldots,k$ as defined in equation \eqref{invarianceApprox}.
 In \cite{GFnonauto2, GFKPG} an explicit model of the perturbation was introduced and analytically formulated. It consists of "shaking" 
every box $B_i$ before and after applying the flow map $S^t$. As a consequence, only those invariant sets that resist perturbations will continue to exist as robust almost-invariant sets and are, thus, relevant in real world settings. 
Under the perturbed dynamics, the transition matrix is given by
\begin{equation}\label{matautodiff}
 (P^{t}_{N}(\epsilon))_{ij} =\frac{m(B_\epsilon(S^t(B_\epsilon(B_i))) \cap B_j)}{m(B_\epsilon(S^t(B_\epsilon(B_i)))}.
 \end{equation}
  $B_\epsilon$ is the ball centered at zero with radius $\epsilon$, which can be thought of as the perturbation amplitude. 
 $P^{t}_{N}(\epsilon)$ is actually the finite rank approximation of the explicitly diffused Perron-Frobenius operator; see \cite{GFKPG} for 
 more details 
 and the numerical implementation.
 
As in section \ref{pertubedMotiv}, the added perturbation yields a reversible row stochastic transition matrix 
 $Q^{t}_{N}(\epsilon)$ from  ${P}^{t}_{N}(\epsilon)$ analogously to \eqref{introReversible}, where $\pi(\epsilon)$ denotes the unique stationary density of $P^{t}_{N}(\epsilon)$.
Hence, $Q^{t}_{N}(\epsilon)$ has $k$ eigenvalues $\{\lambda^t_i(\epsilon)\}_{i=1}^k$ that satisfy the properties \textbf{(a)},  \textbf{(b)} and  \textbf{(c)} 
outlined in section \ref{pertubedMotiv}. 
The corresponding eigenvectors, denoted as $\{U^t_i(\epsilon)\}_{i=1}^k$, can be expressed as in equation \eqref{peigvec}.

Let $\{X^t_i(\epsilon)\}_{i=1}^k$ be the right eigenvectors of $Q^{t}_{N}(\epsilon)$ corresponding to the eigenvalues $\{\lambda^t_i(\epsilon)\}_{i=1}^k$ . 
Then due to the self-adjoint property of $Q^{t}_{N}(\epsilon)$ with respect to the inner product $\langle\cdot,\cdot \rangle_{\pi^t(\epsilon)}$, we have for $j=2,\ldots,k$
\begin{equation}\label{deflationTechnik}
\lambda_j^t(\epsilon)=\max_{ x\neq 0,\,x\in \mathbb{R}^N} \left\{ \frac{\langle Q^{t}_{N}(\epsilon) x,x\rangle_{\pi^t(\epsilon)}}{ \|x\|^2_{\pi^t(\epsilon)}}\right\}, 
\end{equation}
under the $\pi^t(\epsilon)$-orthogonal constraint $$\langle x, \textbf{1}\rangle_{\pi^t(\epsilon)} =\langle x, X^t_2(\epsilon)\rangle_{\pi^t(\epsilon)}=...
                     = \langle x, X^t_{j-1}(\epsilon)\rangle_{\pi^t(\epsilon)}=0.$$
Note that $\textbf{1}=X^t_1(\epsilon) = [1,1,\ldots,1]$ denotes the right stationary distribution of $Q^{t}_{N}(\epsilon)$.
In \cite{GFnonauto2, GFKPG}, the eigenvalue $\lambda_2^t(\epsilon)$ and the corresponding left eigenvector 
$U^t_2(\epsilon) = \mathcal{D}_N X^t_2(\epsilon), \,\, \mathcal{D}_N = \text{diag}([\pi^t(\epsilon)_1,\pi^t(\epsilon)_2,\ldots,\pi^t(\epsilon)_N])$ were used to approximate 
two robust maximal almost-invariant sets. Indeed, due to the $\pi^t(\epsilon)$-orthogonality relations among the right eigenvectors $\{X^t_i(\epsilon)\}_{i=1}^k$ and the positive sign of 
$X^t_1(\epsilon)$, the sign structure of $U^t_1(\epsilon)$ is given as 
$$ \text{sgn}(U^t_1(\epsilon)) = (+,+,+,+,\ldots,+,+,+,+,+,\ldots,+).$$
We can therefore predict the sign structure of $U^t_2(\epsilon)$ as follows
$$ \text{sgn}(U^t_2(\epsilon)) = (+,+,+,+,\ldots,+,-,-,-,-,-,\ldots,-),$$
subject to a convenient box reordering. Hence, it follows that 
$U^t_2(\epsilon)$ yields positive and negative level sets, which partition the phase space into two dominant almost-invariant sets, whenever 
$\lambda_2^t(\epsilon)\approx 1$. Similarly, since further $k-2$ eigenvalues are clustered near $1$,  each eigenvector 
$U^t_j(\epsilon)$ yields a sign structure that may be sorted so that $j$ almost-invariant sets are obtained. In \cite{Deuflard}, all $k-1$ leading eigenvectors are used to compute 
almost-invariant sets. This method does not need the corresponding eigenvalues, but only the sign structures of the eigenvectors. However, in this work we use the eigenvectors separately, because we ultimately need to study the trends of the corresponding eigenvalues to understand bifurcation of patterns. 

Given \eqref{sys2}, we can numerically compute and visualize the eigenvector patterns $U^t_j(\epsilon),\,\,j=1,2,\ldots,k$, as well as their corresponding eigenvalues 
$\lambda_j(\epsilon)$. For this we use GAIO \cite{Gaioo}, which is a MATLAB-based software package for set-oriented numerics in dynamical systems. 
\begin{figure}[!htb]
\begin{center}
\begin{tabular}{cccc}
  \includegraphics[width=0.2\textwidth]{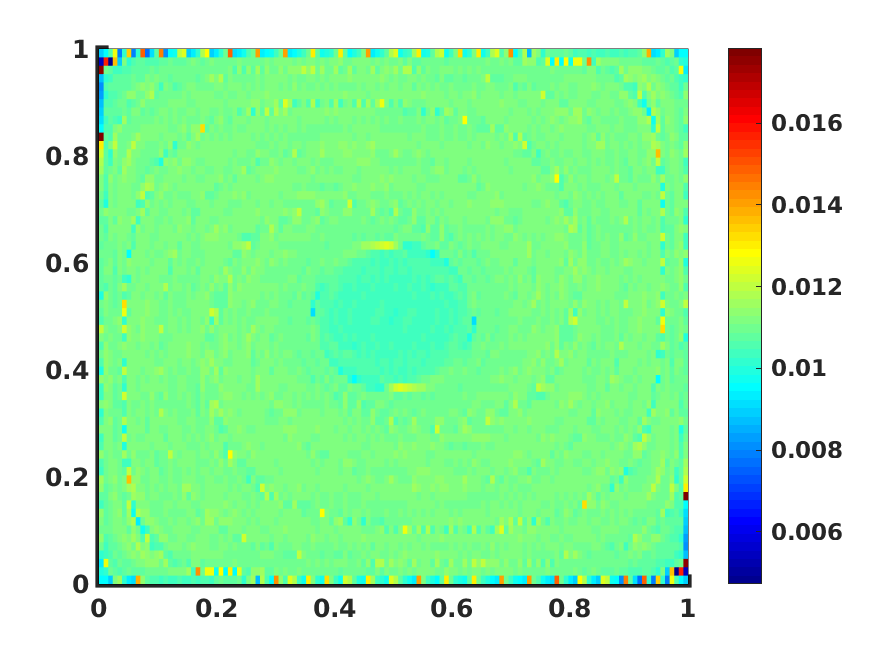} &
   \includegraphics[width=0.2\textwidth]{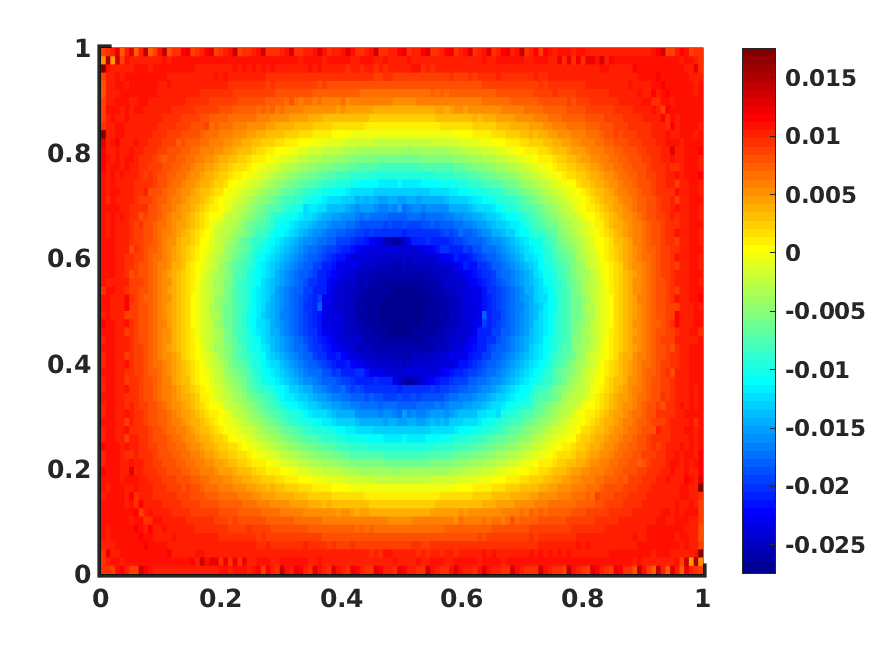} &
  \includegraphics[width=0.2\textwidth]{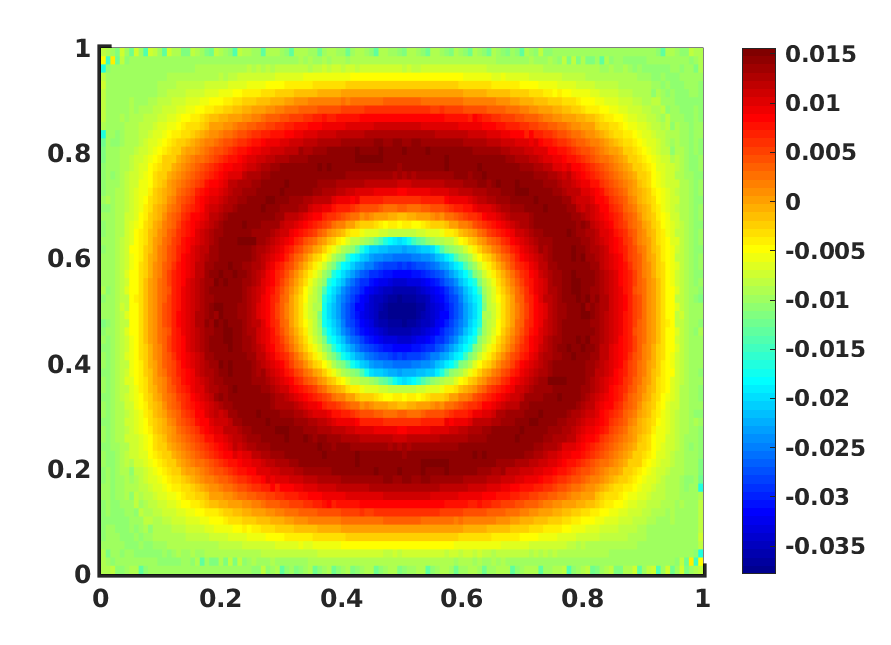} &
 \includegraphics[width=0.2\textwidth]{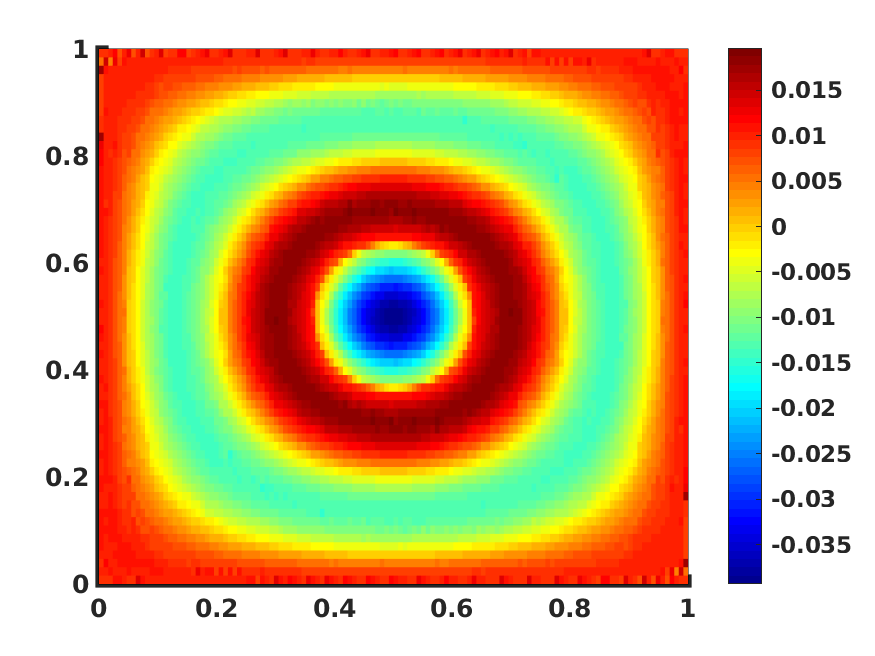} \\
 {\scriptsize $U_1^t(\epsilon)$} & {\scriptsize $U_2^t(\epsilon)$} & {\scriptsize $U_3^t(\epsilon)$} &{\scriptsize $U^t_4(\epsilon)$}\\
\end{tabular}
\end{center}
\caption{\label{fig:pEV} First 4 dominant eigenvectors of $Q^{t}_{N}(\epsilon)$ for model \eqref{sys2}.}
\end{figure}  
We approximate the flow map by using a fourth order Runge Kutta ODE solver with a time interval of length $1$ and step size of $h=0.01$, i.e $100$ time steps. The domain is subdivided into $2^{depth}$ rectangular grid sets (boxes). Here, we use $depth =13$, which gives $N=2^{13}= 8192$ boxes $B_i$ that partition the phase space $M$. In each box $900$ test points are uniformly samples as initial data for constructing 
the transition probabilities of the $N \times N$-transition matrix.

\begin{figure}[!htb]
  \begin{center}
    \includegraphics[width=0.3\textwidth]{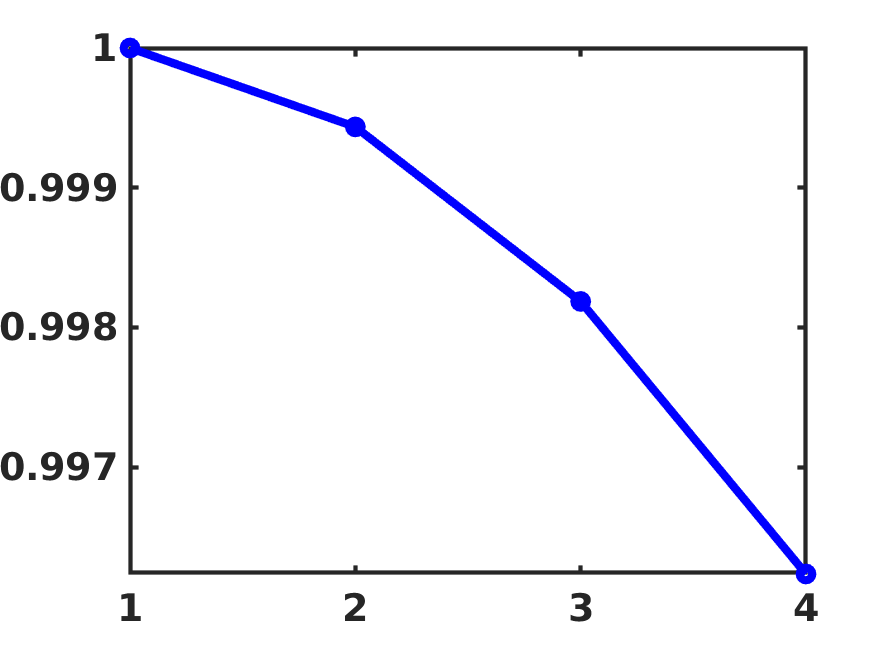}
 \end{center}
\caption{\label{fig:first4eigs} First 4 dominant eigenvalues of $Q^{t}_{N}(\epsilon)$ from system \eqref{sys2}.}
\end{figure}

In figure \ref{fig:pEV}, the $k=4$ dominant eigenvectors are plotted, with the corresponding eigenvalues shown in figure \ref{fig:first4eigs}. These are all clustered near $1$, as the result 
of additional external perturbations. Note that the numerical discretization induces a small amount of noise in the order of magnitude 
of the box diameters (\cite{numGary}, Lemma 2.2). 
That is, the numerical discretization directly yields an approximation of \eqref{matautodiff}, and, hence, it is not necessary to add explicit diffusion in practice, although it is required on the theoretical level. Also note that the leading eigenvector $U_1(\epsilon)$ is approximately constant due to area preservation of the underlying system \eqref{sys2}, with some small numerical artefacts at the boundary of the domain.

A set-oriented bifurcation analysis of a dynamical system will be 
exclusively based on studying the changes of the spectral data as a response to qualitative changes in the underlying dynamics. That is, one needs to focus on both the 
eigenvectors and their corresponding eigenvalue.
In this way, trends of the spectral data can be used to understand whether or not  
there is hint of any qualitative changes of patterns generated by the corresponding eigenvectors. 
But, beforehand, we will first consider some toy models and investigate bifurcations of patterns in an experimental manner.

\section{Numerical experiments of bifurcation}\label{numExp}
Now, we start to dive into the main purpose of this work through an experimental approach. We study the changes in the trends of the dominant spectrum when the almost-invariant patterns 
undergo different qualitative changes. This may be understood as a "bifurcation analysis" of the stationary macroscopic dynamics of the Markov chain.\\
The process resulting in qualitative changes of a pattern can only occur in two ways: Either it starts from the inside towards the outside of the pattern or the other way round. 
\begin{example}\label{firstExample}
We revisit the $60$-states Markov chain with the perturbed $3$ invariant patterns as introduced in example \ref{exsect3}. In this experiment, we want to understand how the spectrum 
behaves when the change of the pattern starts from its boundary. 
\begin{figure}[!htb]
\begin{center}
\begin{tabular}{cccc}
 \includegraphics[width=0.2\textwidth]{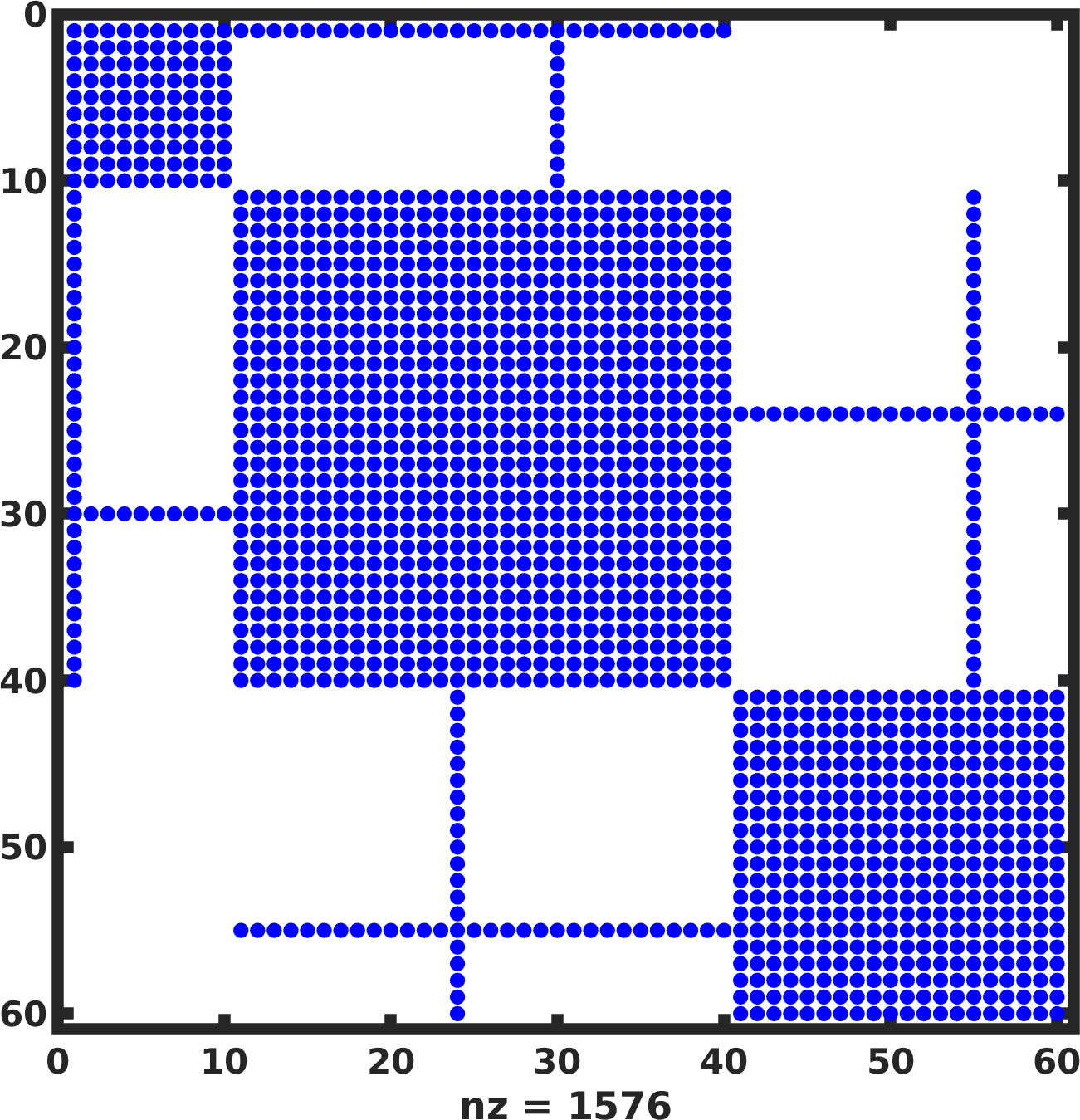}  &
  \includegraphics[width=0.2\textwidth]{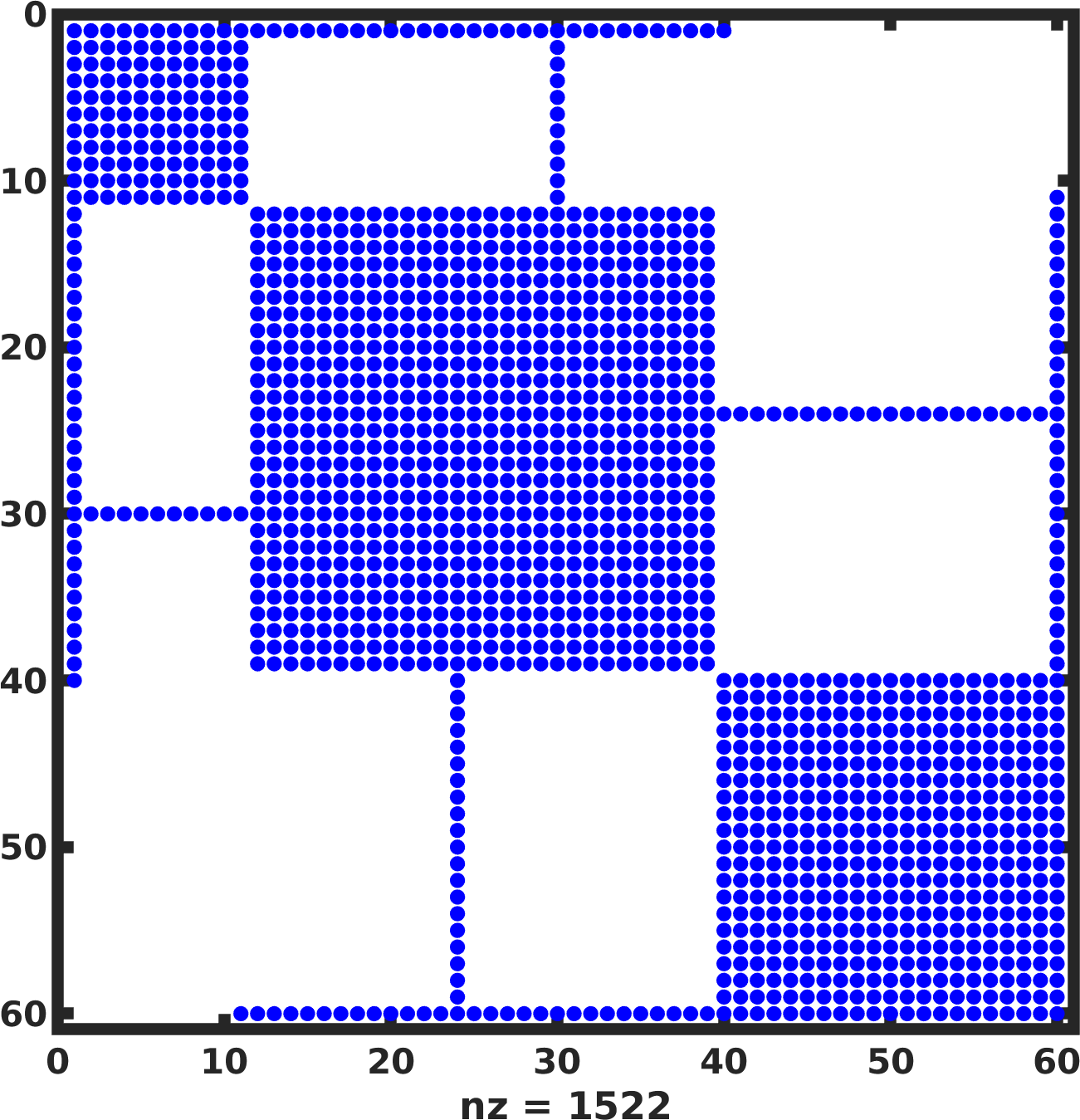} &
   \includegraphics[width=0.2\textwidth]{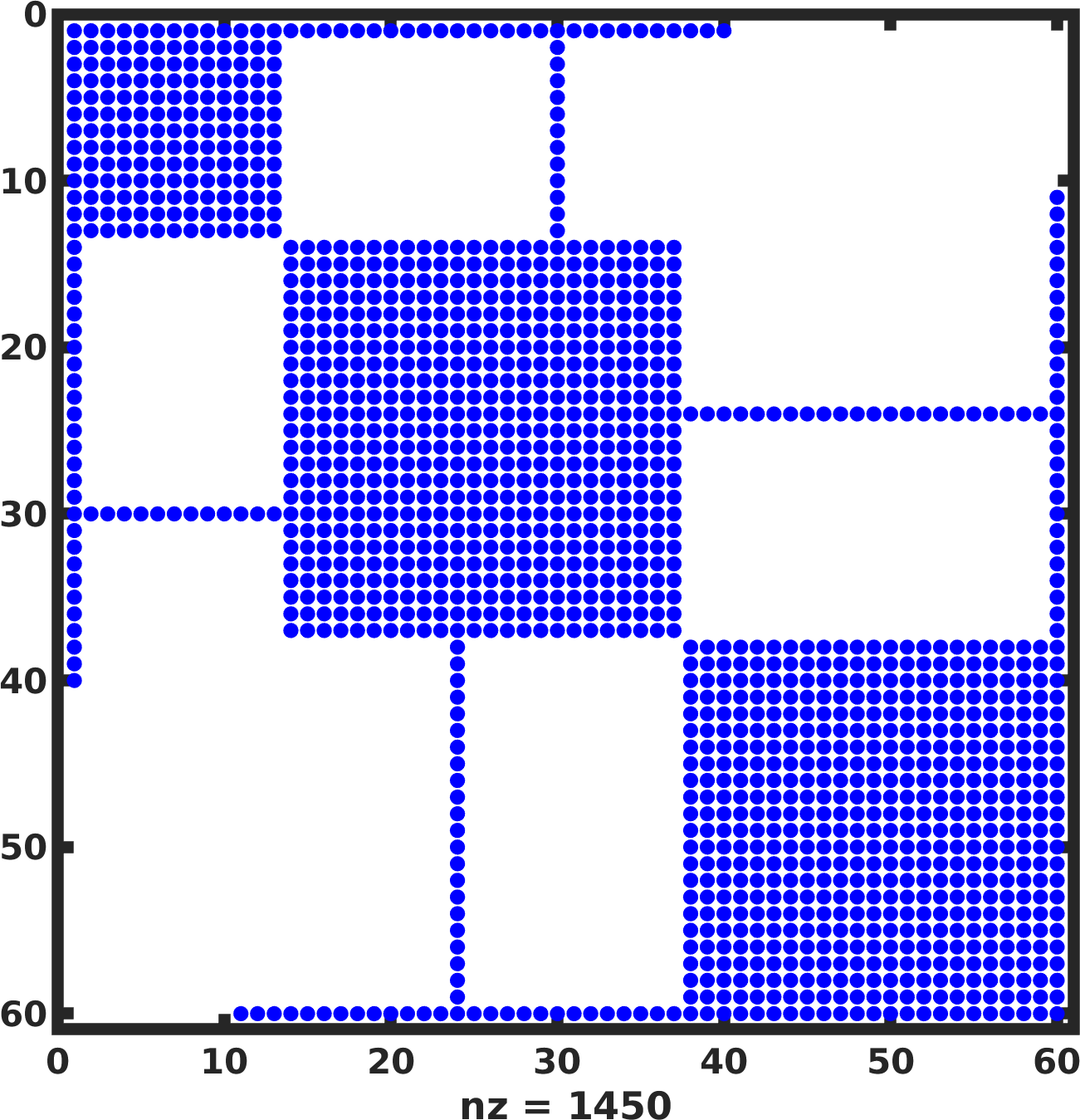} &
    \includegraphics[width=0.2\textwidth]{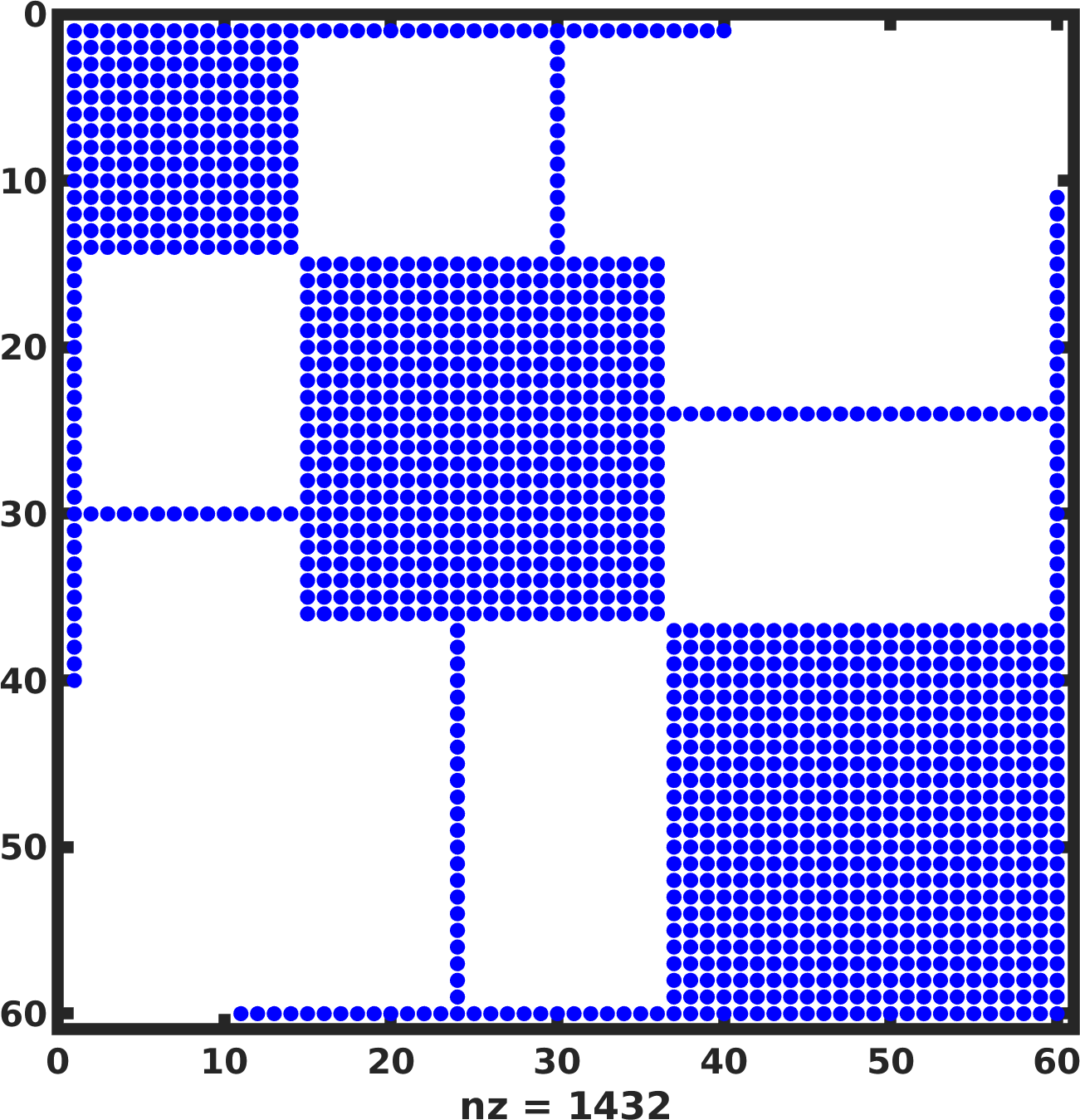} \\
     \includegraphics[width=0.2\textwidth]{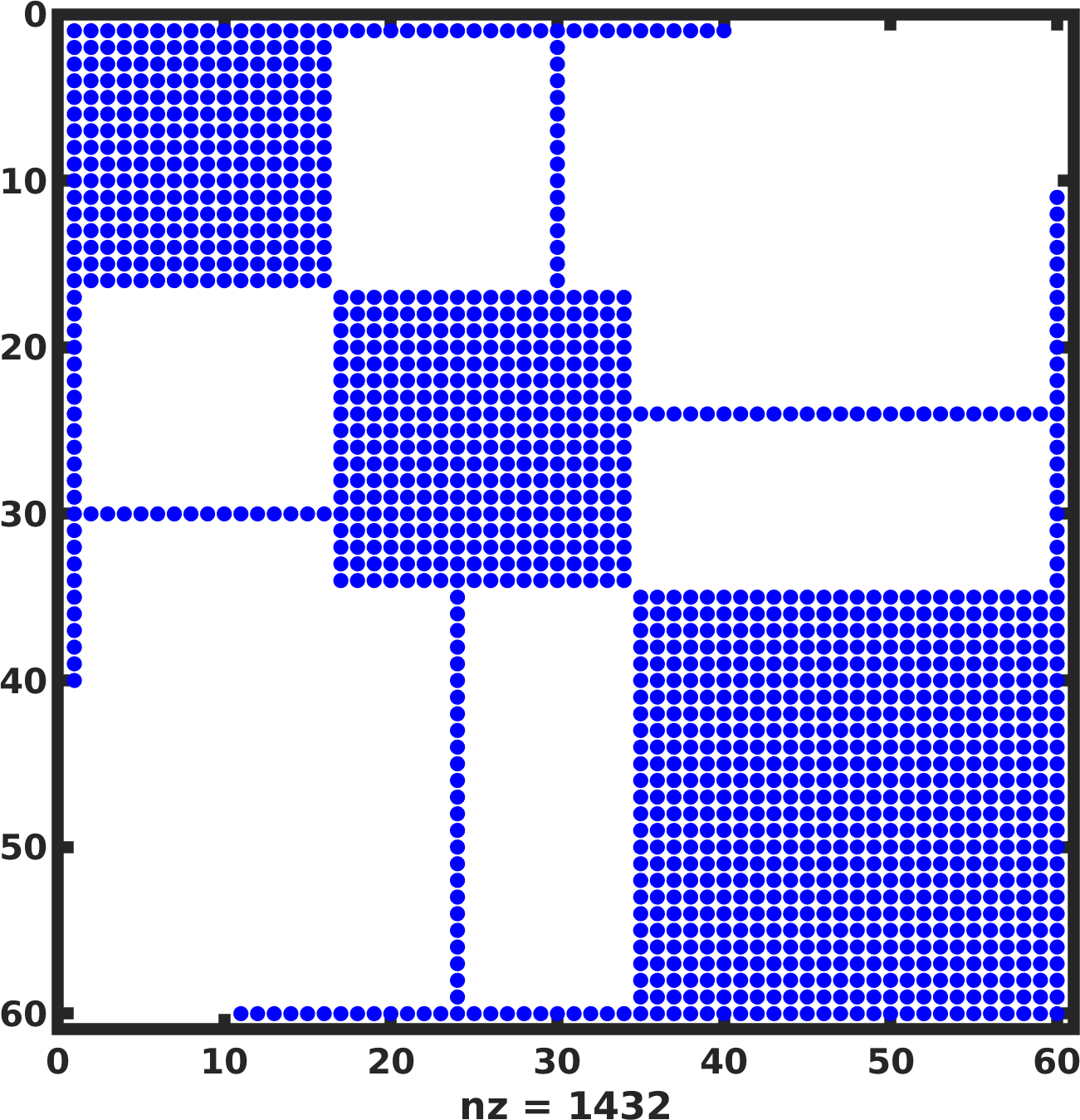}  &
  \includegraphics[width=0.2\textwidth]{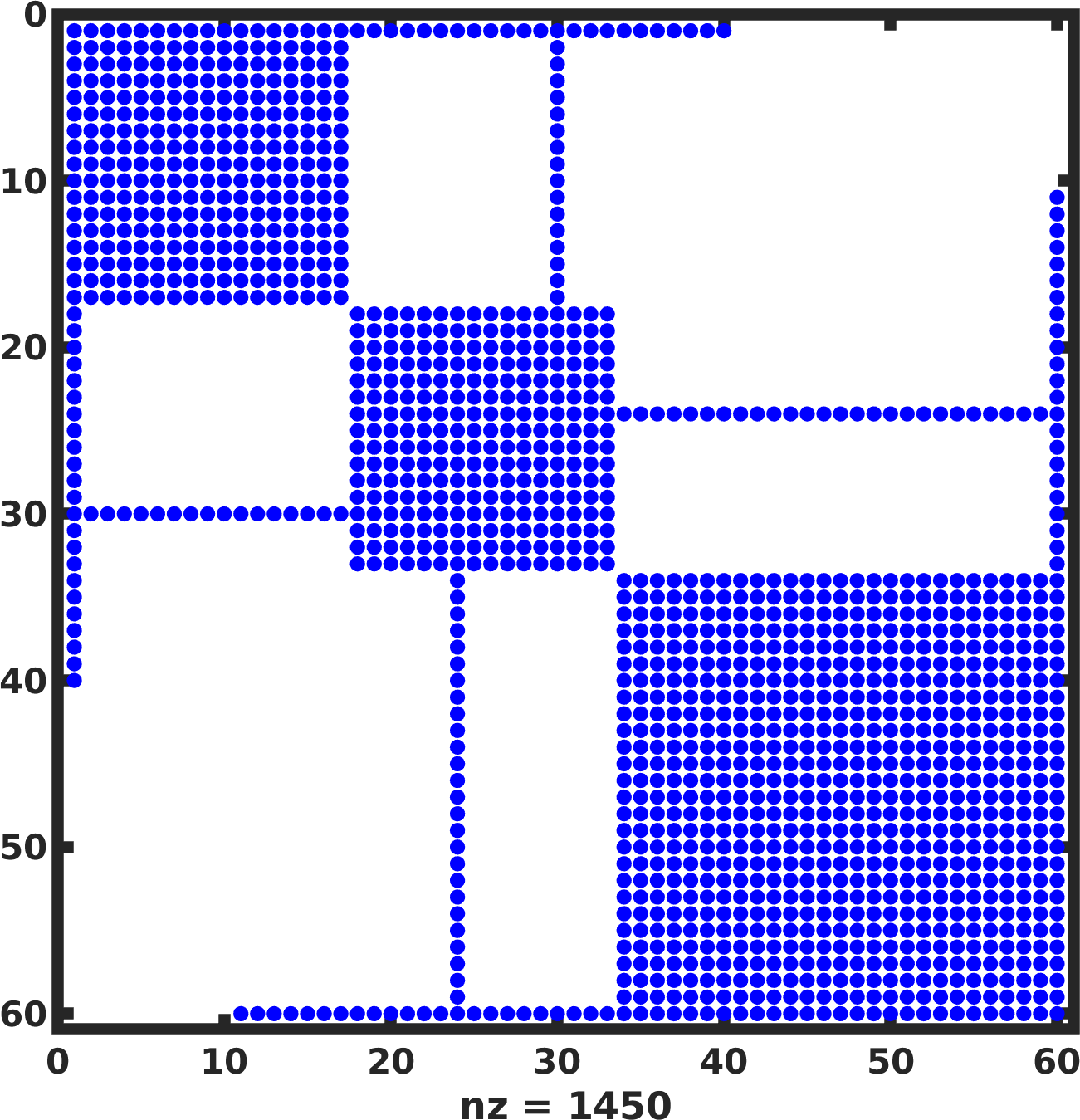} &
   \includegraphics[width=0.2\textwidth]{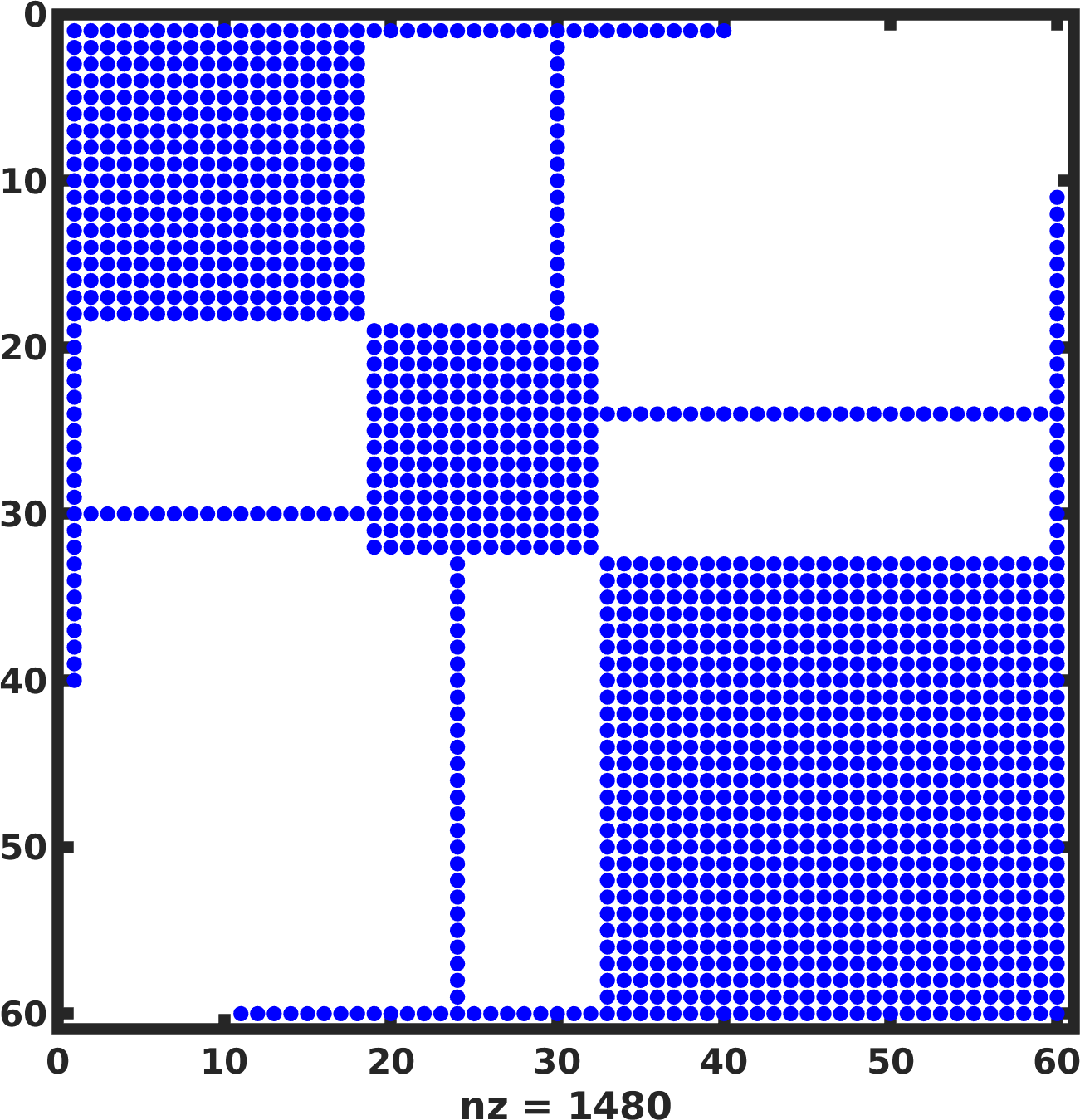} &
    \includegraphics[width=0.2\textwidth]{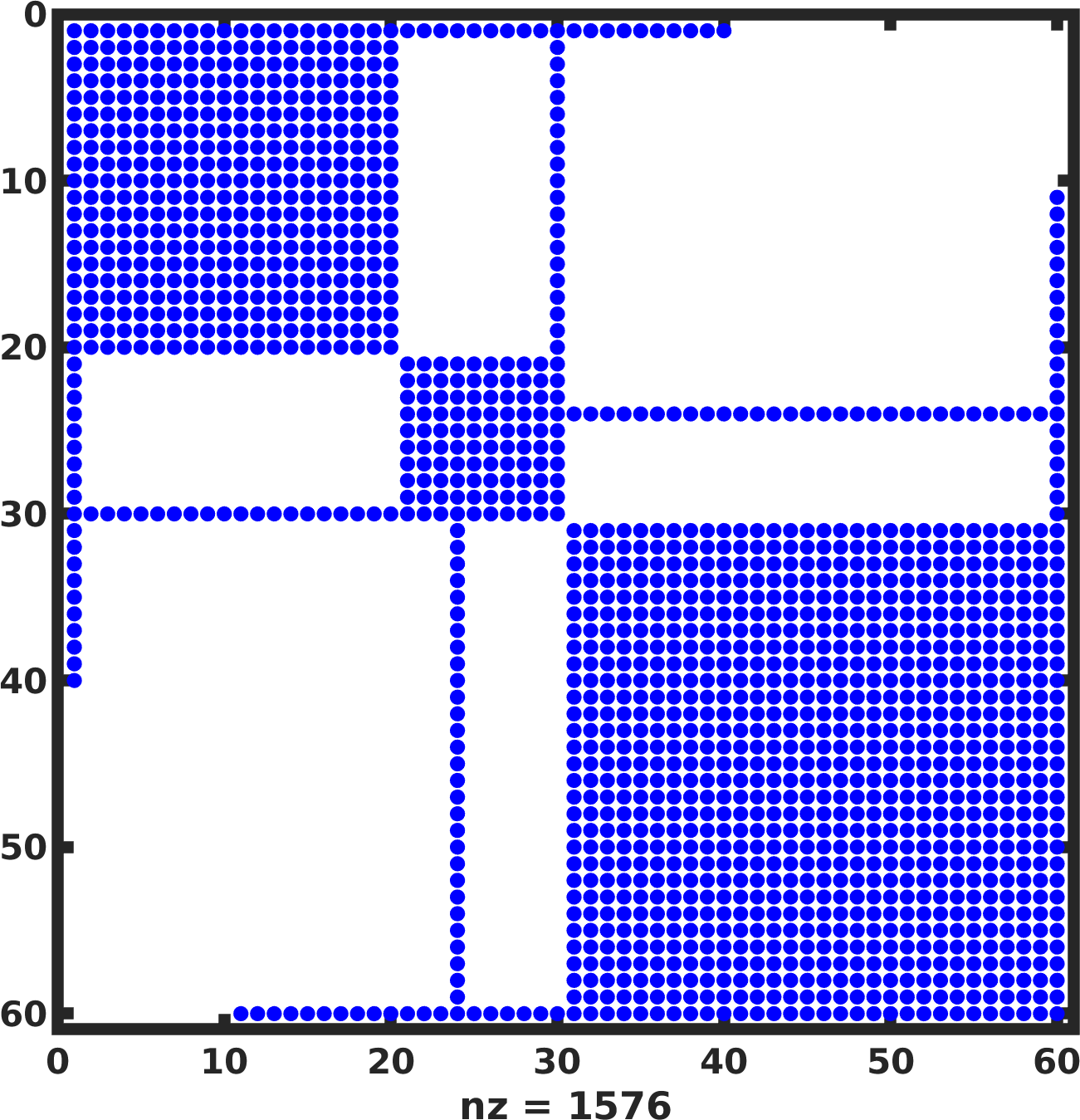} \\
\end{tabular}
\end{center}
\caption{\label{fig:threefold} Different transition matrices of example \ref{firstExample}, where the outer two almost-invariant patterns grow at the expense of the center one. }
\end{figure}  
Thus, as shown in figure \ref{fig:threefold}, we manually decrease uniformly the size of the middle invariant pattern, $A_2$, 
while increasing the size of both $A_1$ and $A_3$, simultaneously. These changes are captured by the dominant spectrum as illustrated in figure \ref{fig:threefold_eigs}. 
\begin{figure}[!htb]
\centering
\includegraphics[width=.4\textwidth]{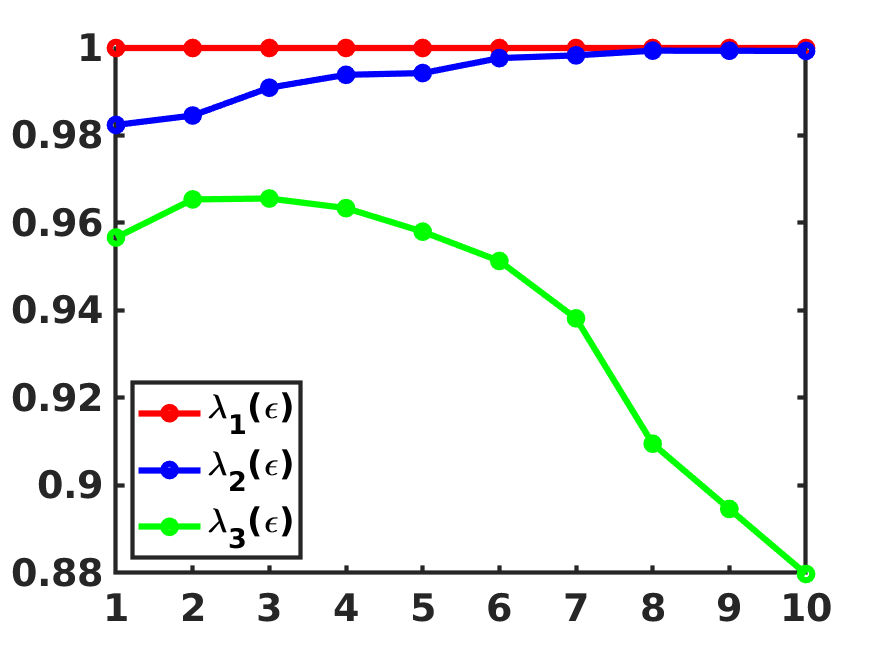}
\caption{\label{fig:threefold_eigs} Spectral signature of the shrinking of one almost-invariant pattern in example \ref{firstExample}, with two other patterns growing and  becoming more invariant.  }
\end{figure}
Eigenvalue $\lambda_3(\epsilon)$ decreases in magnitude as the middle pattern $A_3$ shrinks in size. This shrinking process is captured in the eigenvector $U_3(\epsilon)$, where
the support of $U_3(\epsilon)$ in $A_3$ is becoming smaller and smaller, as demonstrated in figure \ref{fig:threefold_evects}.
The opposite is noticed in the changing process of $U_2(\epsilon)$. The corresponding eigenvalue  $\lambda_2(\epsilon)$ approaches $1$ as $\lambda_3(\epsilon)$ decreases. 
In this process, one can clearly see that the system tends to become nearly reducible with two growing lumped states $A_1(\epsilon)$ and $A_3(\epsilon)$.  That explains the growth of  
$\lambda_2(\epsilon)$ towards $1$.
\begin{figure}[!htb]
\begin{center}
\begin{tabular}{cccc}
\includegraphics[width=0.2\textwidth]{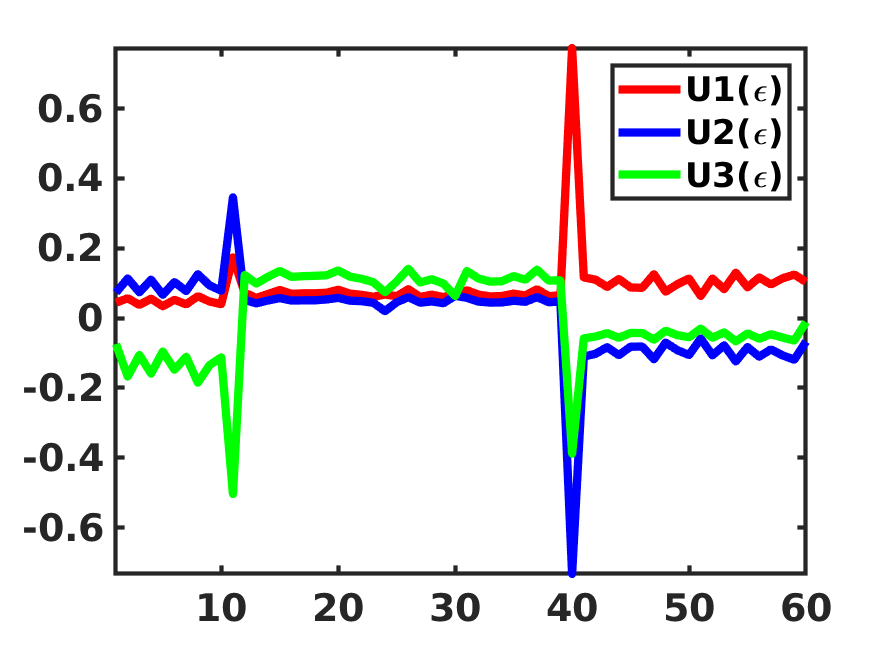}  &
  \includegraphics[width=0.2\textwidth]{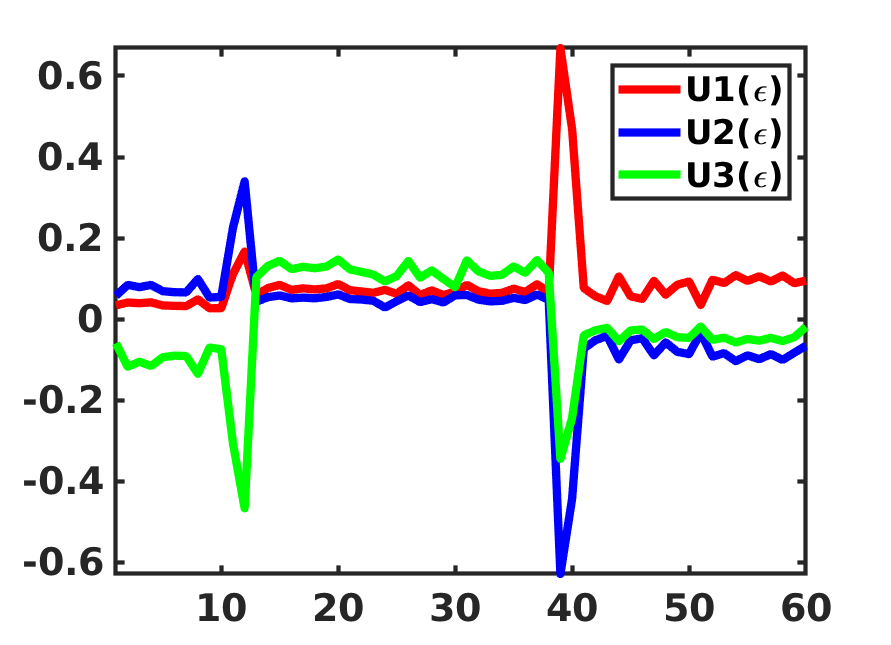} &
   \includegraphics[width=0.2\textwidth]{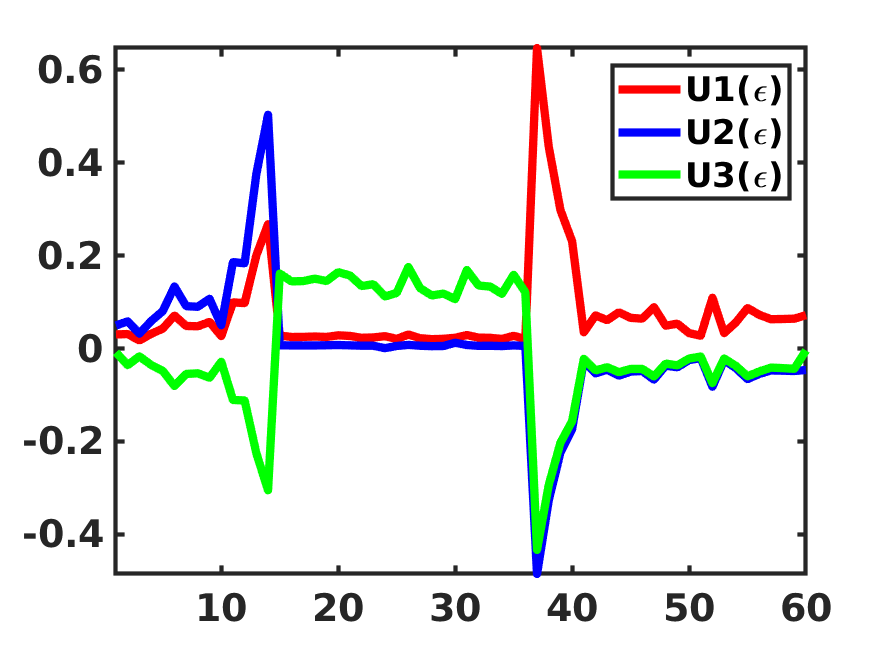} &
    \includegraphics[width=0.2\textwidth]{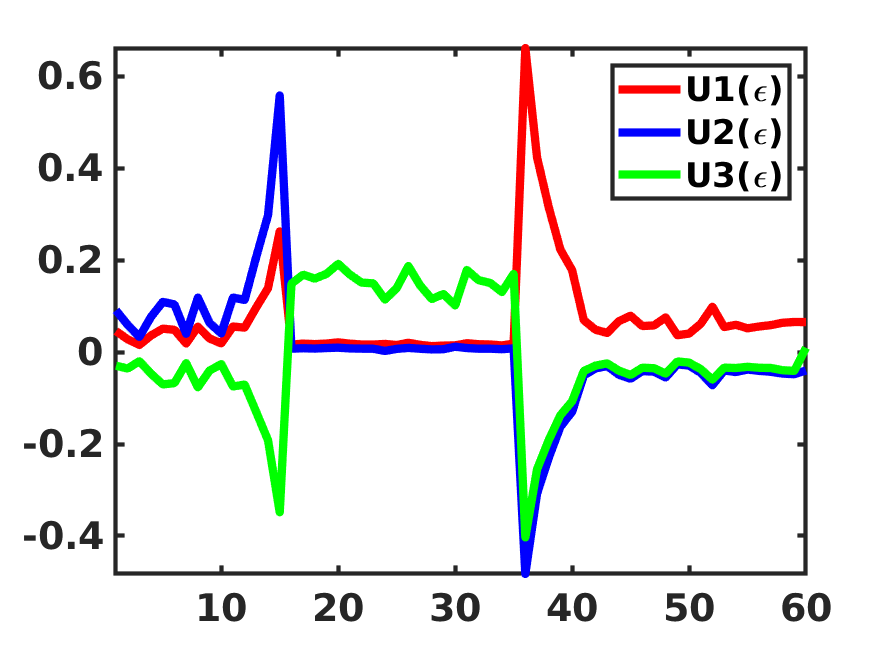} \\
     \includegraphics[width=0.2\textwidth]{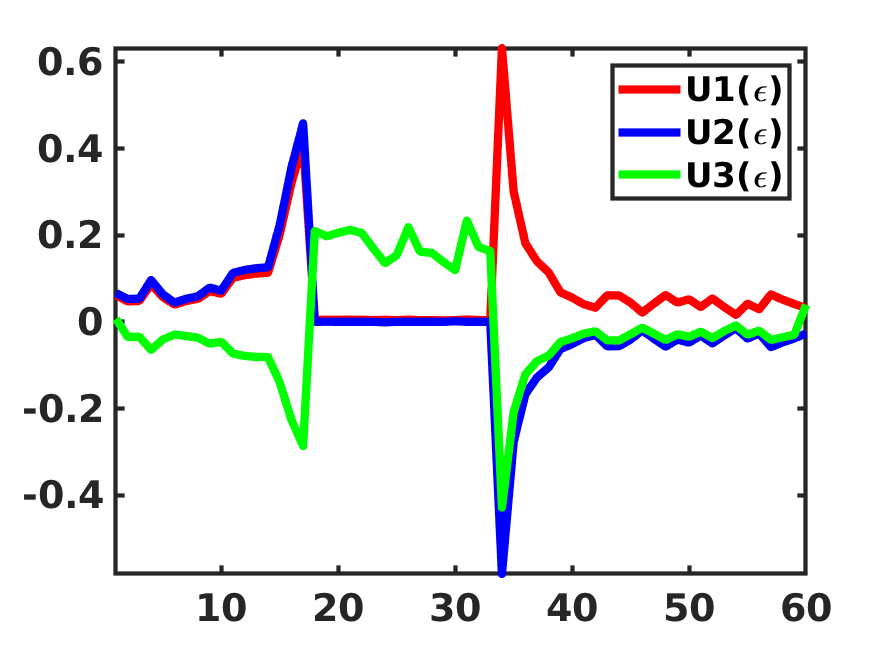}  &
  \includegraphics[width=0.2\textwidth]{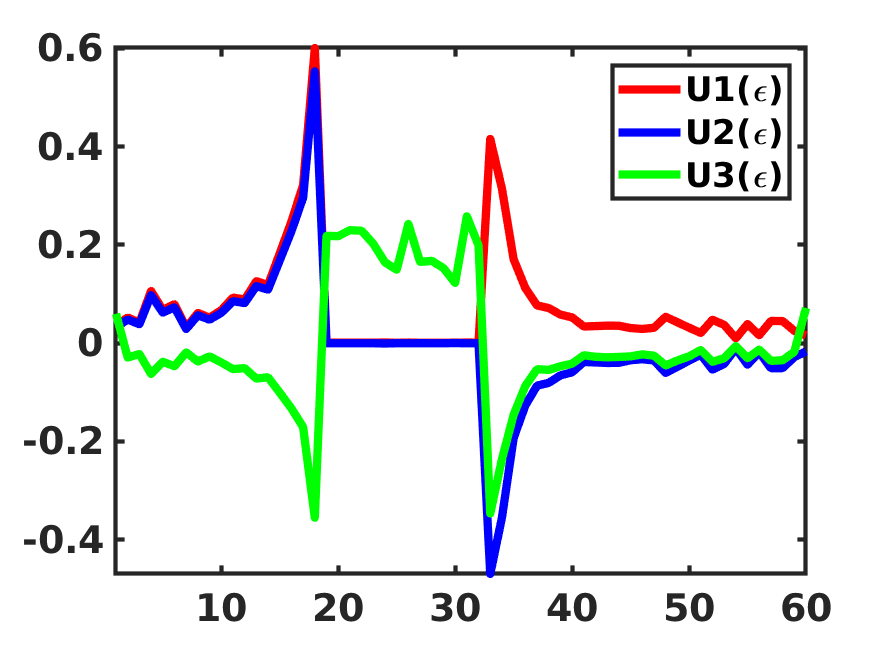} &
   \includegraphics[width=0.2\textwidth]{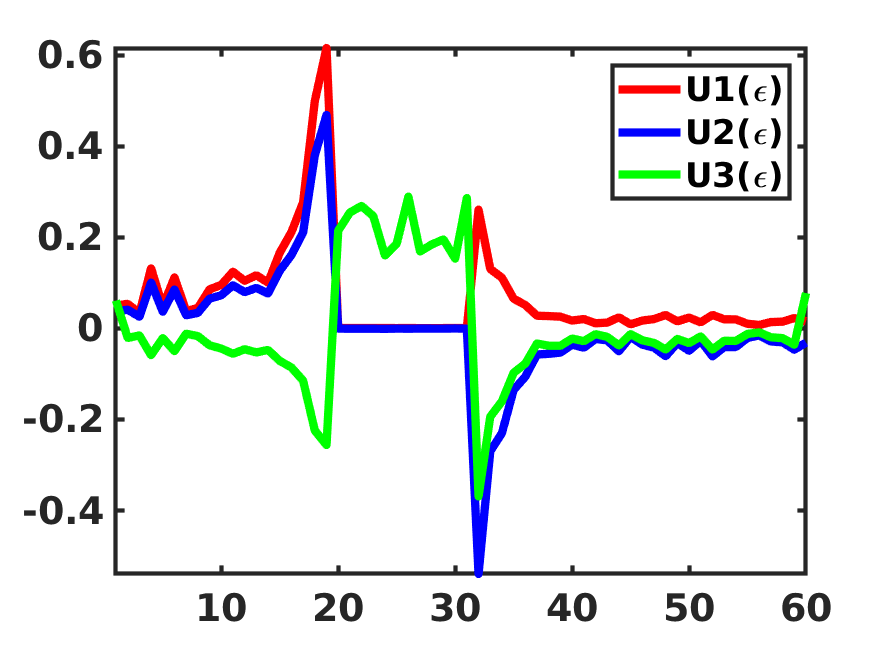} &
    \includegraphics[width=0.2\textwidth]{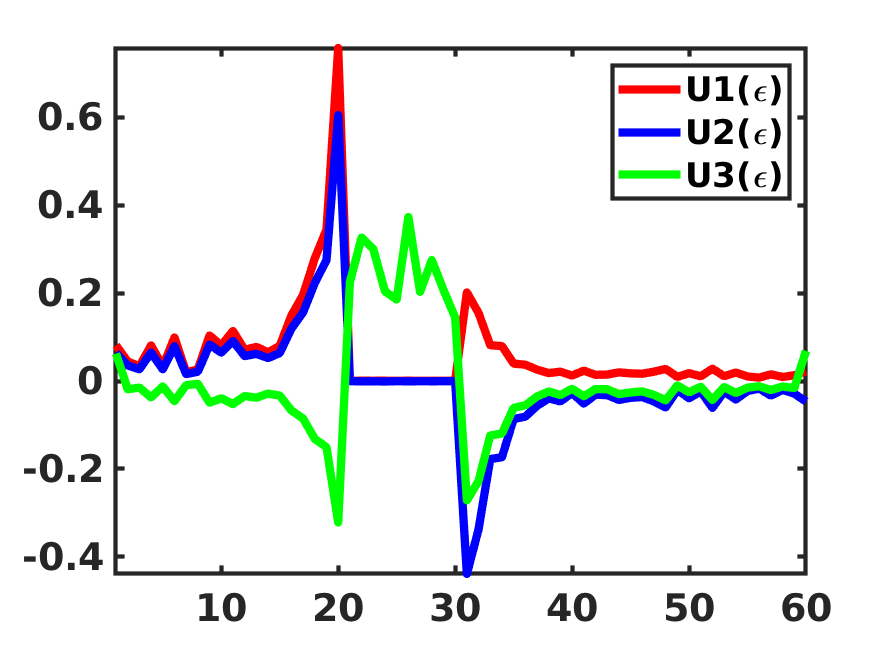} \\
    \end{tabular}
    \end{center}
\caption{\label{fig:threefold_evects} Changes in the three dominant eigenvectors for the transition matrices shown in figure \ref{fig:threefold}.}
\end{figure}  
It is necessary to understand the behavior of the eigenvalues and their correct interpretation with respect to the dynamics of the almost-invariant patterns. Indeed, this experiment 
clearly suggests a relationship between the eigenvalues and the size of the patterns. 
\end{example}

\begin{example}\label{seconExample}
Here the qualitative change is provoked from the interior of the middle pattern $A_2$. The aim is to experiment the behavior of the spectrum 
with respect to a sudden growing change from a local region. The corresponding transition matrices of the gradually changed Markov chain are shown in figure \ref{fig:picthfork}. 
The evolution of the dominant eigenvalues shown in figure \ref{fig:pitch_eigs} indicates the importance of the eigenvalue $\lambda_4(\epsilon)$, which is not part of the dominant spectrum at first. It increases very quickly in magnitude until it 
crosses $\lambda_3(\epsilon)$. The corresponding 
eigenvector, $U_4(\epsilon)$, is supported on the newly born almost-invariant pattern as illustrated in figure \ref{fig:pitch_evects}.
\begin{figure}[!htb]
\begin{center}
\begin{tabular}{cccc}
 \includegraphics[width=0.2\textwidth]{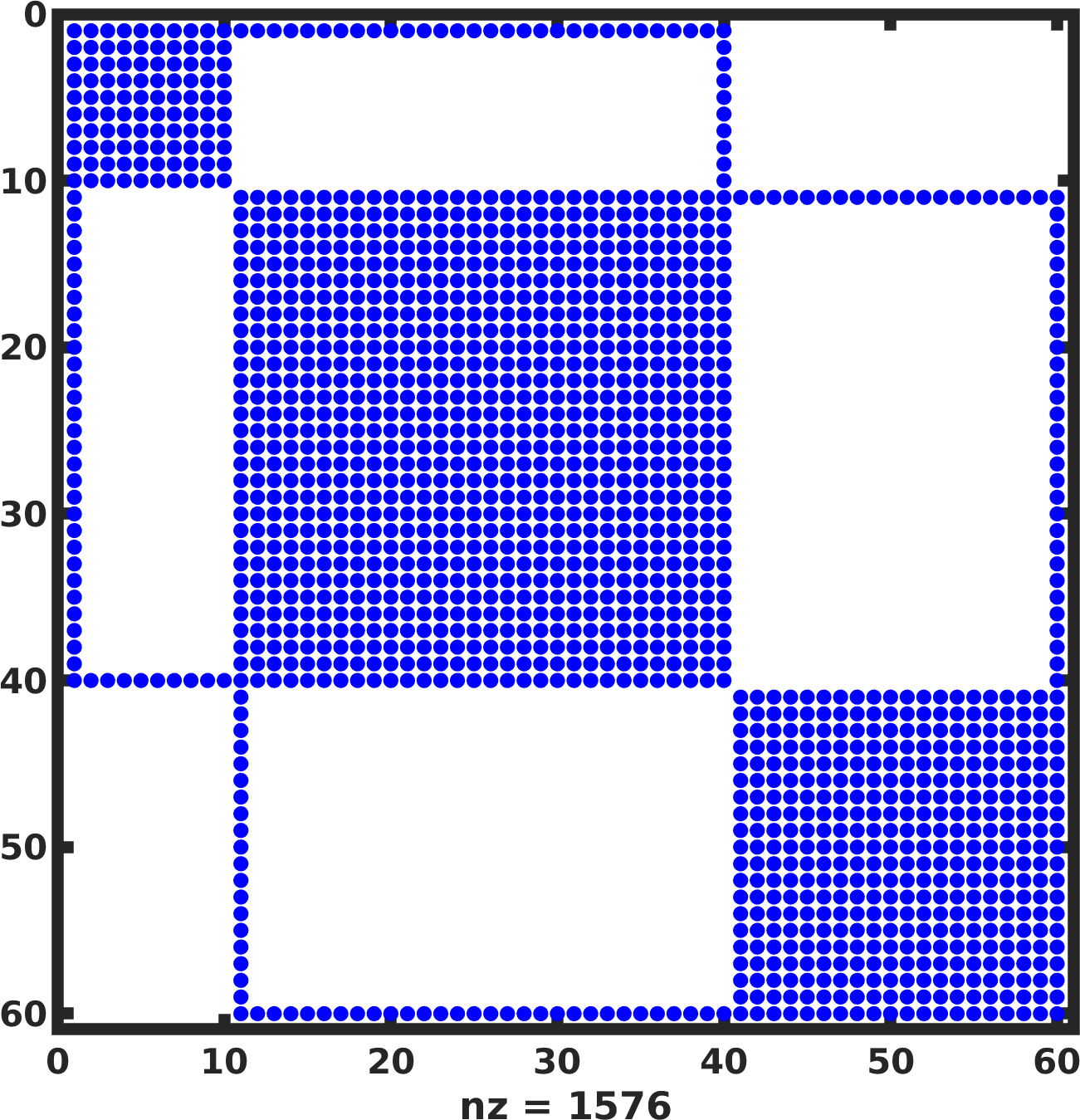}&
 \includegraphics[width=0.2\textwidth]{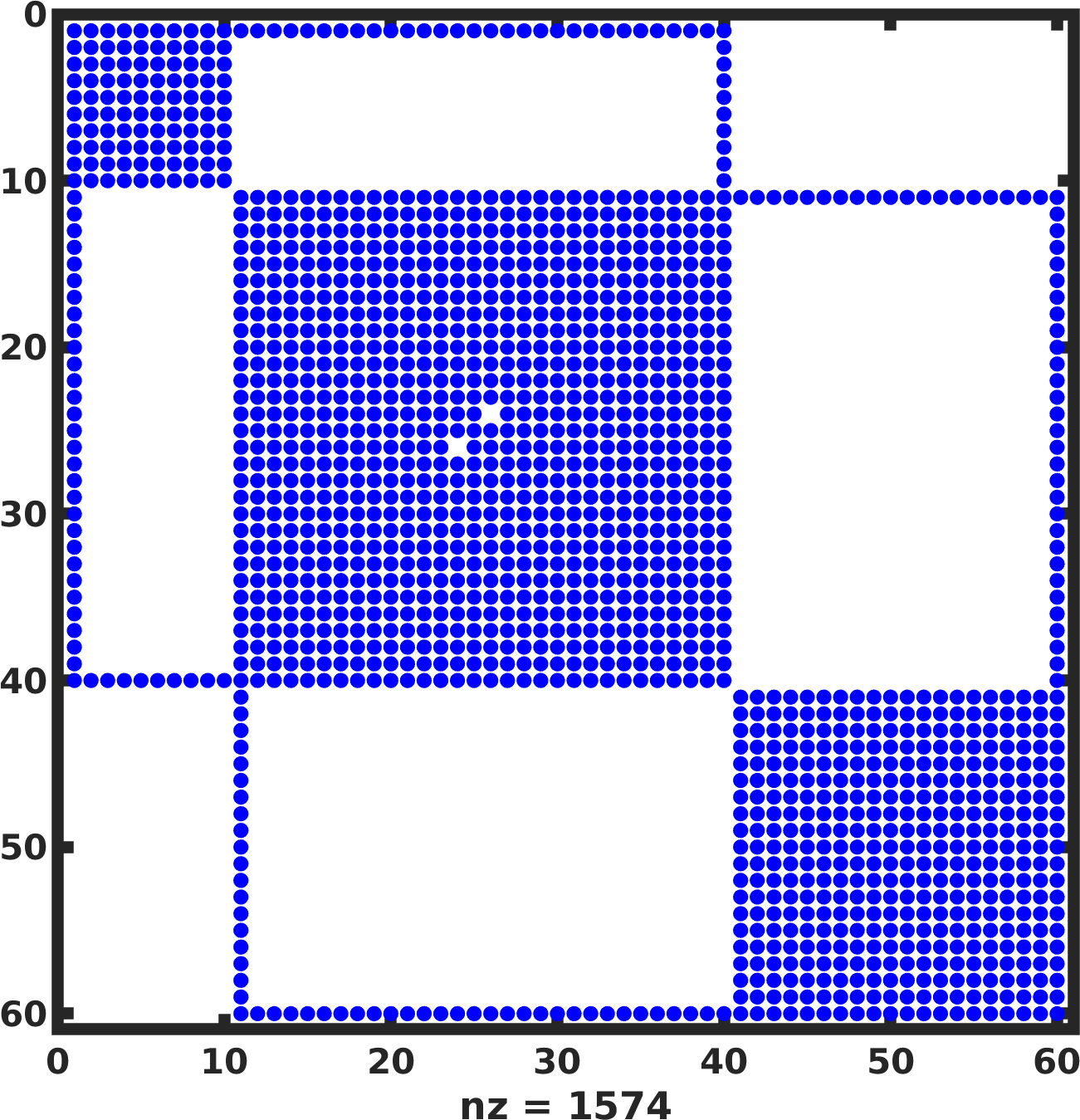} &
   \includegraphics[width=0.2\textwidth]{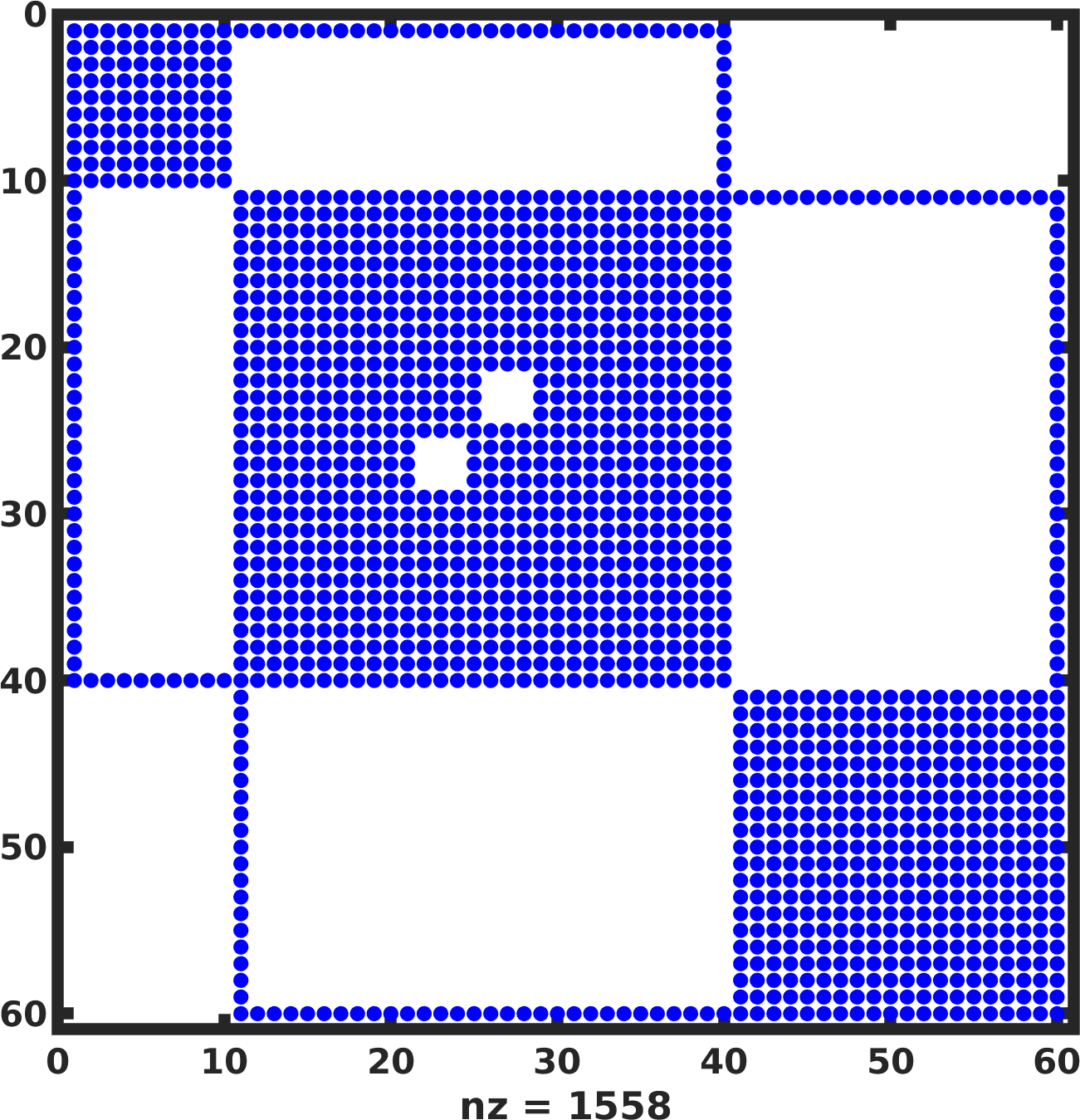} &
   \includegraphics[width=0.2\textwidth]{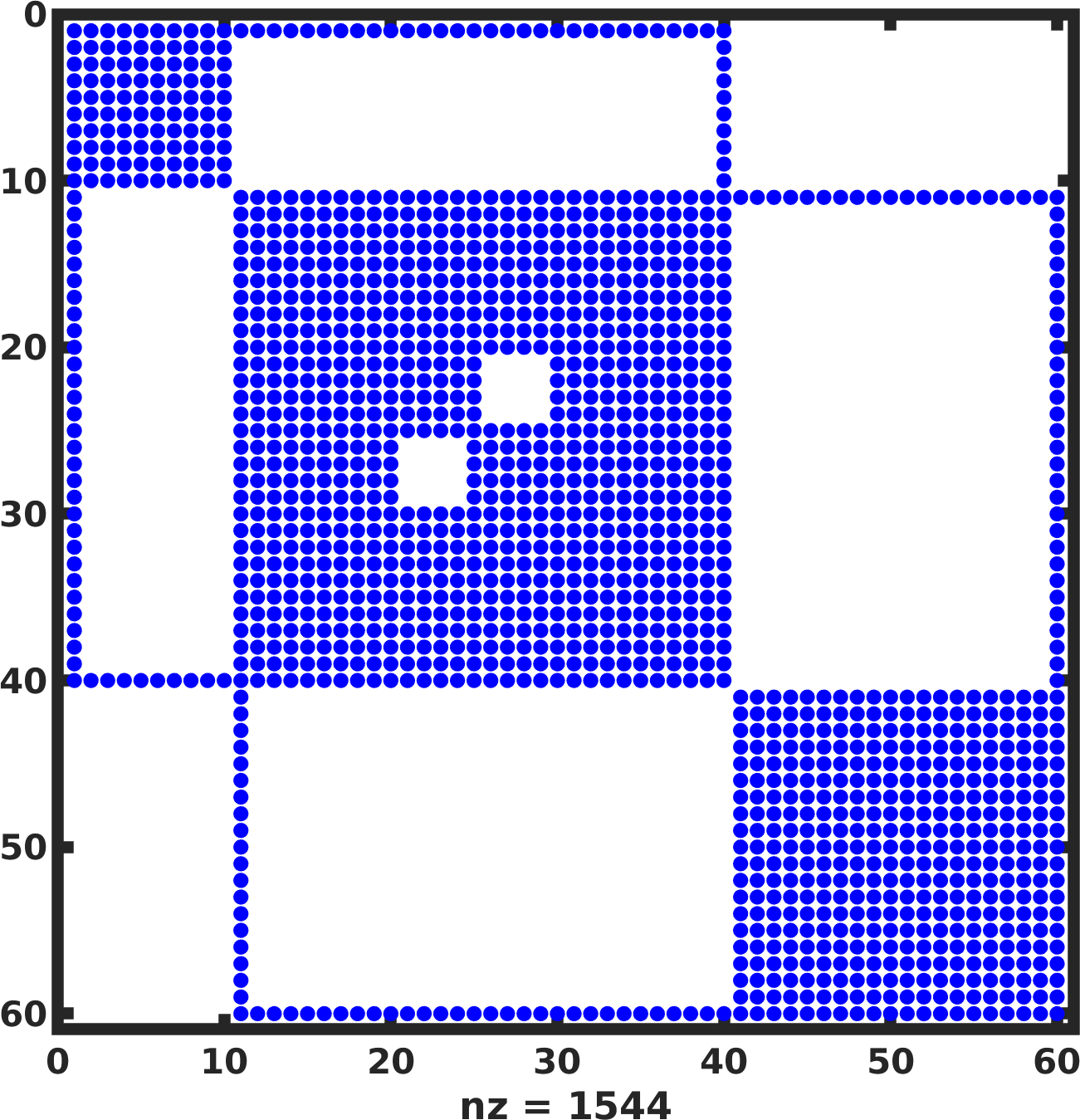}\\
    \includegraphics[width=0.2\textwidth]{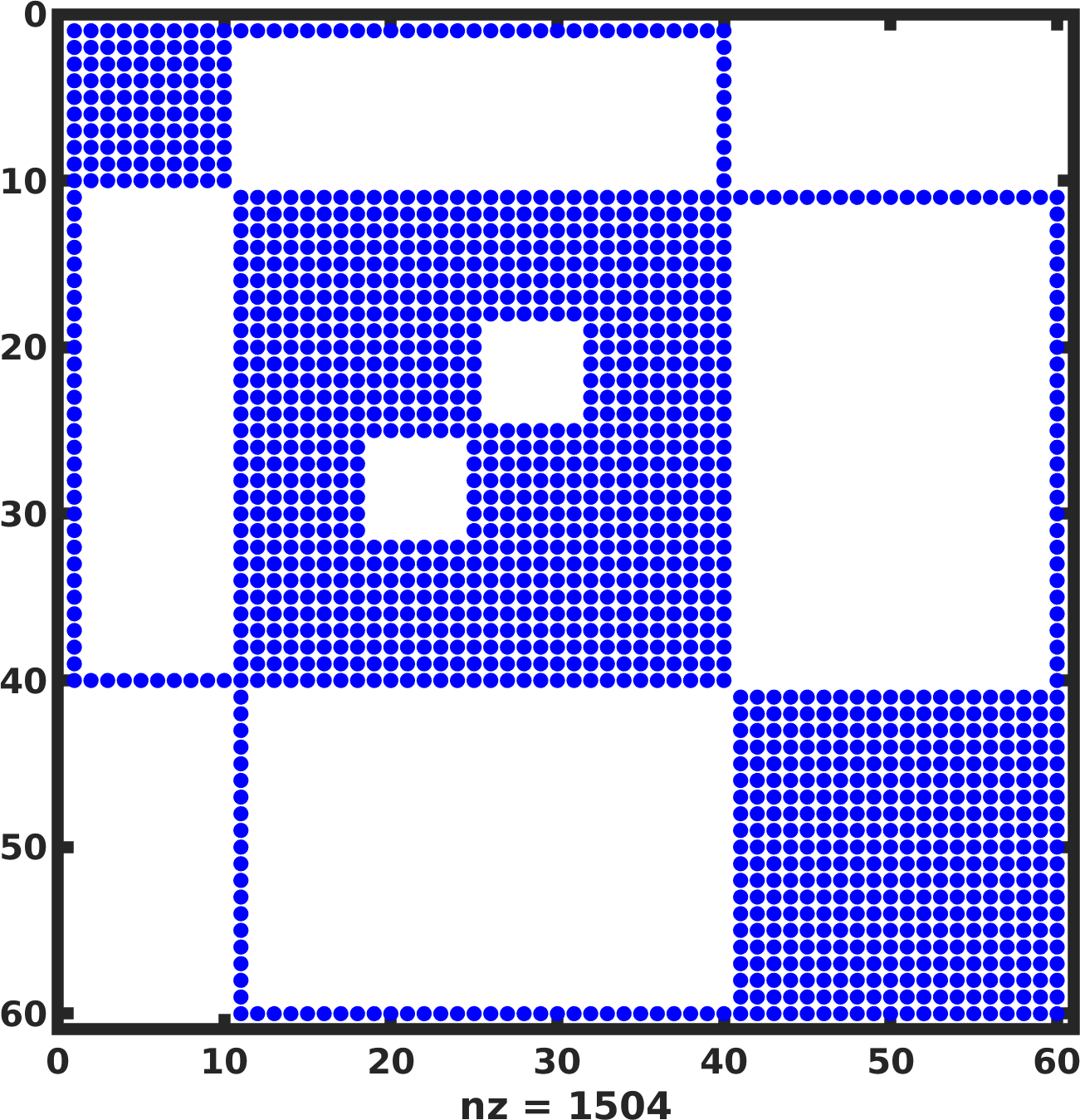}&
     \includegraphics[width=0.2\textwidth]{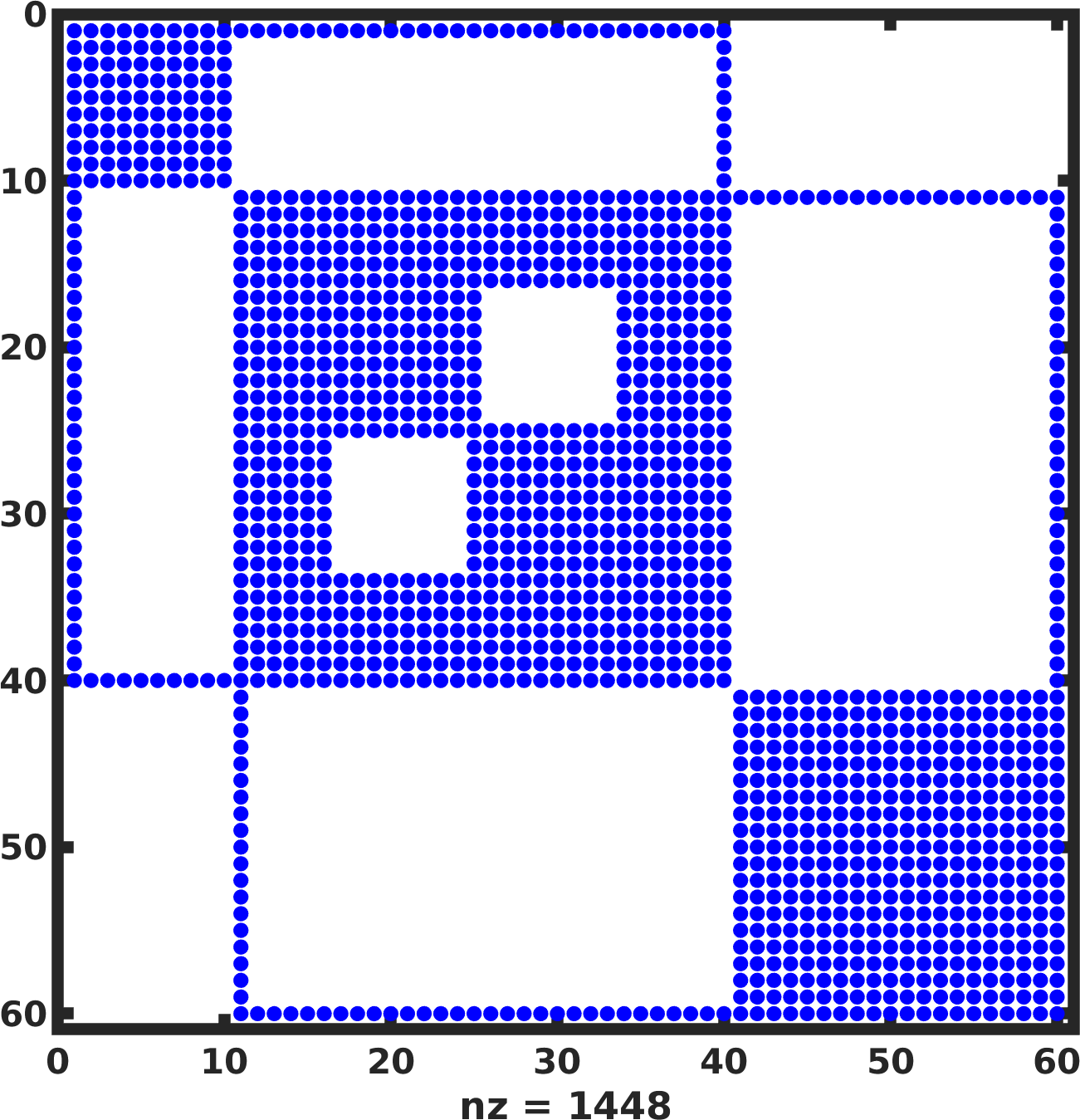}&
      \includegraphics[width=0.2\textwidth]{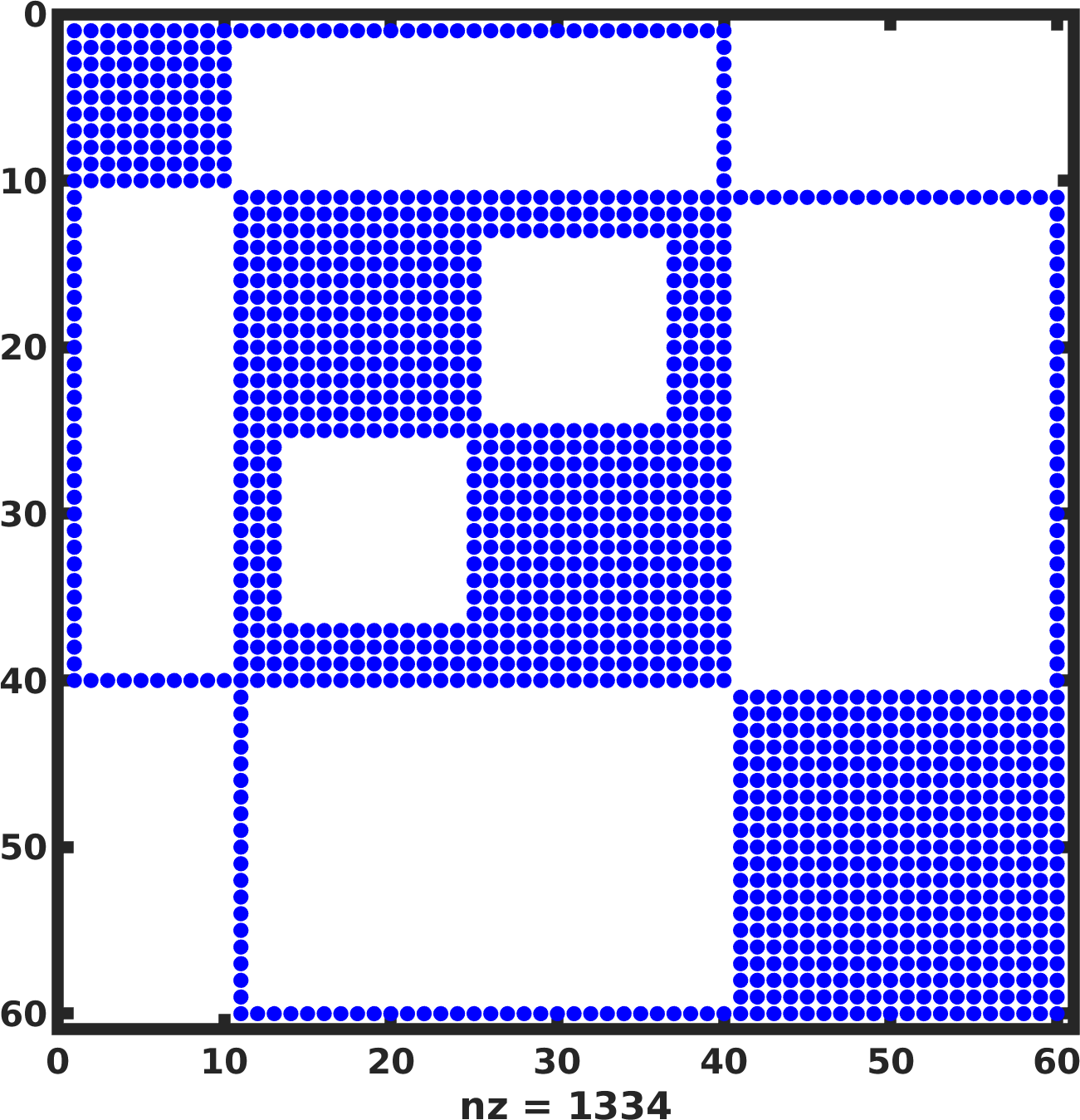}&
       \includegraphics[width=0.2\textwidth]{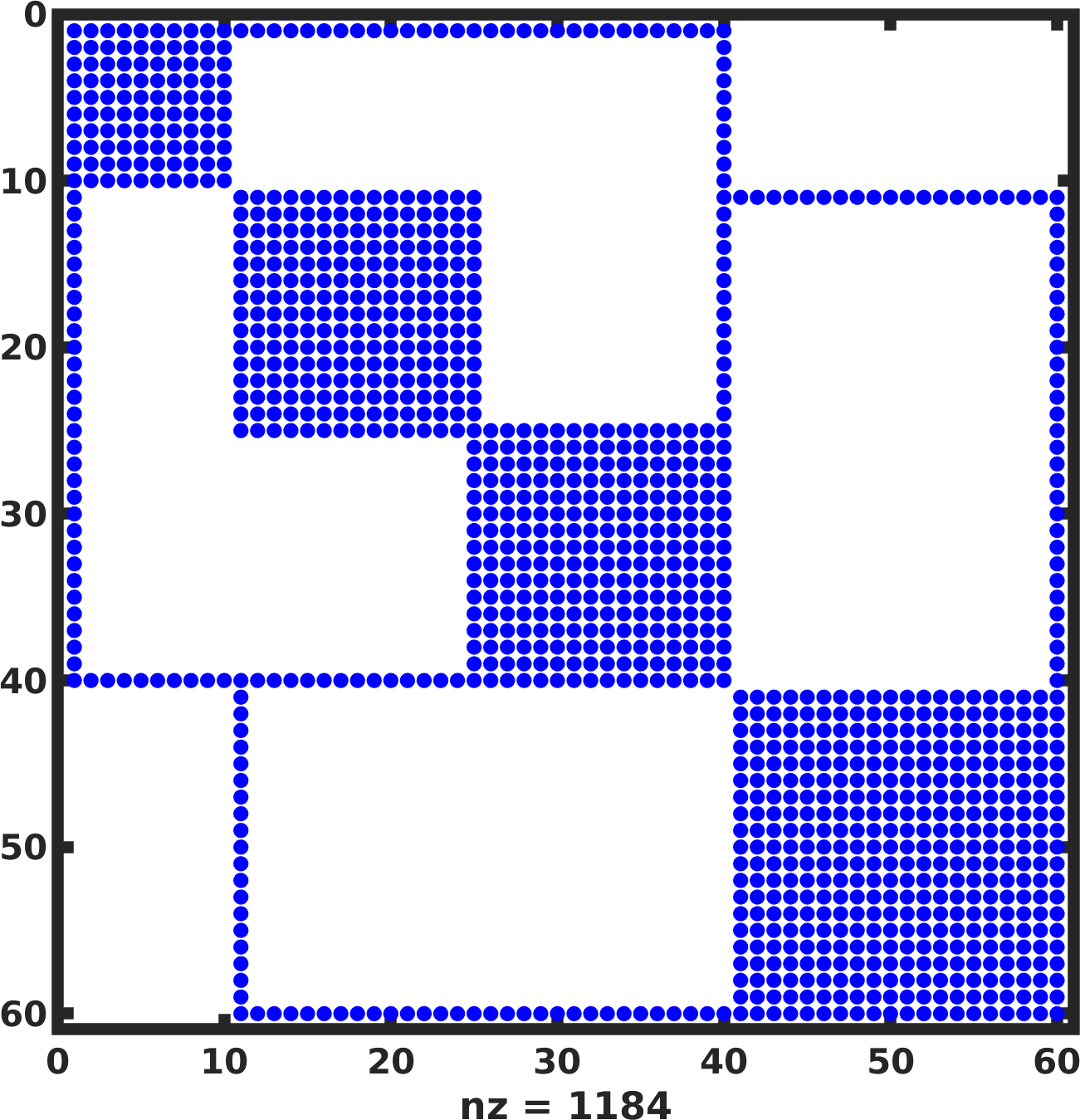}\\
\end{tabular}
\end{center}
\caption{\label{fig:picthfork} Different transition matrices of example \ref{seconExample}, where the split of the central almost-invariant pattern is provoked locally in its interior.}
\end{figure}  

\begin{figure}[!htb]
\centering
\includegraphics[width=.4\textwidth]{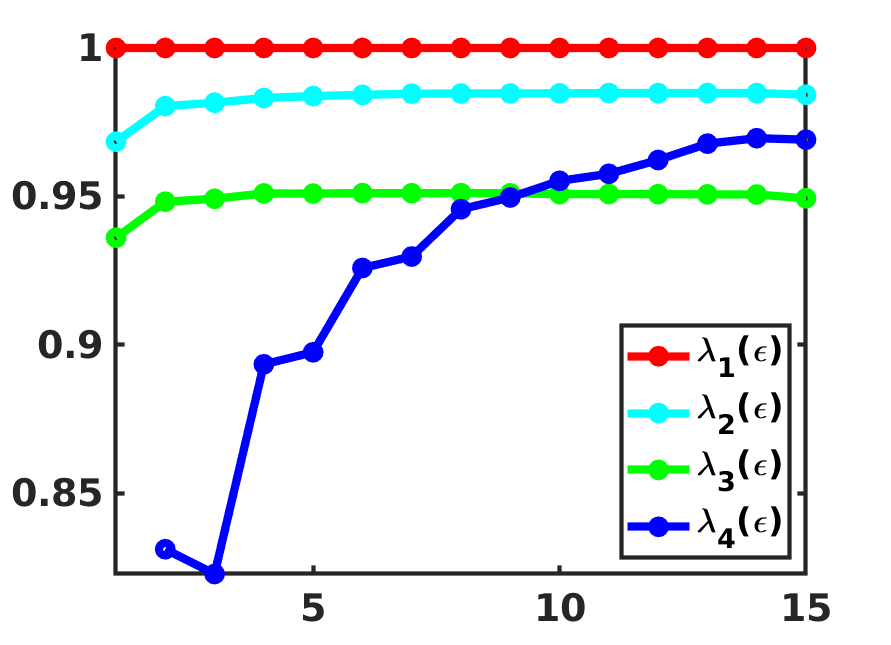}
\caption{\label{fig:pitch_eigs} Spectral signature of the splitting of an almost-invariant pattern in example \ref{seconExample}.}
\end{figure}

\begin{figure}[!htb]
\begin{center}
\begin{tabular}{cccc}
 \includegraphics[width=0.2\textwidth]{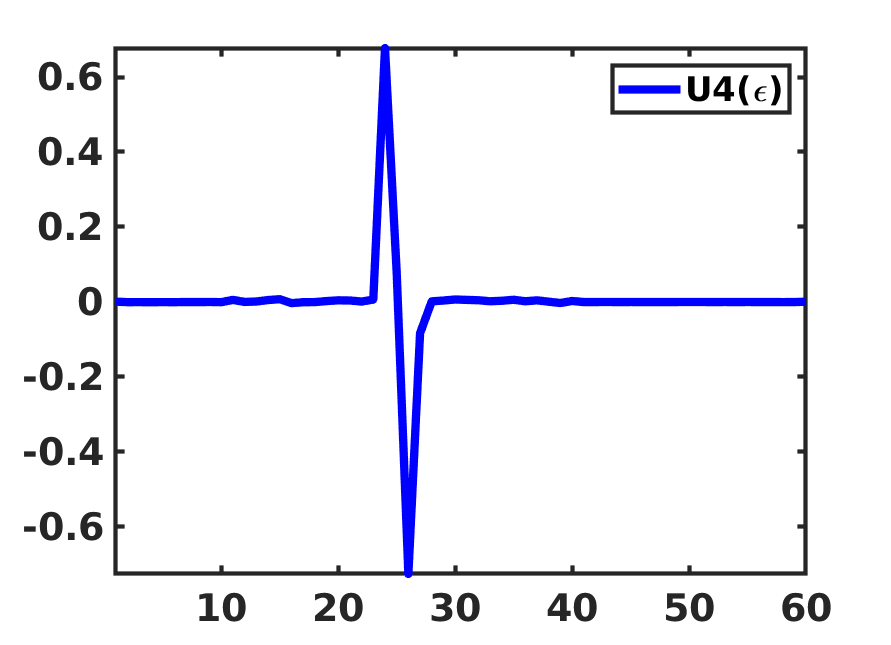}&
  \includegraphics[width=0.2\textwidth]{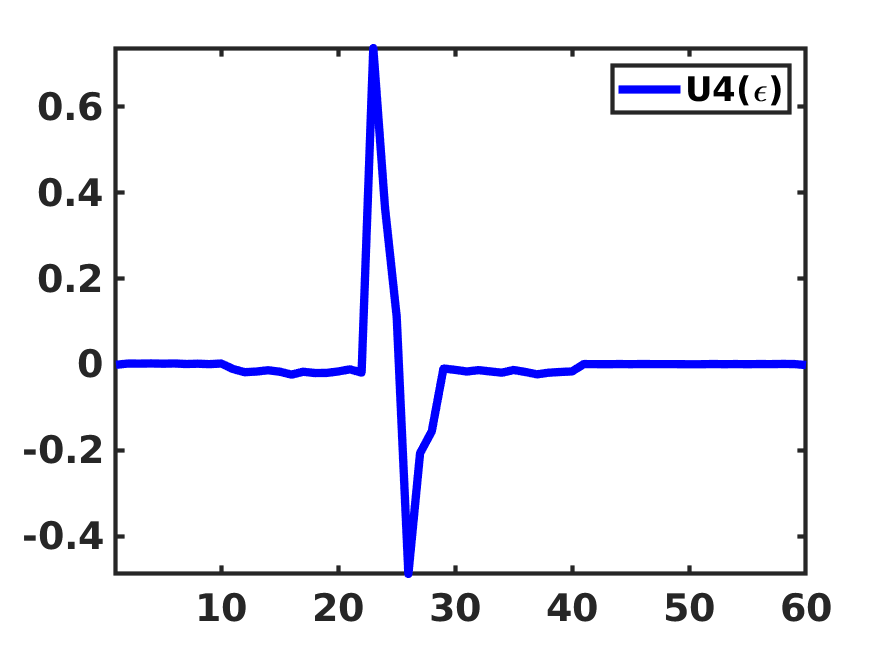}&
 \includegraphics[width=0.2\textwidth]{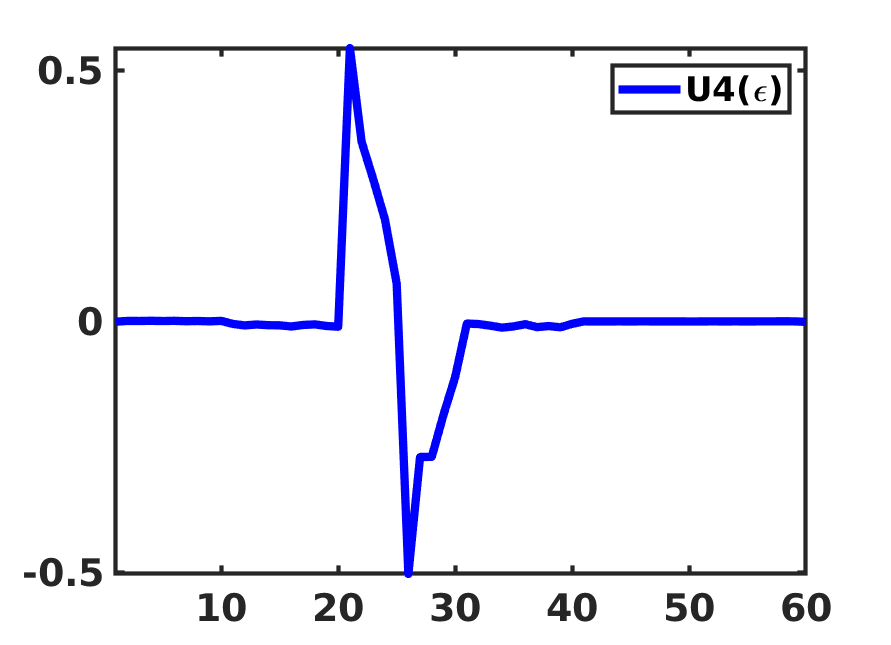}&
 \includegraphics[width=0.2\textwidth]{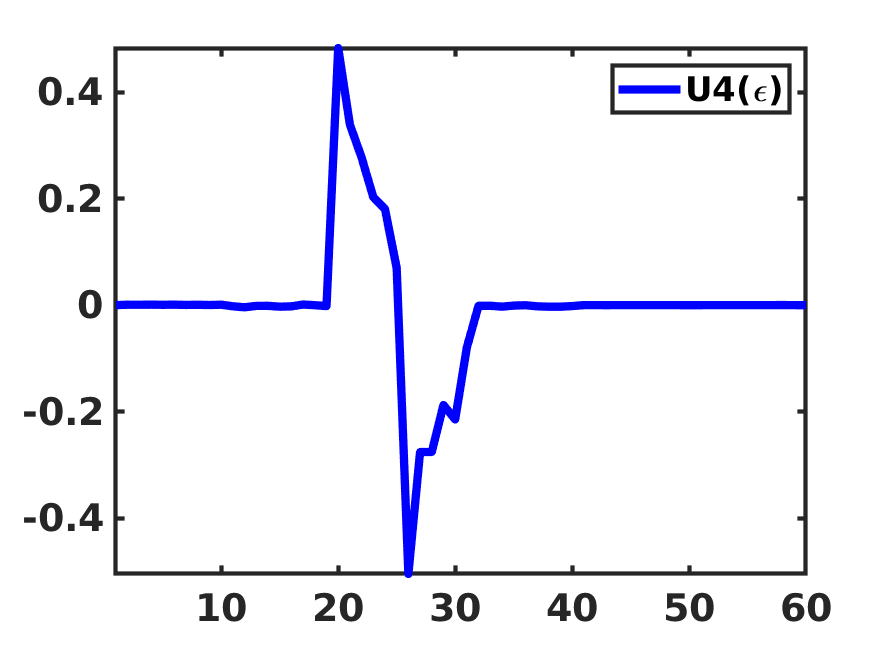}\\
 \includegraphics[width=0.2\textwidth]{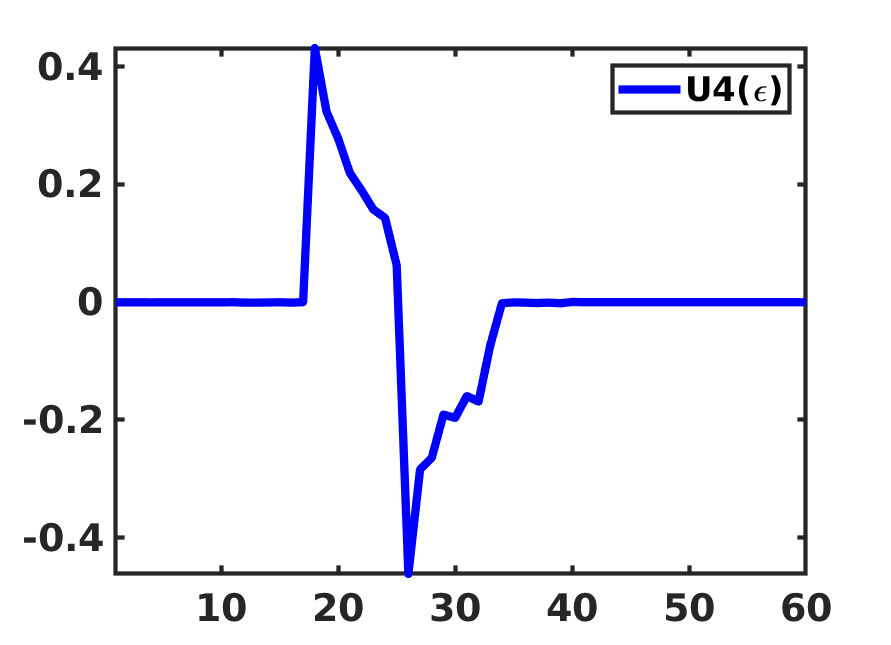}&
  \includegraphics[width=0.2\textwidth]{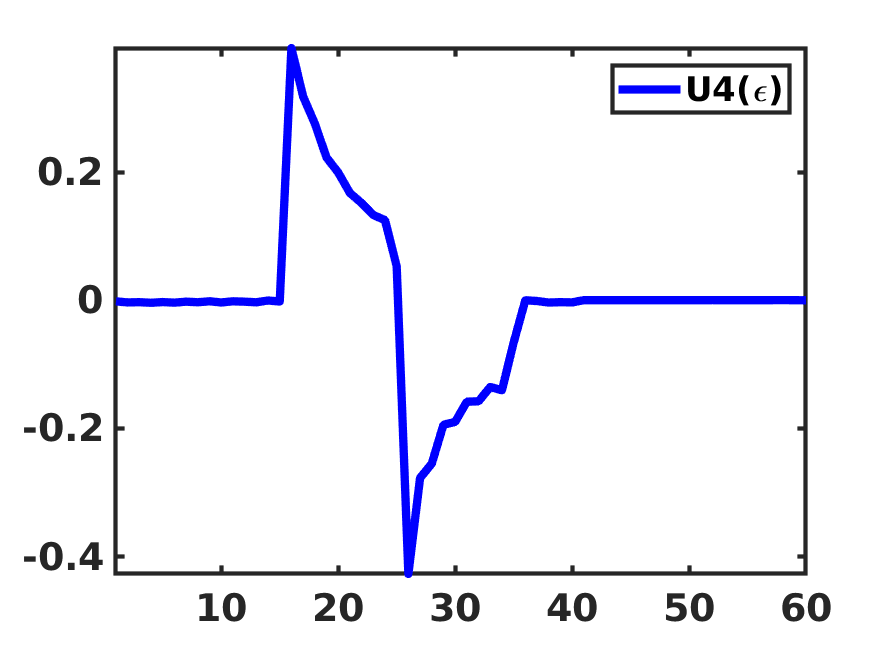}&
 \includegraphics[width=0.2\textwidth]{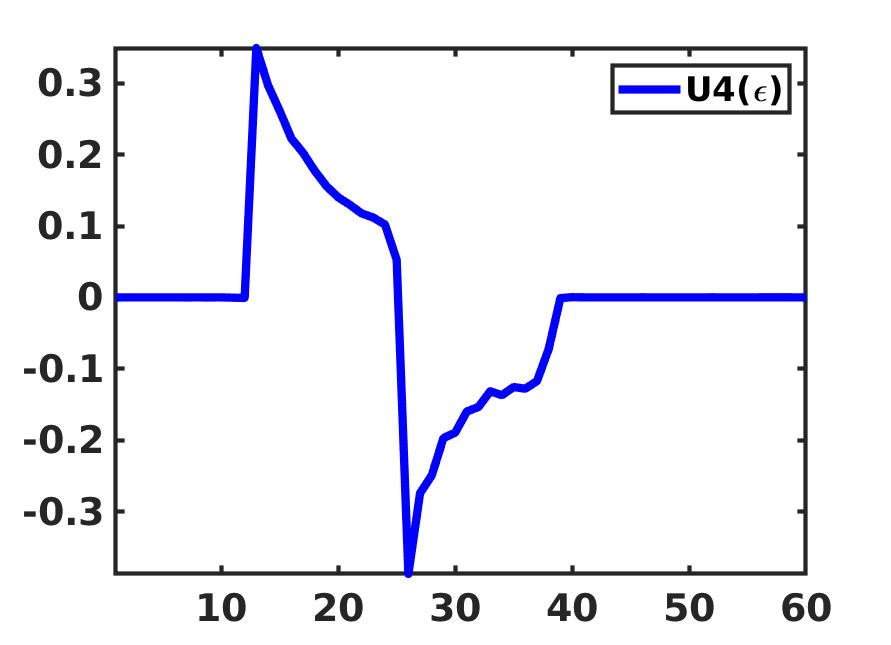}&
 \includegraphics[width=0.2\textwidth]{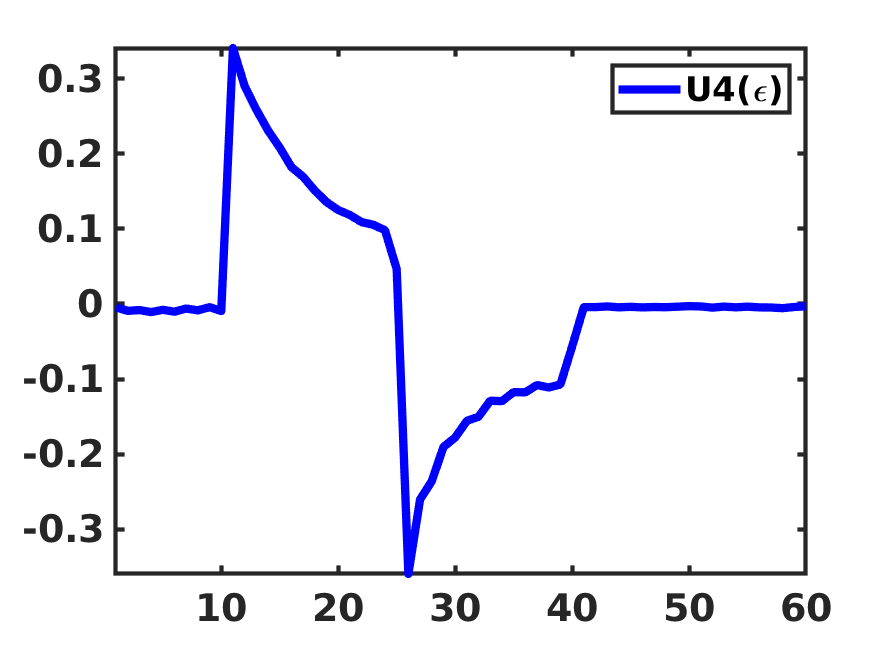}\\
\end{tabular}
\end{center}
\caption{\label{fig:pitch_evects} Changes of the previously subdominant eigenvector $U_k(\epsilon)$ (i.e. $k=4$) in example  \ref{seconExample}. }
\end{figure}  
Unlike the first example, there is no variation in the trends of the eigenvalues $\lambda_i(\epsilon),\,i=1,2,3$. This is because the 
shapes of the first three invariant patterns $A_i$, $i=1,2,3$ have not been affected by the sudden birth of the new pattern $A_4$. Therefore, 
in this experiment we clearly see that the trends of the three dominant eigenvalues are not relevant in order to predict the changes occurring in the dynamics. This can be 
understood by the fact that the change is primarily local and is only happening inside $A_3$. 
Again, the variation in the size of the almost-invariant patterns seems to be a crucial component for understanding the trends of the eigenvalues.
 \end{example}

\begin{example}\label{thirdExample}
Finally, in this experiment we summarize the behaviors observed in the two previous examples \ref{firstExample}-\ref{seconExample} within one toy model. 
At the beginning there are two coexisting almost-invariant patterns. Due to implicitly tuning an 
parameter, which is external to the model, a new pattern arises continuously inside one of these almost-invariant sets. While this new pattern grows, the 
two other almost-invariant sets shrink. This is captured in the behavior of the dominant eigenvalues, see  figure \ref{fig:third_toym_eigs}.
\begin{figure}[!htb]
\centering
\includegraphics[width=.4\textwidth]{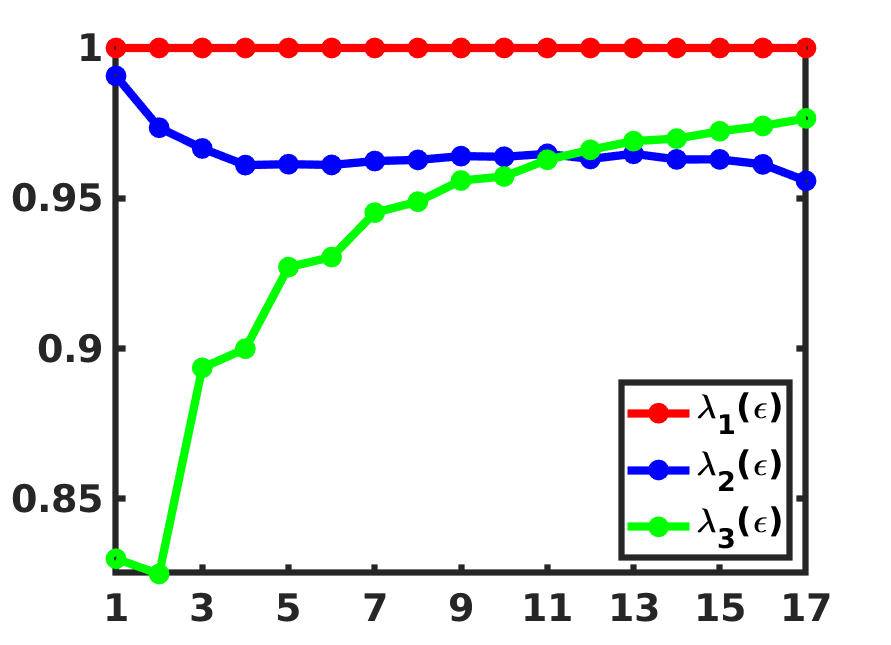}
\caption{\label{fig:third_toym_eigs} Spectral signature of a splitting of an almost-invariant pattern while another pattern is decreasing in size as described in example \ref{thirdExample}.}
\end{figure}
Indeed, as in equation \eqref{peigvec}, $\lambda_3(\epsilon)$ appears to rise from the small magnitude eigenvalues $\{\lambda_j(\epsilon),\, j=k+1,\ldots,N\}$, finally crossing 
$\lambda_2(\epsilon)$, which is decreasing. \end{example}

Although we have only shown very specific settings in examples \ref{firstExample}-\ref{thirdExample}, the spectral effects of the pattern changes that we have illustrated are universal.  In particular, it becomes clear that the study of qualitative changes of patterns, which are visible in the dominant eigenvectors, 
depends strongly on an understanding of the trends of corresponding eigenvalues. Moreover, any 
changing process within the almost-invariant patterns $A_i(\epsilon),\,i=1,\ldots,k$  will be first sensed in the smallest dominant eigenvector $U_k(\epsilon)$ and its corresponding eigenvalue. 
Indeed, the sign structure of the eigenvector $U_k(\epsilon)$ describes the $k$ existing almost-invariant patterns, exhaustively. For instance, in figure \ref{fig:pEV}, 
the $4^{th}$ eigenvector corresponds to the partition of the state space into four almost-invariant patterns. 
In particular, if $k=2$, then $U_2(\epsilon)$ partitions the state space into two almost-invariant patterns.
This particular case has been used in many works \cite{GFKPG, GFnonauto1, GFnonauto2} in the context of the numerical computation optimal almost-invariant sets 
from the global evolution of a dynamical system. The ultimate goal is to be able to recognize early warning signals of these critical changes of almost-invariant patterns. 

The trends of the eigenvalues and behavior of the state space as summarized in table \ref{fig:table} will facilitate the understanding of spectral behavior for more realistic systems.
 Note that here, we focused a lot on the splitting and/or shrinking behavior in state space, because we are ultimately interested in understanding such scenarios in real world systems. 
\begin{table}[!htb]
 \begin{center}
\begin{tabular}{|p{0.3\textwidth}|p{0.3\textwidth}|p{0.3\textwidth}|}
        \hline
        \textbf{Behavior in state space}    &    \textbf{Spectrum}   &   \textbf{Eigenvectors} \\
        \hline
        $A_i(\epsilon)$ shrinks  and disappears. At least one $A_j(\epsilon),\,j\neq i$ grows.   &  
$\lambda_i(\epsilon) \searrow$ while $\lambda_j(\epsilon)\nearrow$.        &   
        Support of $U_i(\epsilon)$ decreases, support of  $U_j(\epsilon)$  increases.\\
        \hline
	    $A_i(\epsilon)$ splits inside due to new  growing $A_{k+1}(\epsilon)$. $A_j(\epsilon)$, $j\neq i$ are unaffected.     & 
	    $\lambda_{k+1}(\epsilon)\nearrow$ and transport in $A_i(\epsilon)$ decreases due to increasing barrier  inside $A_i(\epsilon)$.  &             
	     $U_{k+1}(\epsilon)$ is  supported on growing new sets inside $A_i(\epsilon)$.\\
        \hline
        $A_i(\epsilon)$ shrinks because $A_{k+1}(\epsilon)$ increases  from  inside $A_i(\epsilon)$.         & 
         $\lambda_i(\epsilon)\searrow$ while $\lambda_{k+1}(\epsilon)\nearrow$, eventually crossing each other. \newline Then $\lambda_{k+1}(\epsilon)>\lambda_i(\epsilon)$.  &
    $U_{k+1}(\epsilon)$ is supported on growing new sets inside $A_i(\epsilon)$. Support of $U_i(\epsilon)$ decreases.                \\
        \hline
        \end{tabular}
        \end{center}
        \caption{\label{fig:table} Summarized results of the toy model experiments \ref{firstExample}-\ref{thirdExample}.}
        \end{table}
        
\section{Bifurcation of almost-invariant patterns}\label{bfmodels}
Now we study bifurcations of almost-invariant patterns generated by explicit mathematical models. We will consider the setting where there is initially a particular almost-invariant pattern centered at $(0,0)$ and surrounded by ring-like patterns, for each eigenvector of the $k$ dominant eigenvectors such as in figure \ref{fig:pEV}. The motivation for this is that this particular pattern mimics real world vortices.
Following our experiments in examples \ref{firstExample}-\ref{thirdExample}, which are summarized in table \ref{fig:table}, we will track the changes of 
the $k$ dominant eigenvectors and eigenvalues with respect to a concrete external bifurcation parameter $p$. Thus, for the feasibility of this continuation 
task, we assume a fixed perturbation strength $\epsilon$ during all of the process. That is, the variations of the $k$ eigenvalues and eigenvectors 
will only depend on the bifurcation parameter $p \in \mathbb{R}$. 

\begin{figure}[h!]
\begin{center}
 \includegraphics[width=0.4\linewidth]{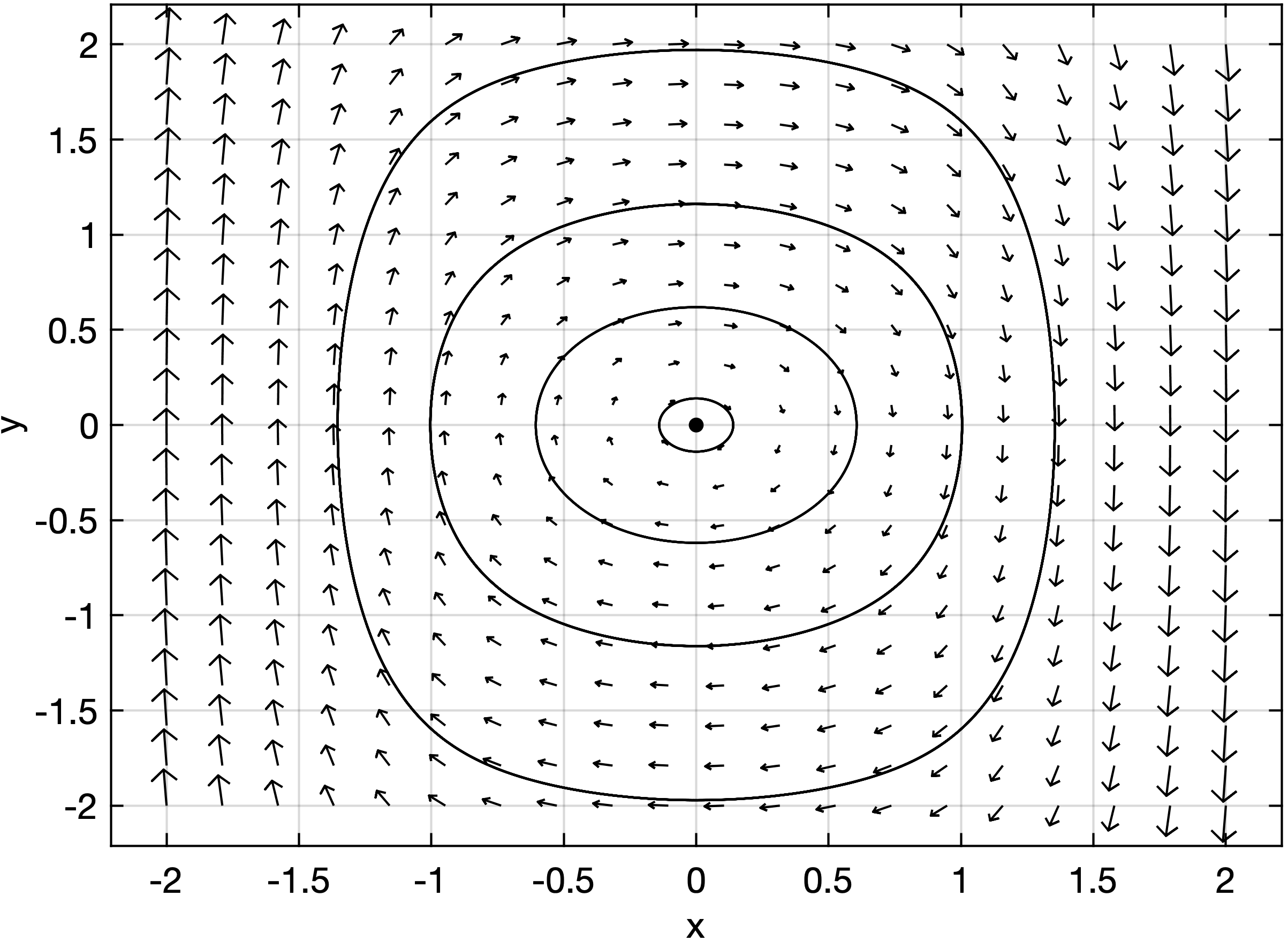}
 \includegraphics[width=0.4\linewidth]{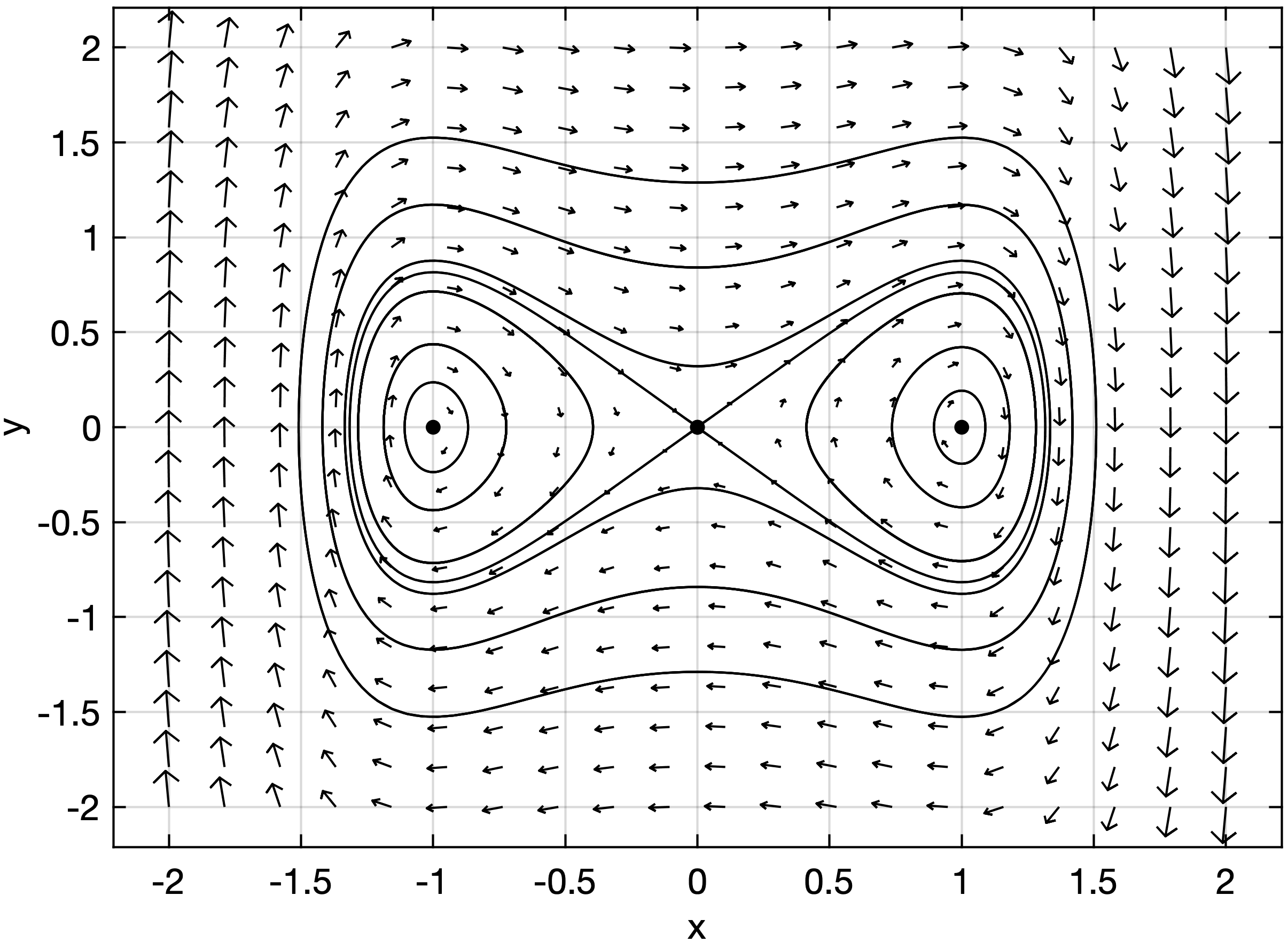}\\
 \end{center}
  \caption{Phase plane of system \eqref{psys1} for parameters $p=-1$ (left) and $p=1$ (right). }\label{fig:phaseplot}
 \end{figure}

As a first case study, we consider the $p$-parametrized two-dimensional system 
\begin{equation}\label{psys1}
\begin{aligned}
   \dot{x} & = y\\
   \dot{y} & =p x - x^5
\end{aligned}
\end{equation}
$p\in \mathbb{R}$. System \eqref{psys1} is a conservative Duffing-type oscillator. It is well known that its classical bifurcation consists of the 
qualitative change of the unique elliptic fixed point $(0,0)$, for $p<0$, into a local saddle fixed point, for $p>0$. That is, a pitchfork bifurcation occurs when $p=0$, which has global effects on the dynamics. For $p<0$, the stationary dynamics consists of rotating periodic orbits centered at the unique fixed point $(0,0)$. 
These are destroyed, for $p>0$, with the emergence of two symmetric elliptic fixed points at $(\sqrt[\leftroot{-2}\uproot{2}4]{p},0)$ and $(-\sqrt[\leftroot{-2}\uproot{2}4]{p},0)$, as illustrated 
in figure \ref{fig:phaseplot}. 

To prepare for our spectral analysis, a $p$-parametrized version of equation \eqref{peigvec} may now be restated as follows:
For each $p$, the stationary distribution is given as
   $$ \pi(\epsilon,p) \equiv U_1(\epsilon,p) =[\pi_1(\epsilon,p),\pi_2(\epsilon,p),\ldots,\pi_N(\epsilon,p)],\,\,\pi_i(\epsilon)>0,\,\,\,\forall\,\,\, p,$$  
and for each $i=2,\ldots,k$, 
\begin{equation}\label{ppmueigvec}
\begin{cases}
\begin{split}
    U_i(\epsilon,p) &= \sum_{j=1}^k (\alpha_{ij}+ \epsilon\beta_{ij}) V_j(p)\\
                     &+ \epsilon\sum_{j=k+1}^N \frac{1}{1-\lambda_j(\epsilon,p)}\langle U_j(p), Q^{(1)}U_i(p)\rangle_{\pi(\epsilon,p)} + 
                     O(\epsilon^2),\,\,\alpha_{ij},\beta_{ij} \in \mathbb{R},\\
                     \\
     \lambda_i(\epsilon,p) & > \lambda_j(\epsilon,p),\, j=k+1,\ldots,N.
\end{split}
\end{cases}
\end{equation}
Note that with a fixed $\epsilon$, the additional inequality constraint in \eqref{ppmueigvec}
\begin{equation}\label{degeneracy_eq}
\lambda_i(\epsilon,p)  > \lambda_j(\epsilon,p),\,\, i=1,2,\ldots,k,\,j\geq k+1
\end{equation}
is always satisfied whenever the changes in $p$ 
leave the qualitative behavior of system \eqref{psys1} unaffected. Indeed, due to the perturbation effect,  $\lambda_j(\epsilon,p)<1$, $j=k+1,\ldots,N$ are the 
small magnitude real eigenvalues which 
converge to 0 when $\epsilon$ increases. However, when $\epsilon$ is fixed, the changes in $p$ may qualitatively affect the underlying dynamics.
Thus, it makes sense to 
measure a susceptible radical growth scenario of the $\lambda_j(\epsilon,p)<1$, $j=k+1,\ldots,N$, among many other possible scenarios. 

\subsection{Spectral signature of the classical bifurcation}
Unlike the dominant eigenvectors (as shown in figure \ref{fig:pEV} for system \eqref{sys2}) the remaining $N-k$ eigenvectors $U_j(\epsilon,p)$, $j=k+1,\ldots,N$ 
may not be supported on the whole state space. They are referred as the "weak modes" eigenvectors and may not carry 
dynamically useful information, compared to the $k$ "dominant modes" eigenvectors. However, due to the nature of the global behavior of \eqref{psys1} illustrated 
in figure \ref{fig:phaseplot}, the global classical bifurcation yields a radical change only within a local isolated neighborhood of $(0,0)$. We refer to the latter as 
the critical neighborhood $\mathcal{D}$. Indeed, far from $\mathcal{D}$, closed 
trajectories still remain qualitatively the same before and after the bifurcation; see figure \ref{fig:phaseplot}. Therefore, we will first find a spectral version of the 
classical bifurcation by means of the non-dominating $N-k$ part of the spectrum. That is, we will consider "weak modes" eigenvectors which are only supported on $\mathcal{D}$. 
Note that a special technique to finding those particular eigenvectors is still an open question. Their existence was noticed earlier in \cite{Oliver}, 
but no particular further study about them was made, whatsoever. In this work, we use them to design a spectral bifurcation diagram of the global classical 
bifurcation occurring in \eqref{psys1}. They will also play an import role when studying the bifurcation of "dominant mode" eigenvectors pattern.

\begin{figure}[!htb]
\begin{center}
\begin{tabular}{cccc}
 \includegraphics[width=0.2\textwidth]{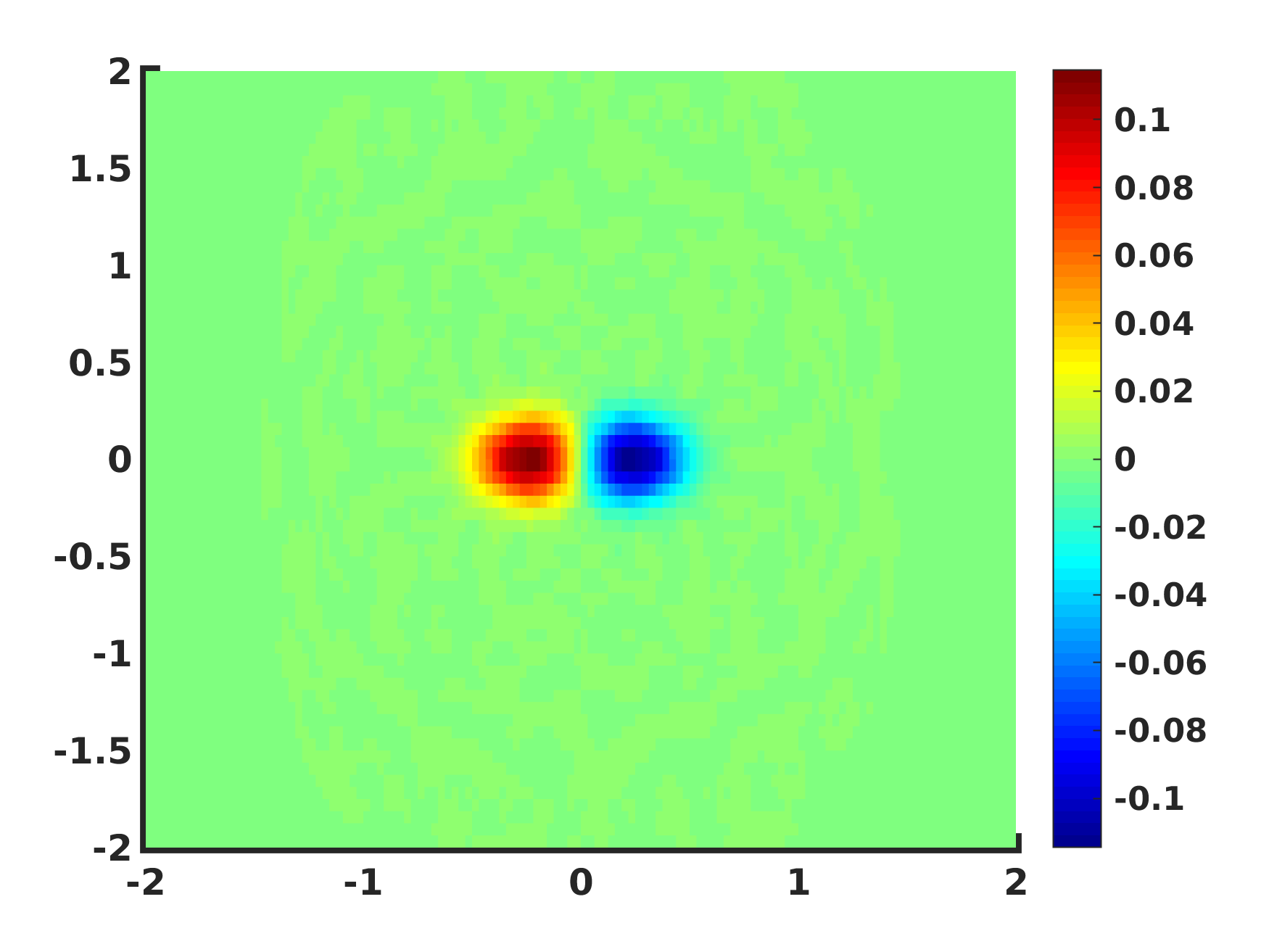}&
 \includegraphics[width=0.2\textwidth]{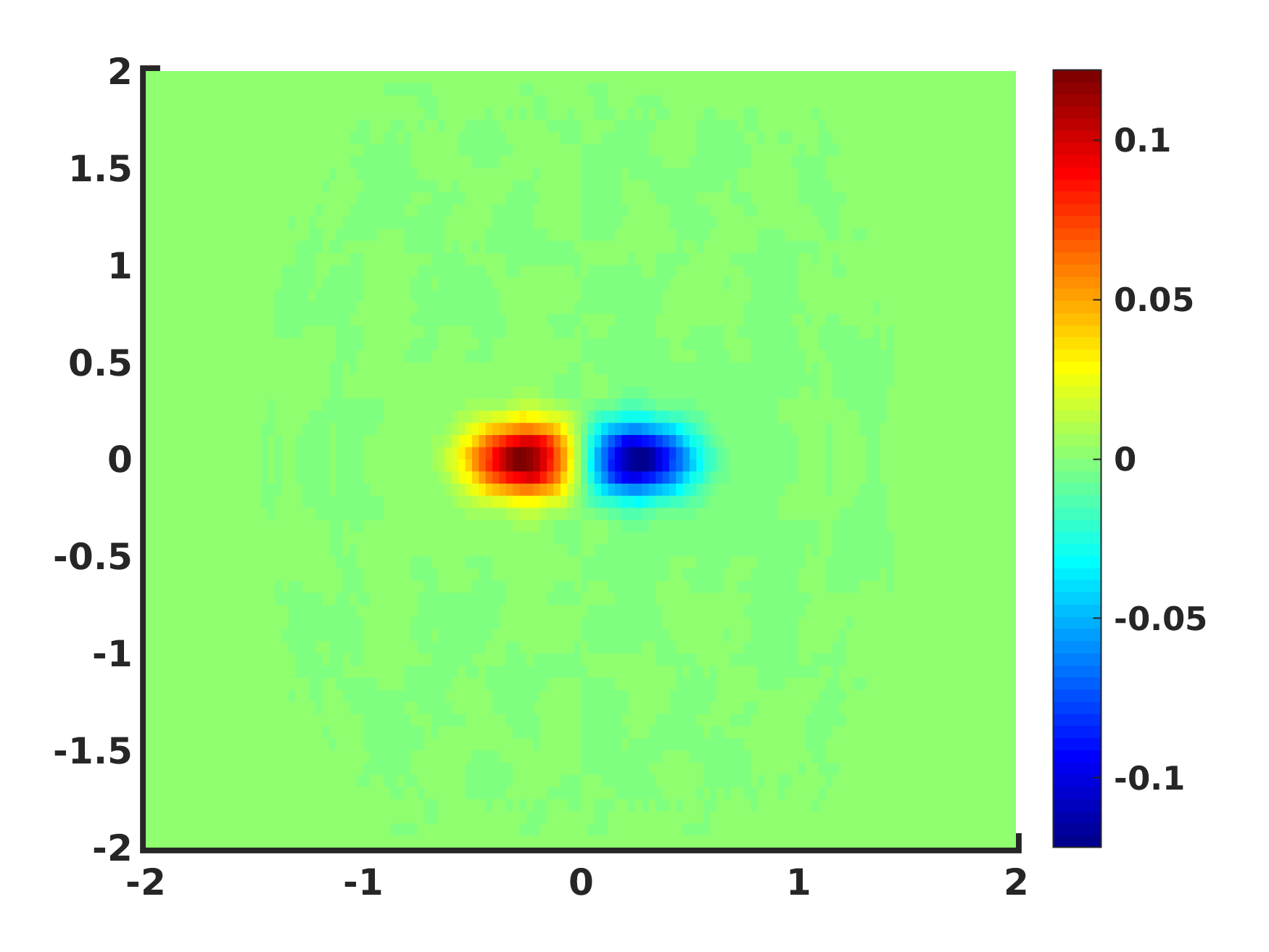} &
  \includegraphics[width=0.2\textwidth]{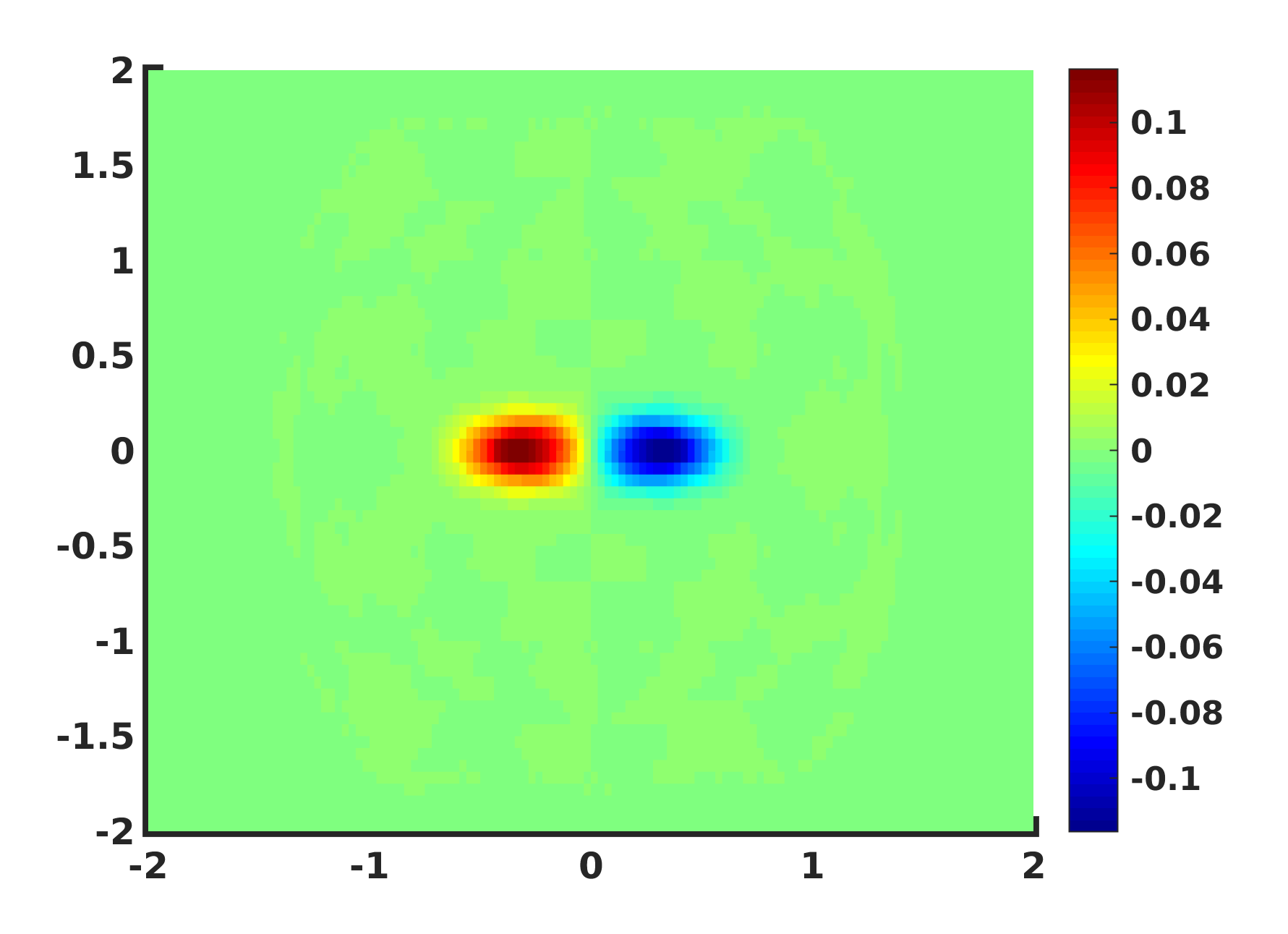} &
   \includegraphics[width=0.2\textwidth]{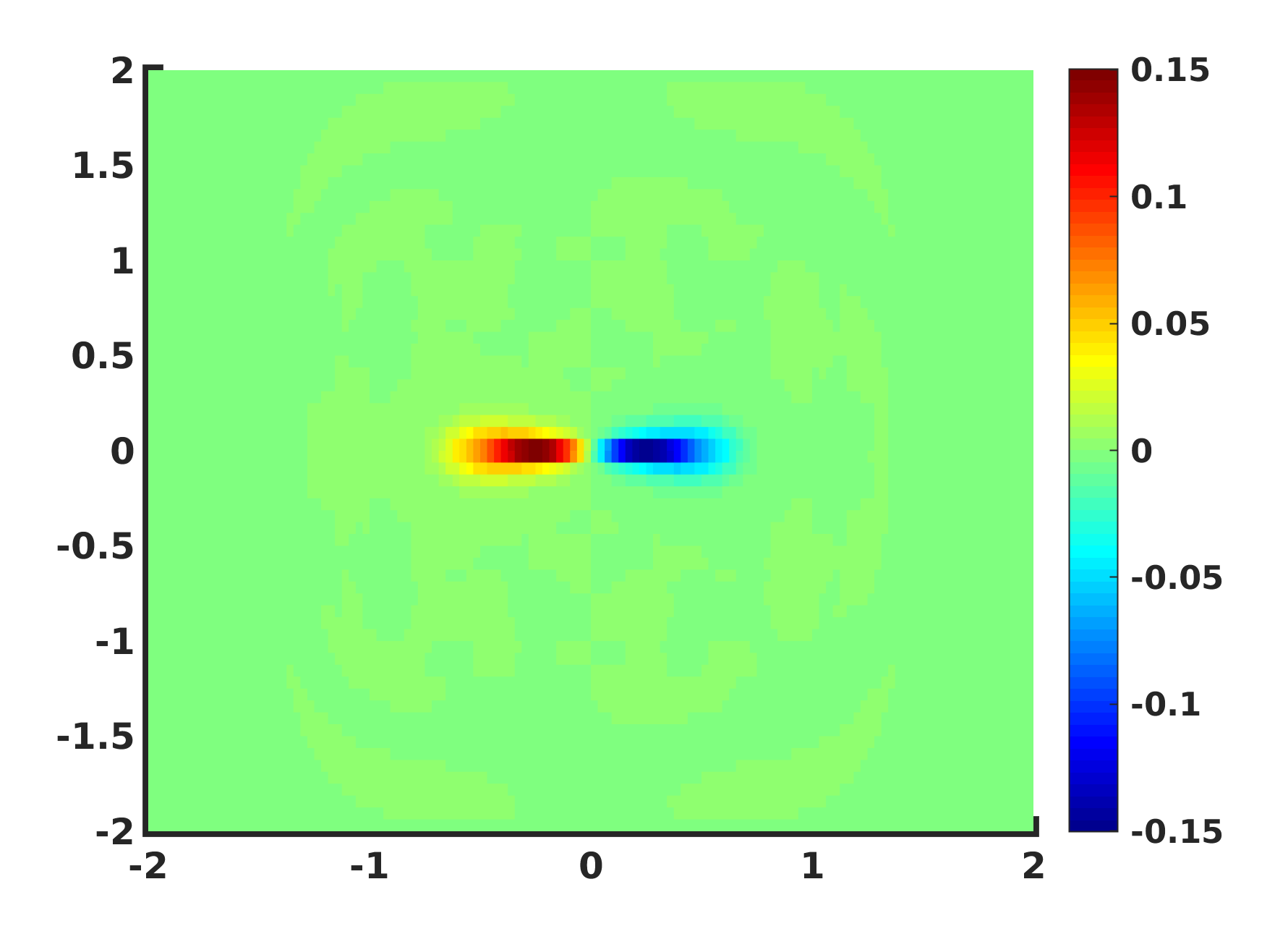} \\
   {\scriptsize (a) } &   {\scriptsize (b) } &   {\scriptsize (c) } &   {\scriptsize (d) } \\
    \includegraphics[width=0.2\textwidth]{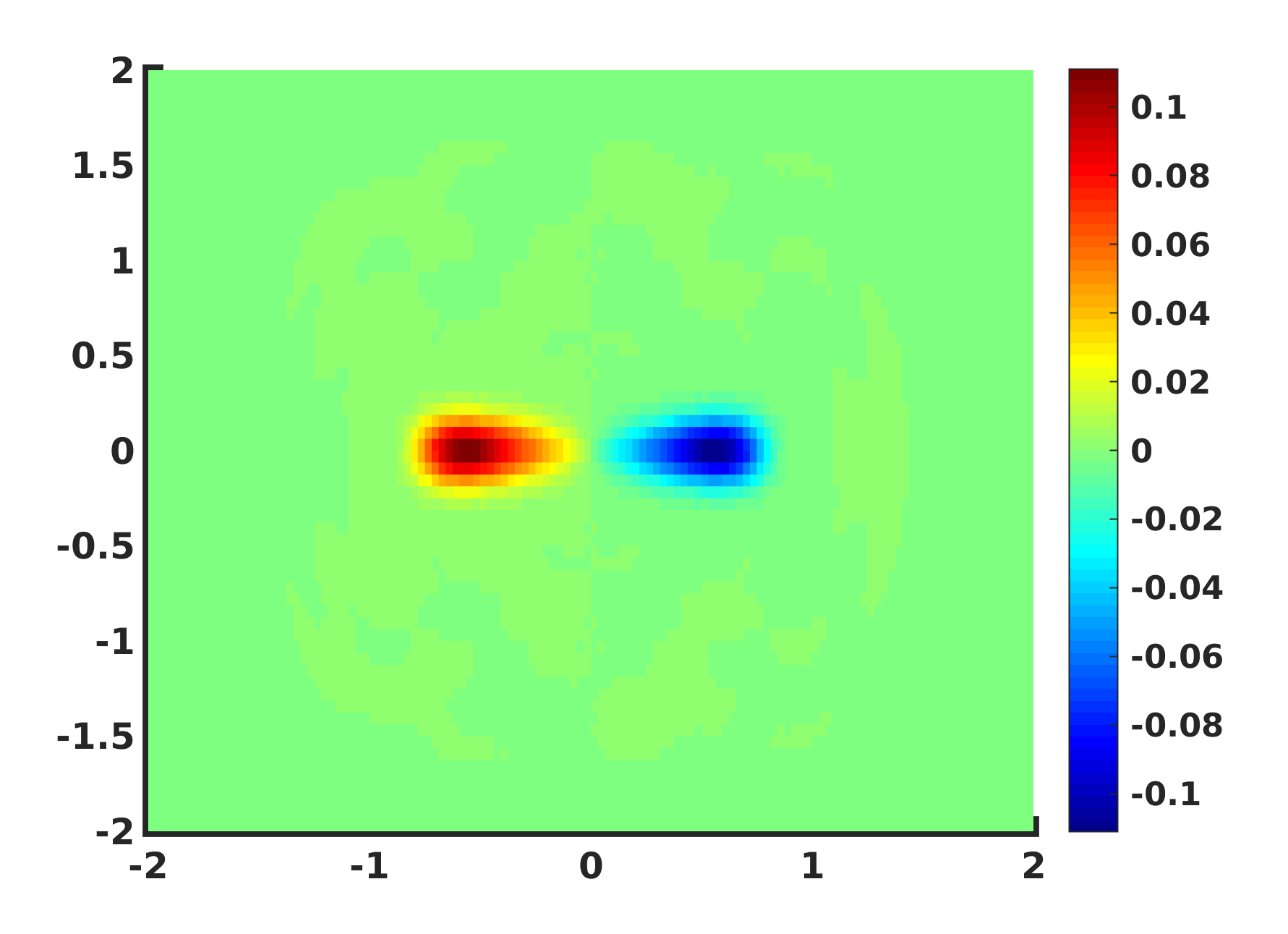}&
     \includegraphics[width=0.2\textwidth]{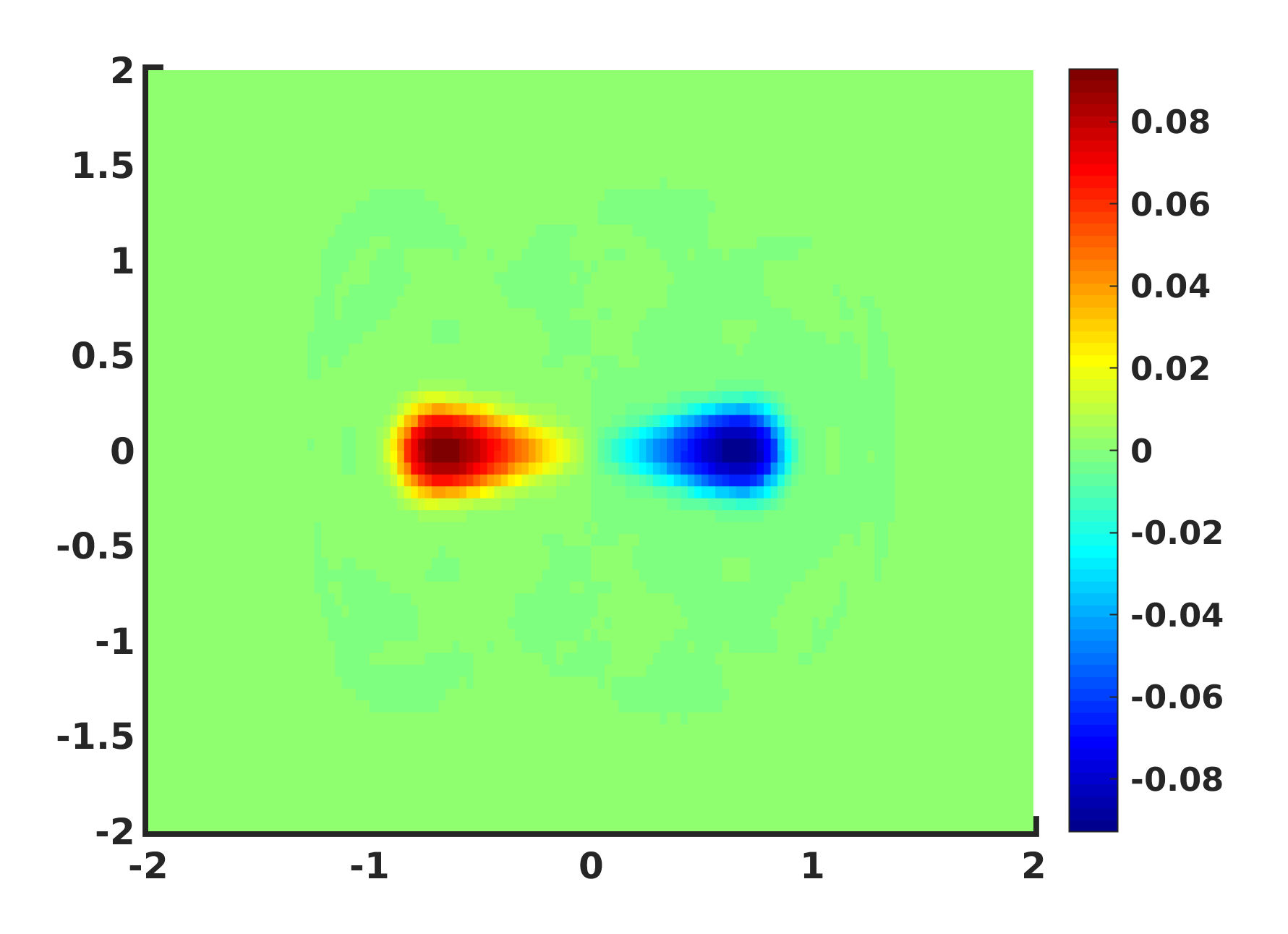}&
      \includegraphics[width=0.2\textwidth]{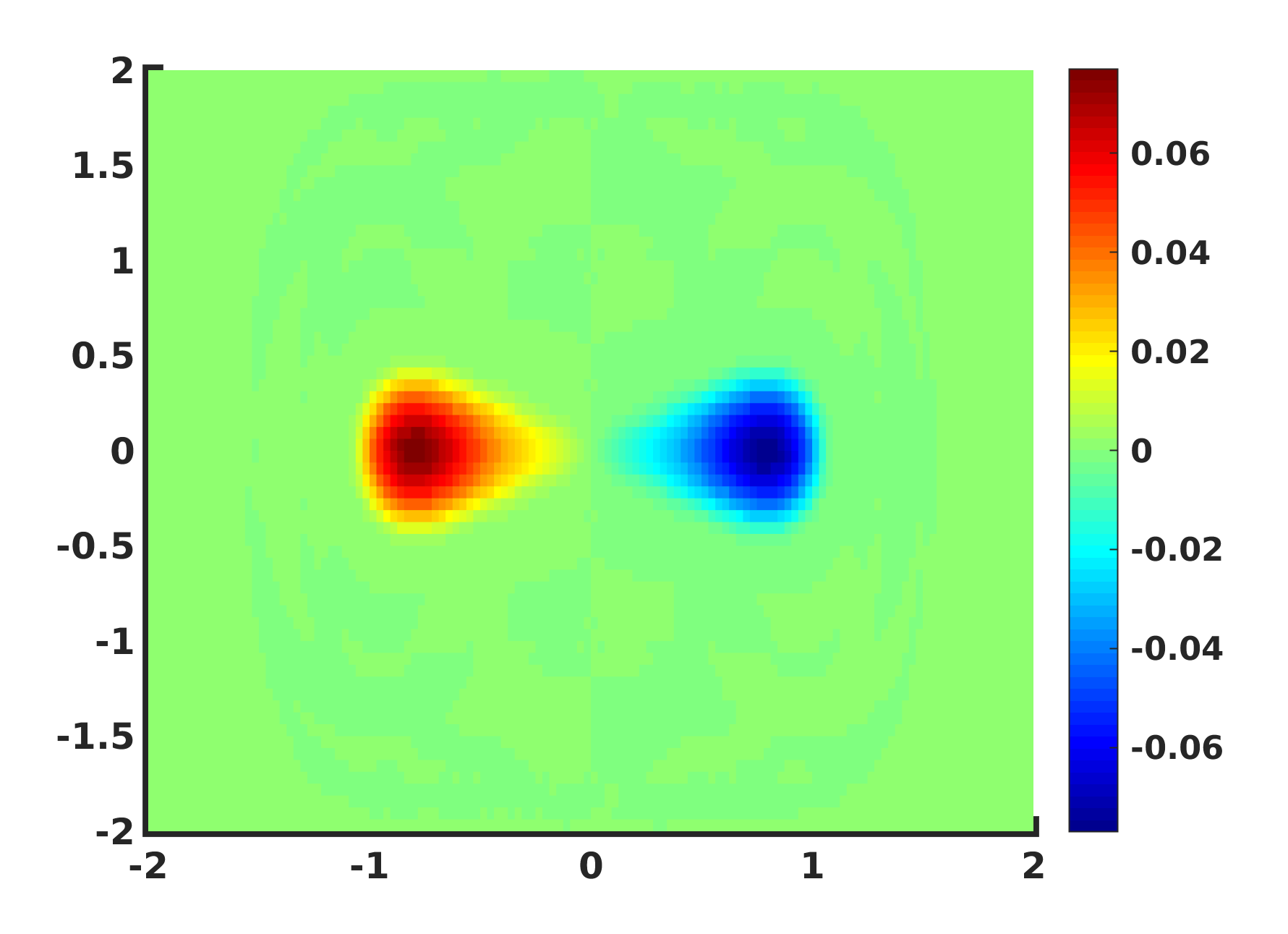}&
       \includegraphics[width=0.2\textwidth]{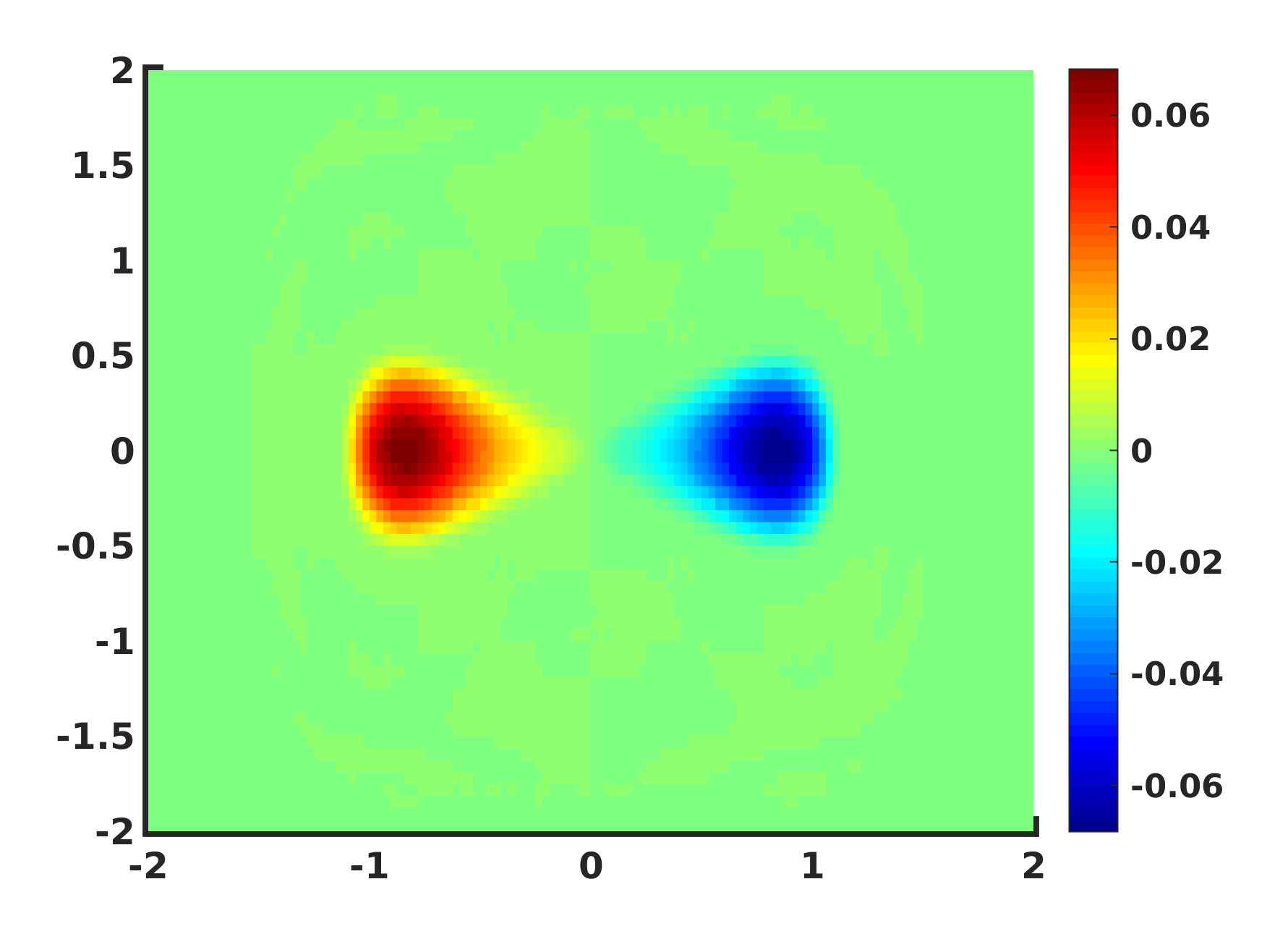}\\
        {\scriptsize (e) } &   {\scriptsize (f) } &   {\scriptsize (g) } &   {\scriptsize (h) } \\
\end{tabular}
\end{center}
\caption{\label{fig:pt1} Changes of the first "weak modes" patterns from $U_1^\prime(\epsilon,p)$ for increasing $p$.}
\end{figure}

\begin{figure}[!htb]
\begin{center}
\begin{tabular}{cccc}
  \includegraphics[width=0.2\textwidth]{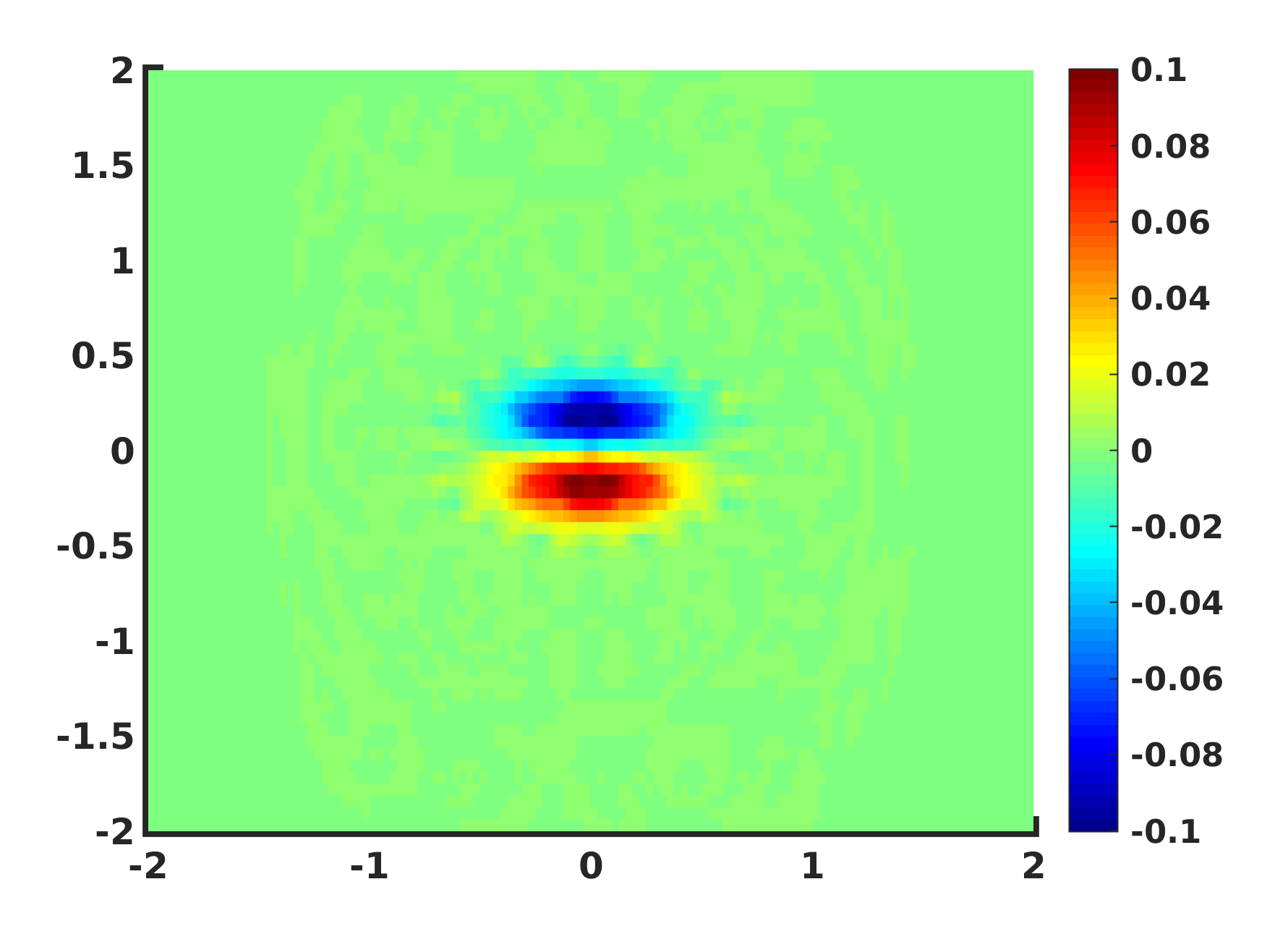} &
   \includegraphics[width=0.2\textwidth]{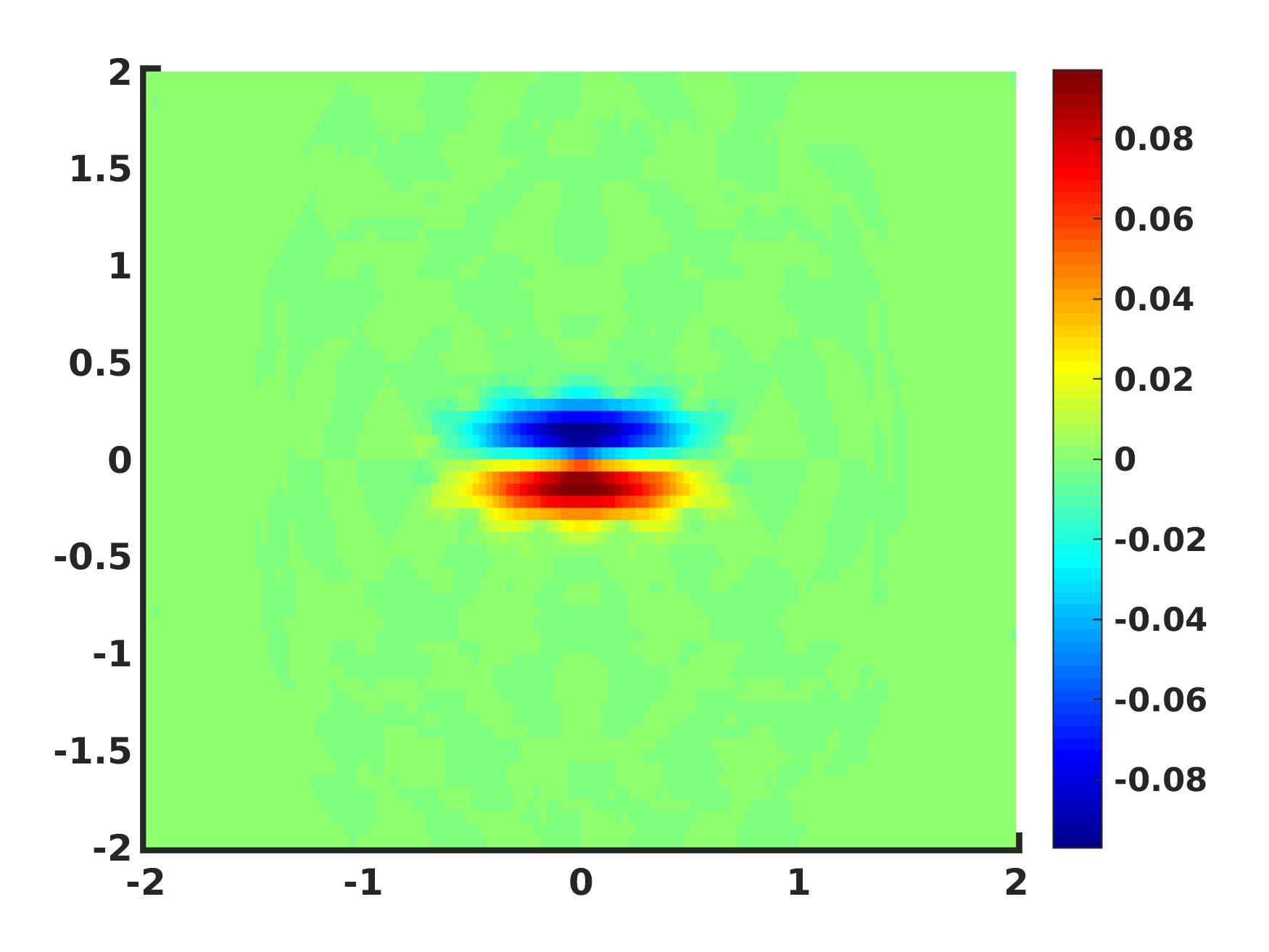} &
 \includegraphics[width=0.2\textwidth]{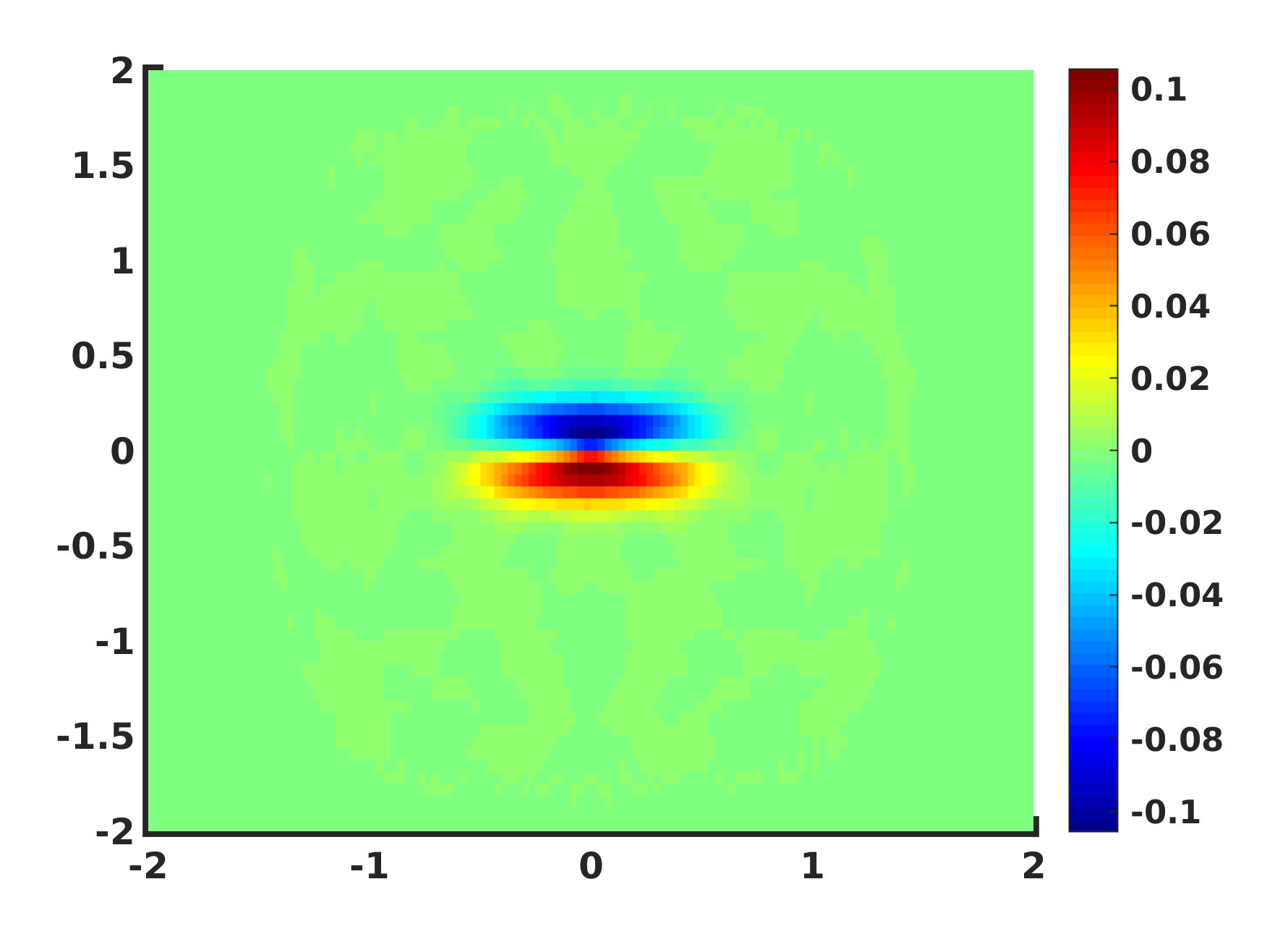} &
 \includegraphics[width=0.2\textwidth]{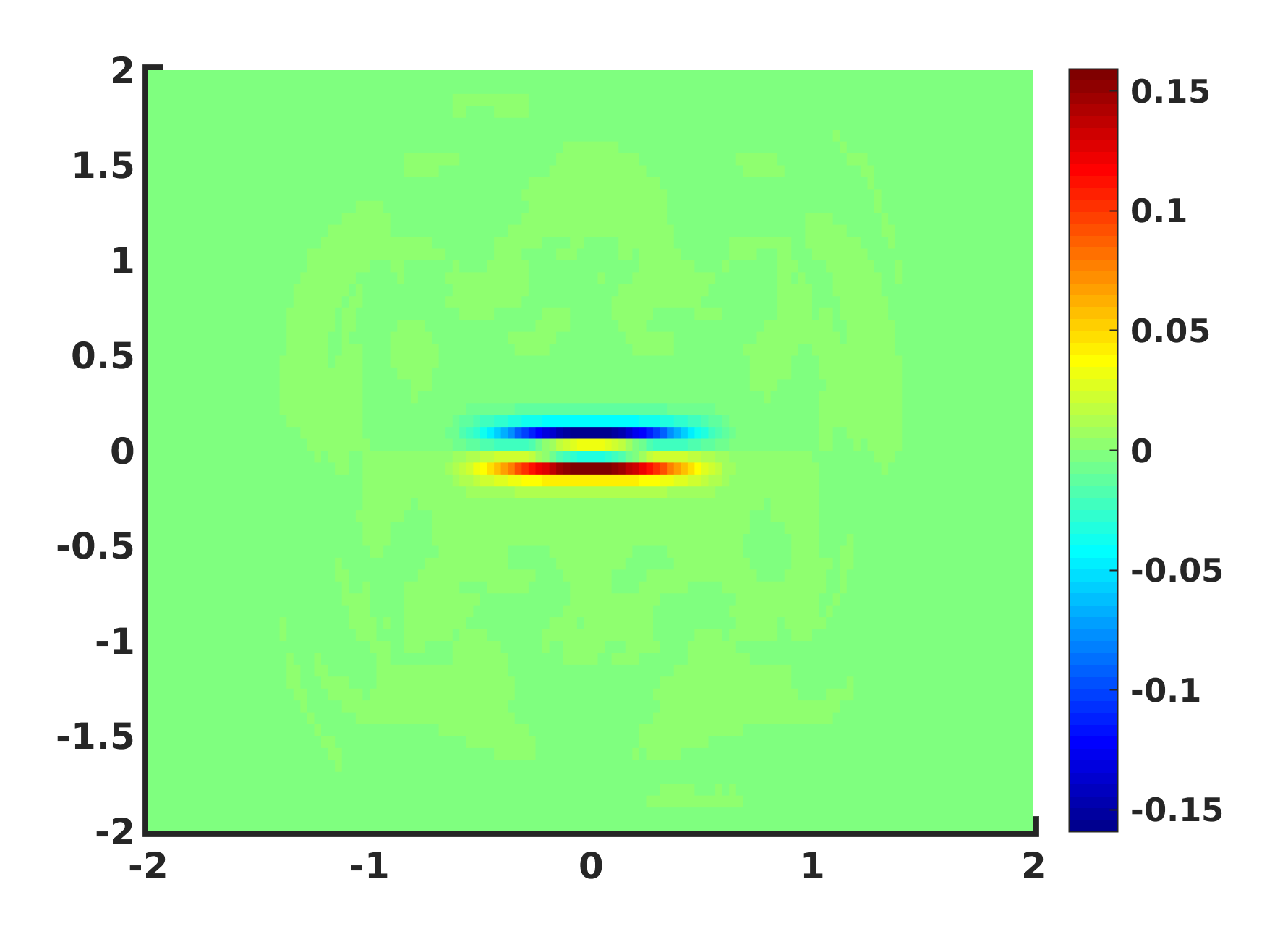}\\
  {\scriptsize (a) } &   {\scriptsize (b) } &   {\scriptsize (c) } &   {\scriptsize (d) } \\
\includegraphics[width=0.2\textwidth]{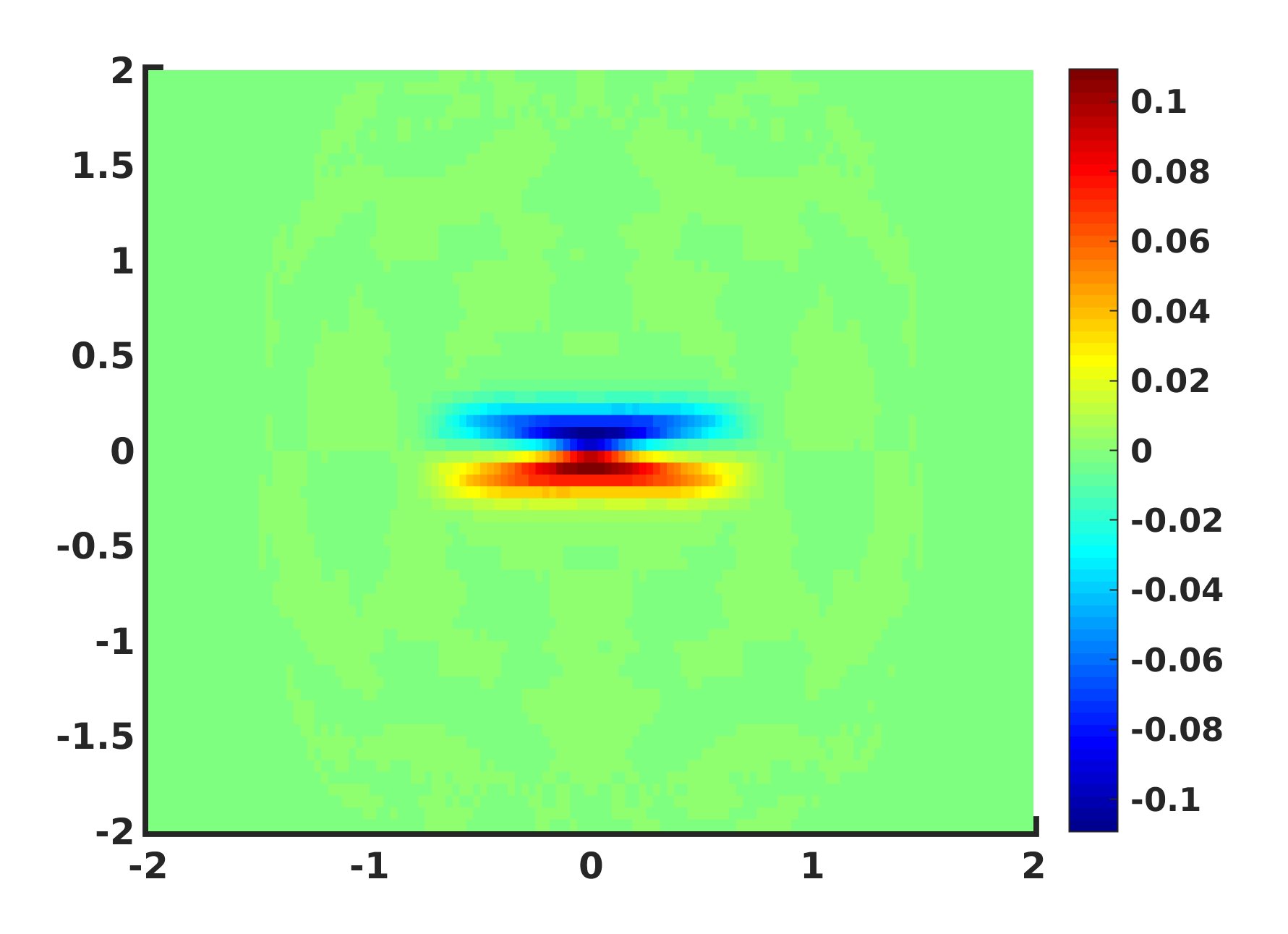} &
   \includegraphics[width=0.2\textwidth]{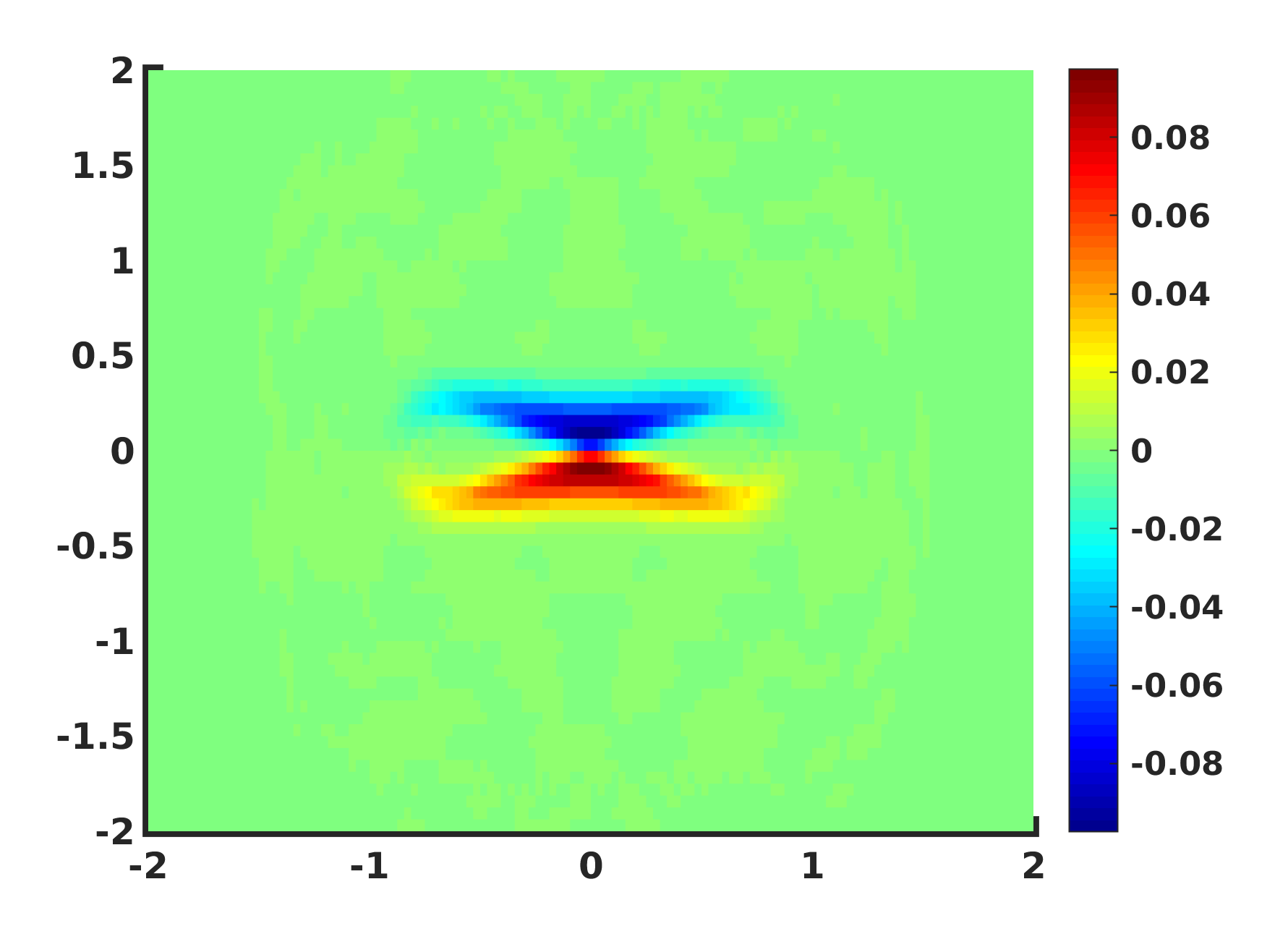} &
 \includegraphics[width=0.2\textwidth]{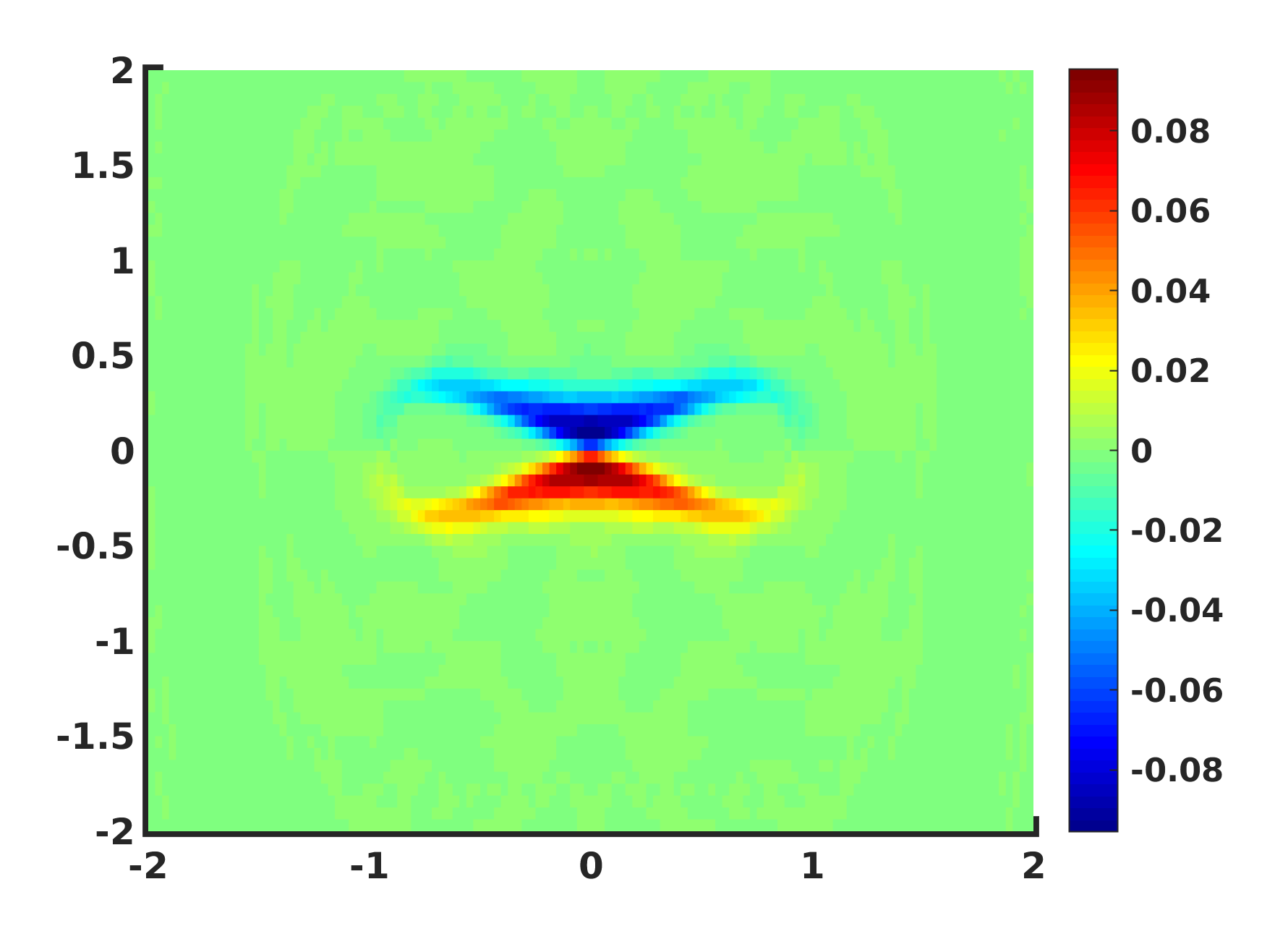} &
 \includegraphics[width=0.2\textwidth]{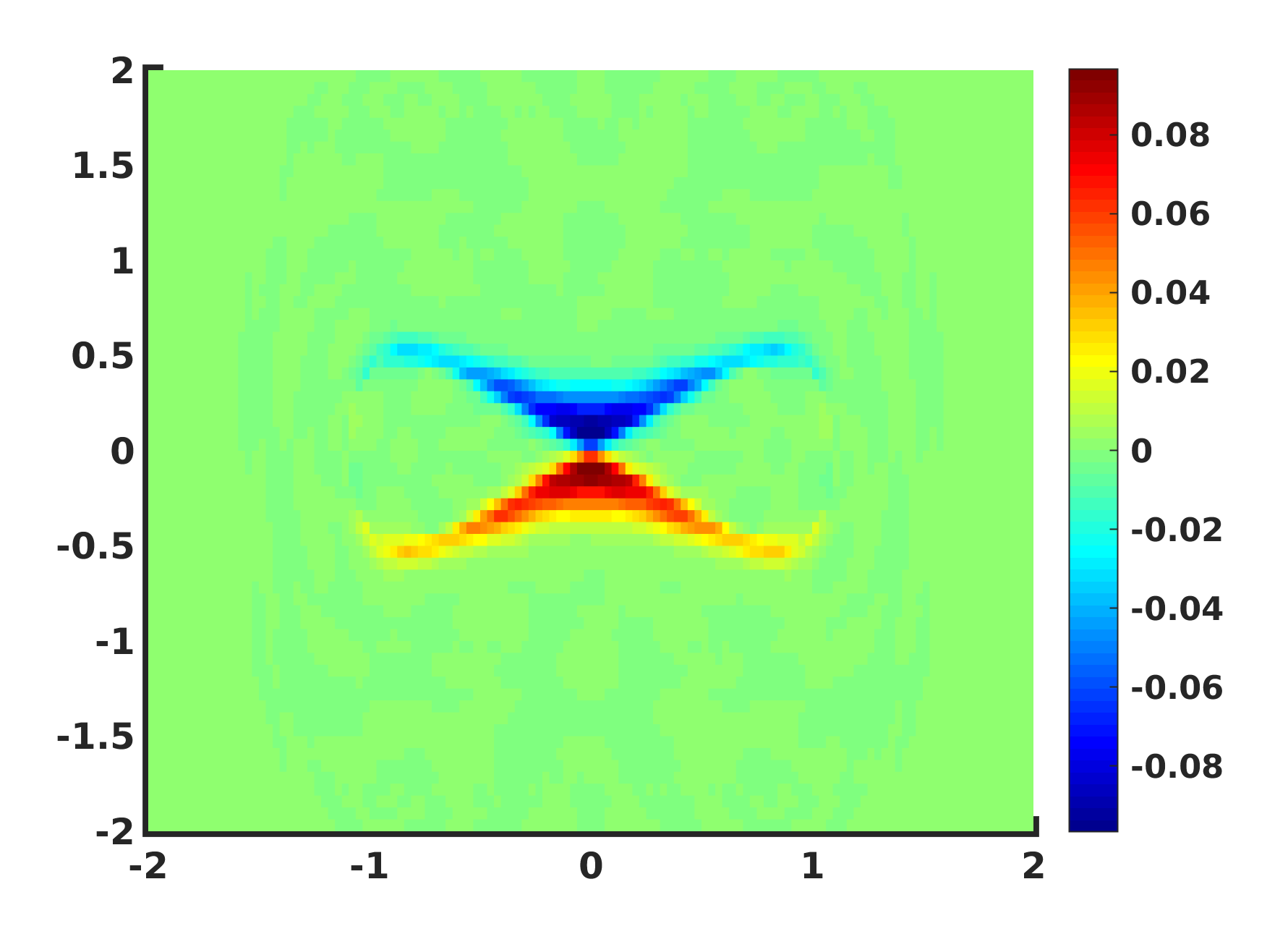}\\
  {\scriptsize (e) } &   {\scriptsize (f) } &   {\scriptsize (g) } &   {\scriptsize (h) } \\
\end{tabular}
\end{center}
\caption{\label{fig:pt2} Changes of the second "weak modes" patterns from $U_2^\prime(\epsilon,p)$ for increasing $p$.}
\end{figure}

The numerical approximation of the spectra is done with exactly the same settings as in section \ref{expModels}. However, the system \eqref{psys1} is open, which 
means that some test points will leave the domain of interest under the evolution of the flow map. To fix this issue, 
an additional box is added in order to capture all the image points that are being mapped out of the initial domain $M$ when computing the transition matrix. Finally, 
this temporary box will be removed from the eigenvector entries by just considering the $2^{depth}$ first entries.
\begin{figure}[!htb]
\begin{center}
\includegraphics[width=.6\textwidth]{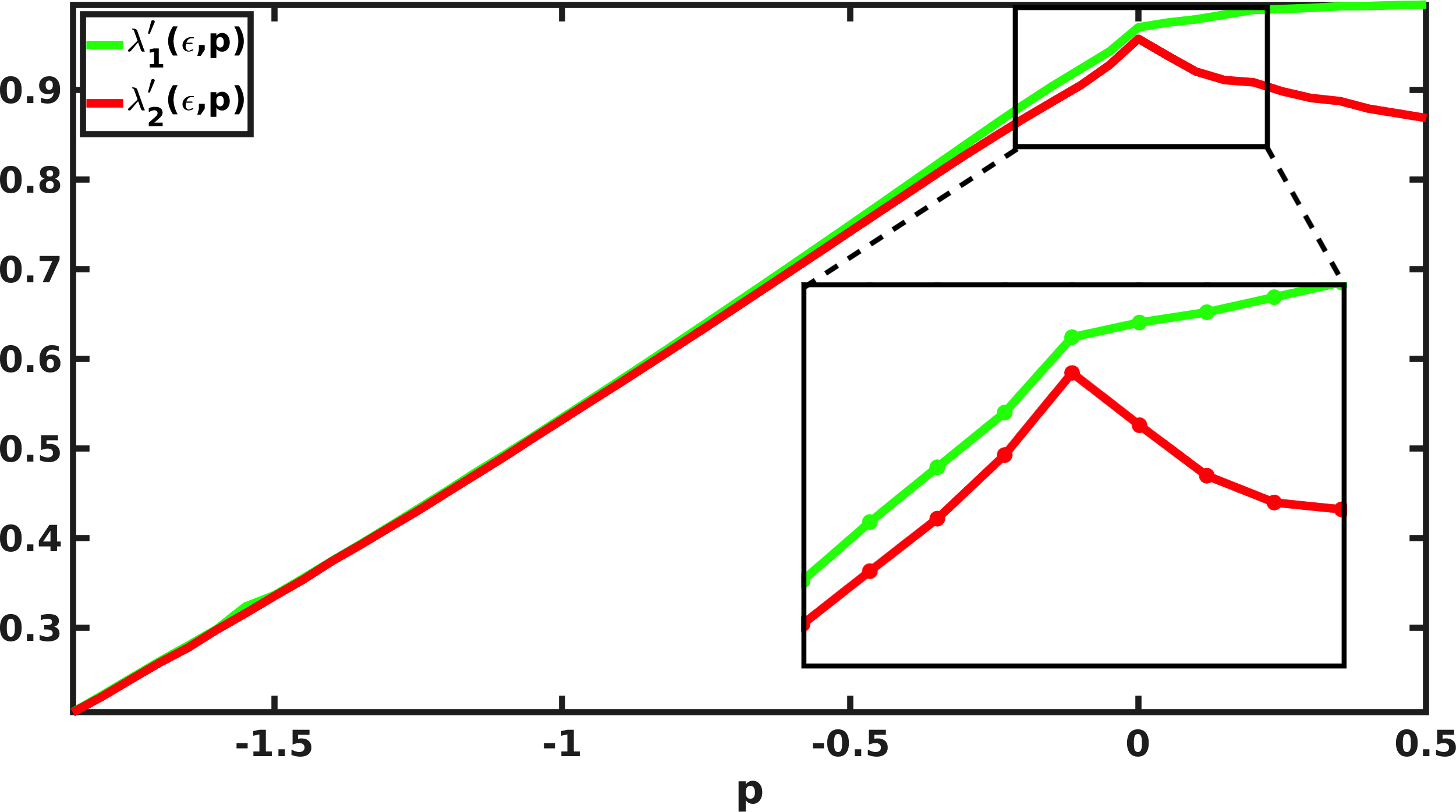} \\[-5mm]
\end{center}
\caption{\label{fig:pteig}  Spectral version of the classical bifurcation diagram with a zoomed diagram in the vicinity of the bifurcation (inlet). Two subdominant eigenvalues $\lambda_1^\prime(\epsilon,p)$ and $\lambda_2^\prime(\epsilon,p)$ rise towards one.}
\end{figure}
Figure \ref{fig:pteig} shows the changes of two small magnitude eigenvalues 
that belong to $\{\lambda_j(\epsilon,p),j=k+1,\ldots,N \}$. We denote by $\lambda_1^\prime(\epsilon,p)$ the green curve of eigenvalues with corresponding 
eigenvectors $U_1^\prime(\epsilon,p)$, in figure \ref{fig:pt1}. Likewise, $\lambda_2^\prime(\epsilon,p)$ corresponds to the red curve in figure \ref{fig:pteig}; their corresponding eigenvectors $U_2^\prime(\epsilon,p)$ are shown in Figure \ref{fig:pt2}.

In figures \ref{fig:pt1}-\ref{fig:pt2}, the patterns of $U_1^\prime(\epsilon,p)$ and $U_2^\prime(\epsilon,p)$ change in size, as the bifurcation parameter $p$ varies. Note that their 
sign structure remains the same. Besides, these eigenvector patterns are only supported on a small isolated neighborhood of $(0,0)$, according to the discretization depth 
of the phase space. 
In figure \ref{fig:pteig}, the eigenvalues $\lambda_1^\prime(\epsilon,p)$ and $\lambda_2^\prime(\epsilon,p)$ are initially very small compared to $1$. They increase linearly fast
together side-by-side until $p=0$. Then they part ways: $\lambda_1^\prime(\epsilon,p)$ continues to increase, while $\lambda_2^\prime(\epsilon,p)$ starts to decrease. That is why 
figure \ref{fig:pteig} is referred as the spectral version of the classical pitchfork bifurcation diagram of system \eqref{psys1} in analogy to the classical 
pitchfork bifurcation diagram (see e.g. \cite{Gheimer}, Chapter $3$, p. $146$). 

In figures \ref{fig:pt1}-\ref{fig:pt2}, one sees that the sign structure of $U_1^\prime(\epsilon,p)$ is symmetric with respect to the $y$-axis, 
while the sign structure of $U_2^\prime(\epsilon,p)$ is symmetric with respect to the $x$-axis. 
As $p<0$ increases towards zero, both eigenvector patterns expand slowly and symmetrically 
along the $x$-axis, but remain nearly constant in the $y$-direction. This is intrinsic to the underlying dynamical system, see figures \ref{fig:pt1}-\ref{fig:pt2} (a)-(c). 
At $p=0$, the fixed point $(0,0)$ bifurcates, which is particularly well observed in figures \ref{fig:pt1}-\ref{fig:pt2} (d). \\
\begin{figure}
\begin{center}
\begin{tabular}{ccc}
 \includegraphics[width=.3\textwidth]{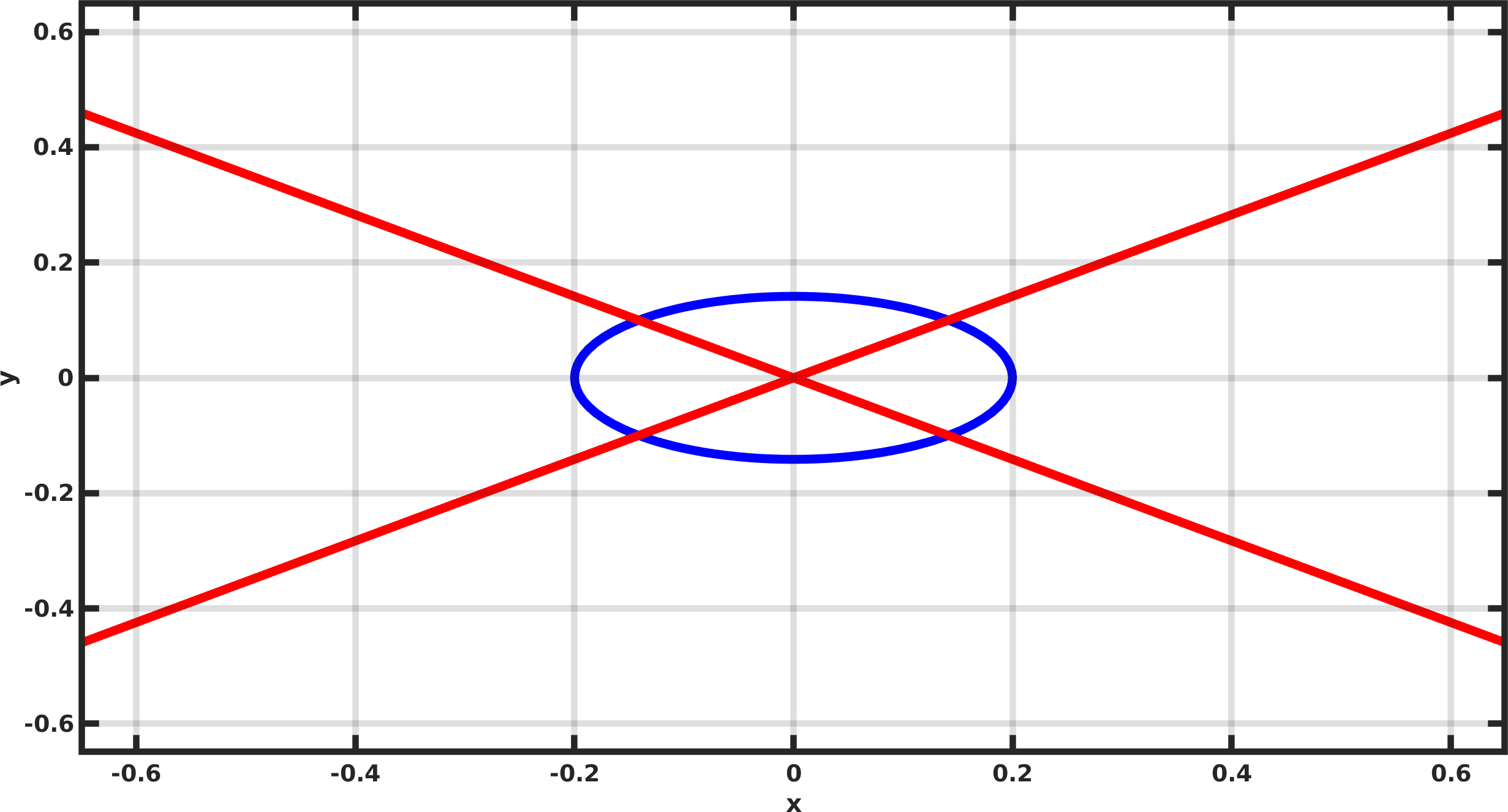} &
  \includegraphics[width=.3\textwidth]{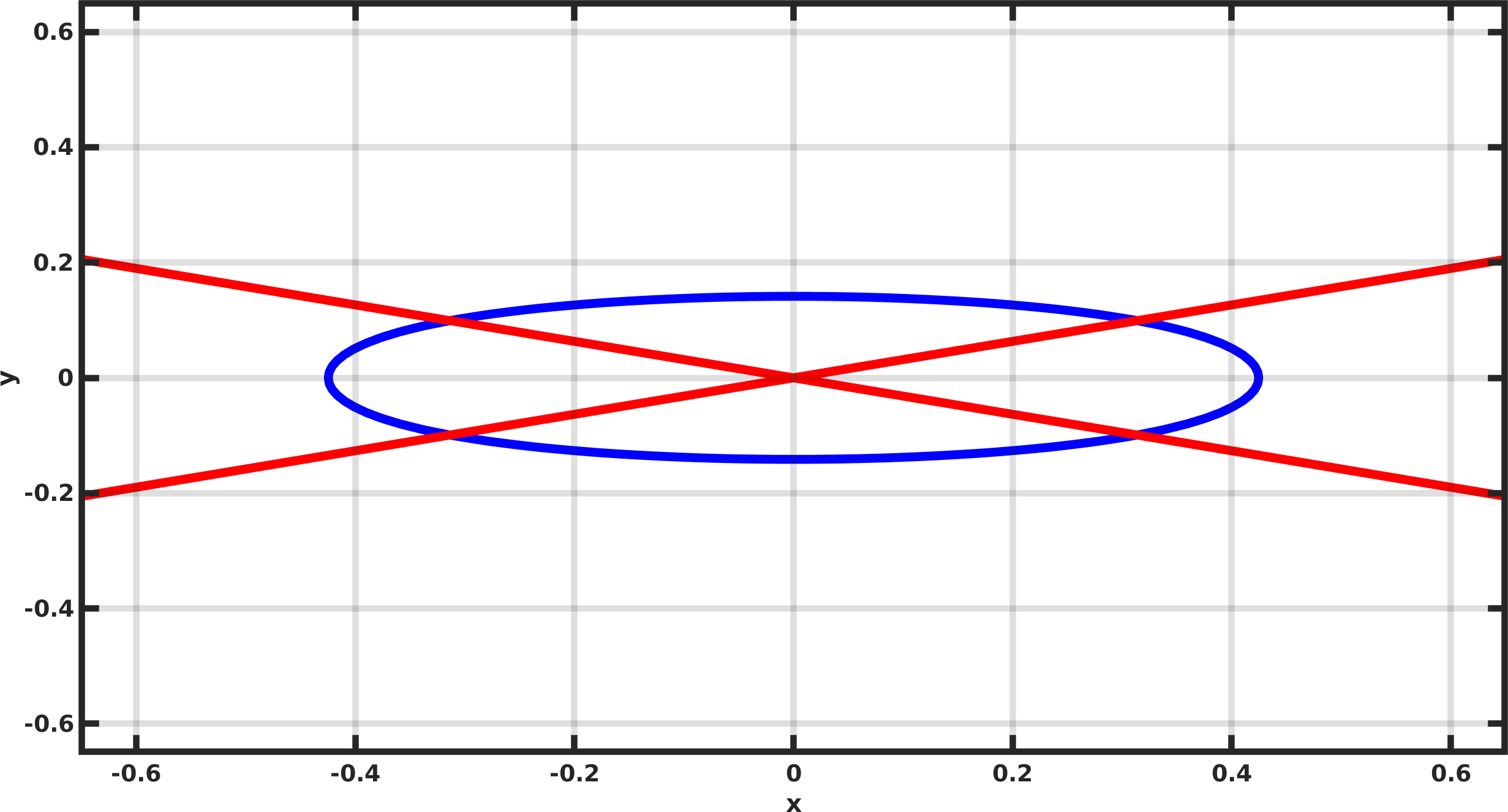} &
   \includegraphics[width=.3\textwidth]{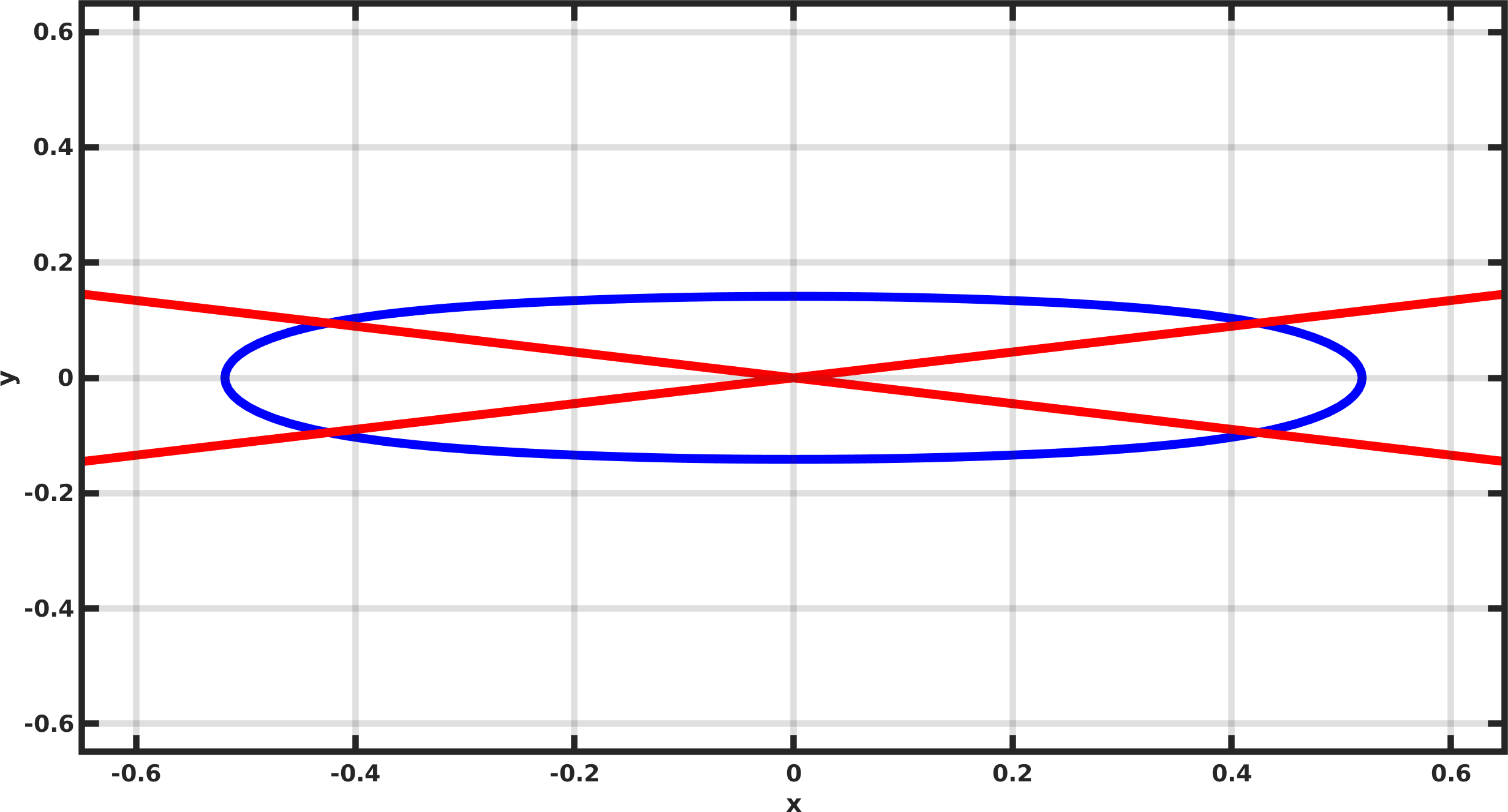} \\
    {\scriptsize (a) $p<0$ } &   {\scriptsize (b)  $p<0$} &   {\scriptsize (c)  $p<0$} \\
   \includegraphics[width=.3\textwidth]{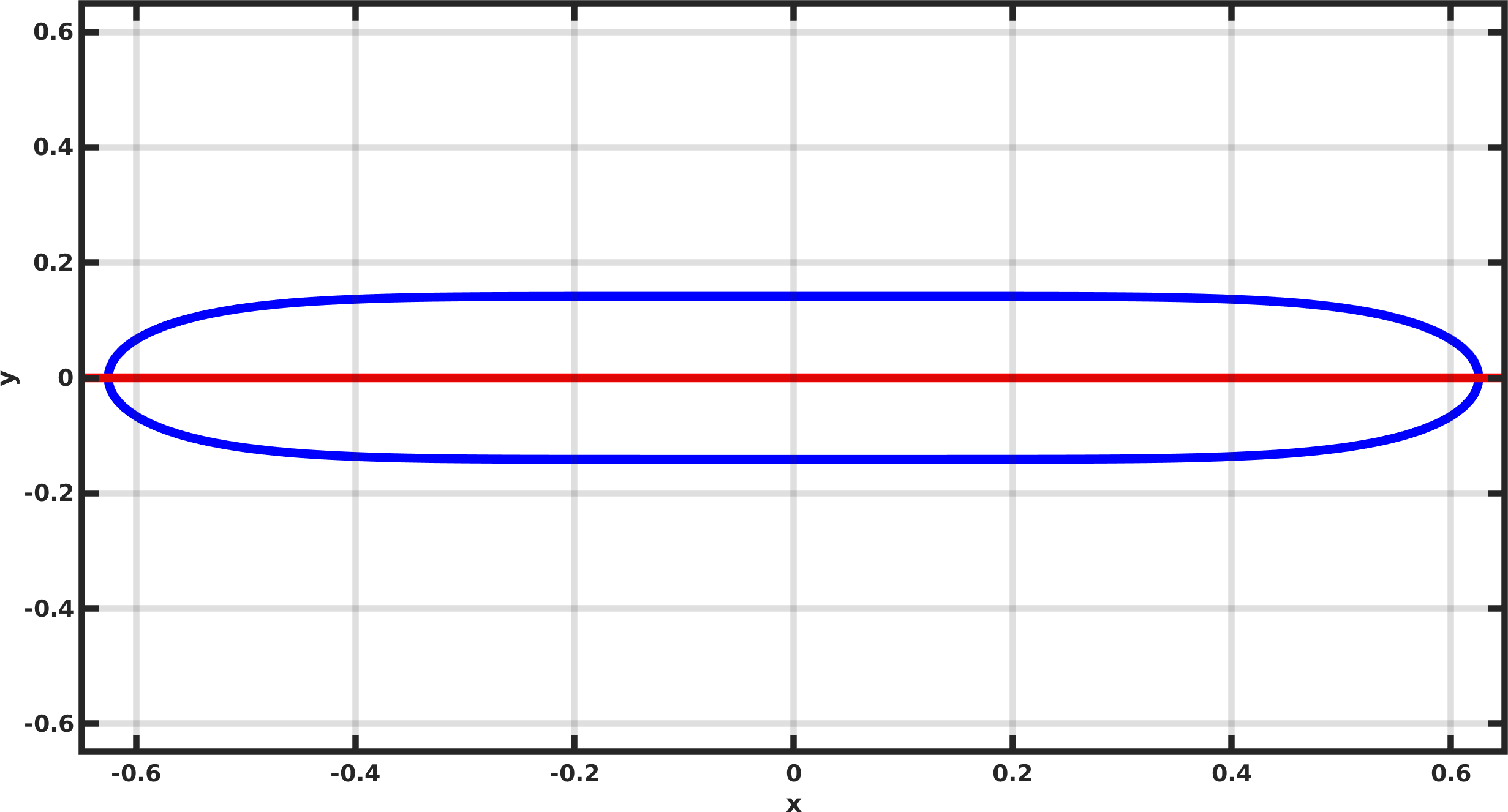} &
  \includegraphics[width=.3\textwidth]{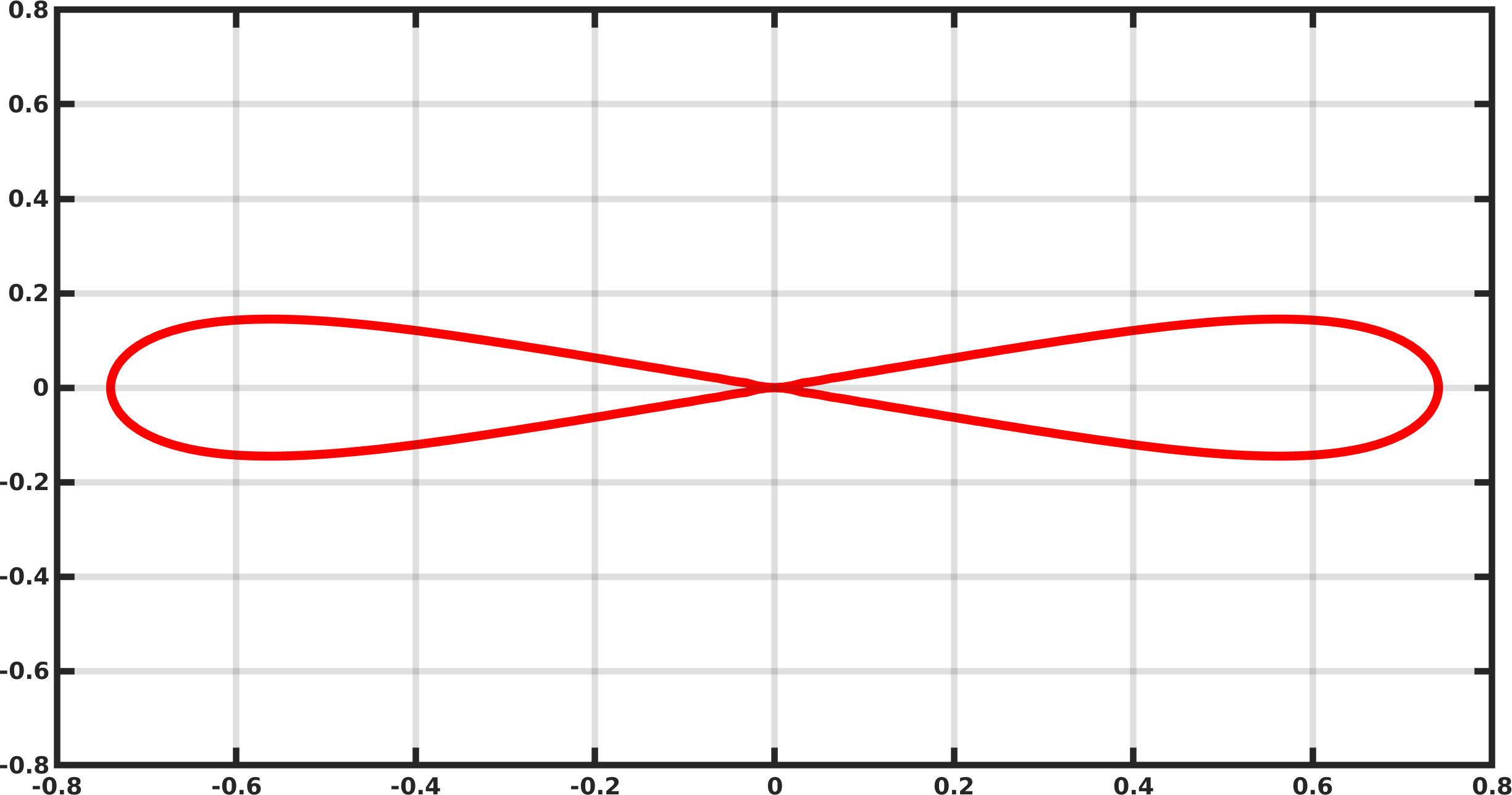} &
   \includegraphics[width=.3\textwidth]{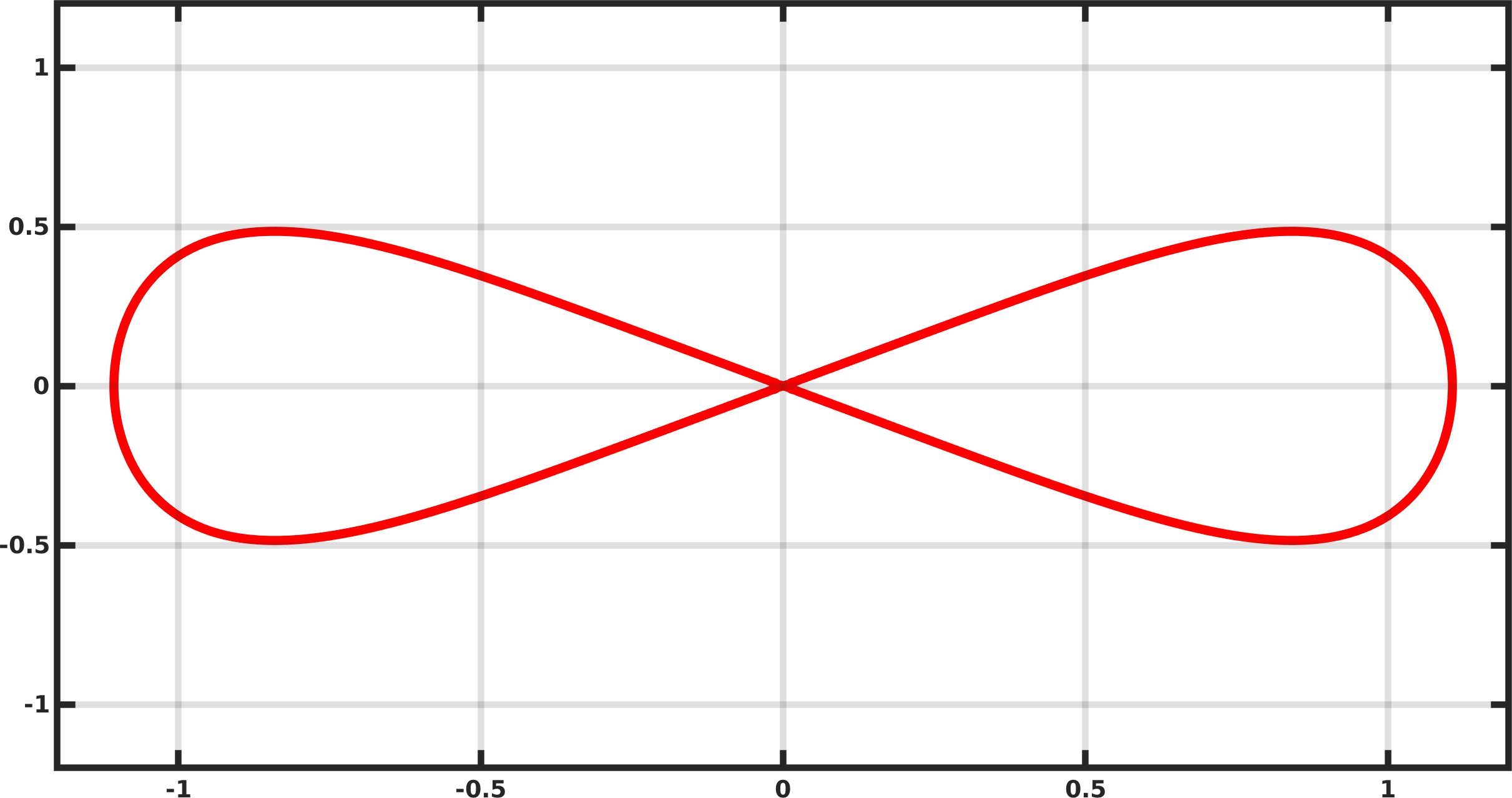} \\ 
     {\scriptsize (d) $p=0$ } &   {\scriptsize (e)  $p>0$} &   {\scriptsize (f)  $p>0$} \\
\end{tabular}
\end{center}
\caption{Local dynamics of system \eqref{psys1} near the origin. For $p\leq 0$ we obtain center manifolds, for $p>0$ homoclinic manifolds.}
\label{fig:cm}
\end{figure}
Another way to see the classical bifurcation is illustrated in figure \ref{fig:cm}. Indeed, the linearization of \eqref{psys1} around the fixed point $(0,0)$ yields a two-dimensional 
matrix whose eigenvalues are purely imaginary. With $\gamma_1 = i\sqrt{-p}$ and $\gamma_2 = -i\sqrt{-p}$ being the eigenvalues of the linearized system, the corresponding 
(generalized) eigenvectors are $v_1 = [x,\,y= \gamma_1 x]^\intercal$ and $v_1 = [x,\,y= \gamma_2 x]^\intercal$, respectively. Let $ E^c_1(p)$ and $E^c_2(p)$  be the (generalized) 
eigenspaces of $v_1$ and $v_2$, respectively, depicted as the two intersecting red lines in figure \ref{fig:cm}(a)-(c). 
Then $E^c(p) = E^c_1(p) \bigoplus E^c_2(p)$ yields the two-dimensional subspace spanned by $E^c_1(p)$ and $E^c_2(p)$. In figure \ref{fig:cm}(a)-(c), $E^c(p)$ is the plane generated by the 
intersection of the two red lines. Therefore, there exists an invariant manifold denoted 
by $W^c(p)$ that is tangent to $E^c(p)$ at $(0,0)$. This is known as the center manifold theorem (\cite{Gheimer}, Chapter $3$, p. $127$) and its main purpose is to isolate the 
complicated asymptotic behavior of the flow by locating such an invariant manifold $W^c$. In system \eqref{psys1}, for negative $p$, every closed 
orbit is a boundary of a center manifold which is tangent to $E^c(p)$ at $(0,0)$. 

Notably, closed orbits are given by $q$-level sets of the derived
Hamiltonian functions $H_q(x,y,p)$ for different $p$. Therefore, there is a constant $\bar{q}$ small enough 
such that $\{H_{\bar{q}}(x,y,p) = \bar{q}\}$ isolates the asymptotic dynamics of the flow near $(0,0)$ from the rest. The blue closed curve shows 
$\{H_{\bar{q}}(x,y,p) = \bar{q}\}$ in figure \ref{fig:cm}(a)-(d). For negative $p$, the eigenspace plane $E^c(p)$ partitions the interior of the closed curve 
$\{H_{\bar{q}}(x,y,p) = \bar{q}\}$ into four regions which are two-by-two symmetric similar to $U_1^\prime(\epsilon,p)$ and $U_2^\prime(\epsilon,p)$ in 
figure \ref{fig:pt1} and figure \ref{fig:pt2}, respectively. As $p$ increases towards $0$, the slopes and the intersection angles of $E^c_1(p) $ and $E^c_2(p)$ decrease and the closed 
curve expands horizontally, while 
remaining constant vertically. Again, this is analogous to the variations of $U_1^\prime(\epsilon,p)$ and $U_2^\prime(\epsilon,p)$, for negative $p$.  
In this way, the local behavior of the flow near the origin is cast into $U_1^\prime(\epsilon,p)$ and $U_2^\prime(\epsilon,p)$ 
for negative $p$. This is another way of providing a better understanding of complicated asymptotic dynamics near $(0,0)$, from a probabilistic approach. 
Note that the choice of $\bar{q}$ is heuristic since there no way to have an exact $\bar{q}$ in order to have the exact local manifold that supports the corresponding eigenvectors. 

For $p=0$,  $E_1^c(p)$ and $E^c_2(p)$ disappear as a consequence of the classical bifurcation. Note that figure \ref{fig:cm}(d) is also in agreement with 
figure \ref{fig:pt1}(d) and figure \ref{fig:pt2}(d).
For $p$ positive, there are three fixed points: Two elliptic fixed points  $(- \sqrt[\leftroot{-2}\uproot{2}4]{p}, 0)$ and  $(\sqrt[\leftroot{-2}\uproot{2}4]{p}, 0)$ and one saddle 
fixed point $(0,0)$. 

The qualitative behavior of the dynamics changes radically with the emergence of two symmetric homoclinic orbits, as shown in \ref{fig:cm}(e). We refer 
to the latter as the homoclinic manifold, since neighboring trajectories are periodic and tangent to it. As $p$ increases, the homoclinic manifold 
increases in size (see figure \ref{fig:cm}(e)-(f)) because $g(p):=\sqrt[\leftroot{-2}\uproot{2}4]{p} $ is an increasing function of $p$.
Nearby solution curves tend to be attracted through the $y$-direction and repelled through the $x$-direction. Indeed, 
figure \ref{fig:cm}(e)-(f) shows that the homoclinic manifold is concave in the $y$-direction and convex in the $x$-direction. Besides, the larger $p$ gets, the more 
does the curvature of the homoclinic manifold grow.

This implies immediately that the global behavior of the dynamics becomes attractive along the $y$-direction. 
Thus, the support of $U_2^\prime(\epsilon,p)$ shrinks symmetrically on both sides of $x$-axis, as shown in figure \ref{fig:pt2}(e)-(h) and, as a matter of fact, 
the eigenvalues $\lambda_2^\prime(\epsilon,p)$ decrease. On the other hand, the dynamics repels along the $x$-direction. Thus, the support of $U_1^\prime(\epsilon,p)$ 
expands symmetrically on both sides of the $y$-axis and is simultaneously enrolled into the two newly co-existing homoclinic orbits; see figure \ref{fig:phaseplot}. 
The latter expand as $p$ increases from zero. As a consequence, the eigenvectors $U_1^\prime(\epsilon,p)$ carry almost-invariant sets bounded by the homoclinic orbits, 
for $p>0$.  Moreover, the corresponding eigenvalues $\lambda_1^\prime(\epsilon,p)$ increase towards $1$, see figure \ref{fig:pteig}.

Recall that these eigenvalues belong to the set of $N-k$ small magnitude spectrum. Hence, the continuous rise of $\lambda_1^\prime(\epsilon,p)$ towards $1$, 
as a consequence of the classical bifurcation, will eventually question the well-definiteness of $k$ dominant 
eigenvectors and their corresponding eigenvalues. Indeed, \eqref{ppmueigvec} is no longer valid if the additional inequality constraint fails.
In the next subsection, this will play a key role for us to characterizing bifurcations of almost-invariant sets.
\subsection{Predicting bifurcation of almost-invariant patterns}
Here, we will characterize bifurcations of almost-invariant patterns and deduce the corresponding generic early warning signs. As mentioned before, we are interested in the changes of the particular 
pattern centered in $(0,0)$ and located in each dominant eigenvector pattern. As shown in figure \ref{fig:pEV}, given any dominant eigenvector $U_j(\epsilon,p),\,j=2,\ldots,k$, 
the particular almost-invariant pattern, denoted $\mathcal{P}_j,\,j=2,\ldots,k$, is the $j^\text{th}$ partition element surrounded  by all ring patterns.  
Indeed, every $U_j(\epsilon,p),\,j=2,\ldots,k$ yields $j$ almost-invariant patterns which partition the state space. In particular, for $j=2$, the second dominant 
eigenvector yields two almost-invariant patterns partitioning the state space. Besides, one of the patterns yields $\mathcal{P}_2$, which is known as the maximal almost-invariant set \cite{GFnonauto2}; 
see figure \ref{fig:pEV}(b). Moreover, $(U_2(\epsilon,p),\mathcal{P}_2)$ is usually a good candidate for modelling real world isolated patterns such as atmospheric vortices. 

\begin{figure}[!htb]
\begin{center}
\begin{tabular}{cccc}
 \includegraphics[width=0.2\textwidth]{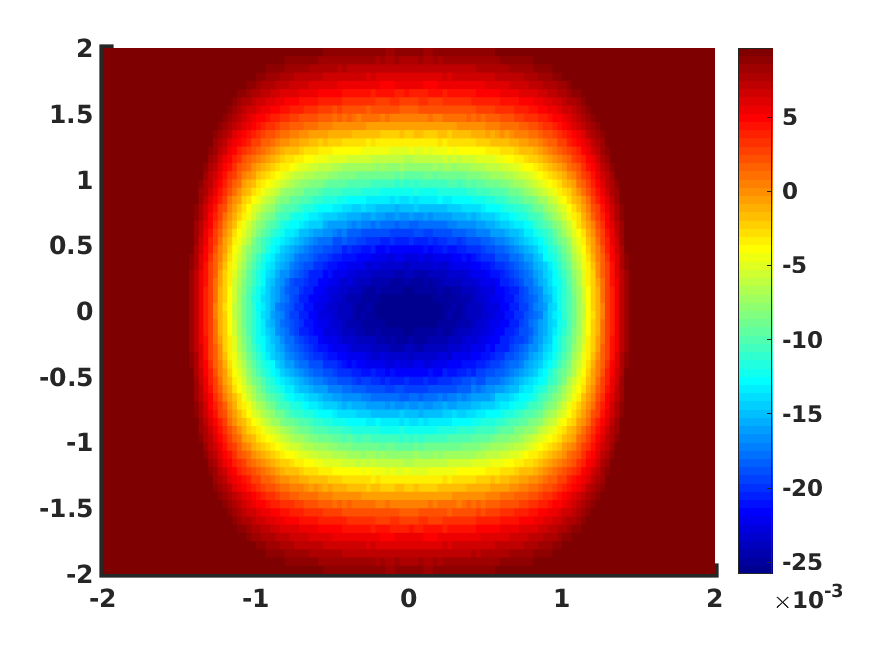} &
  \includegraphics[width=0.2\textwidth]{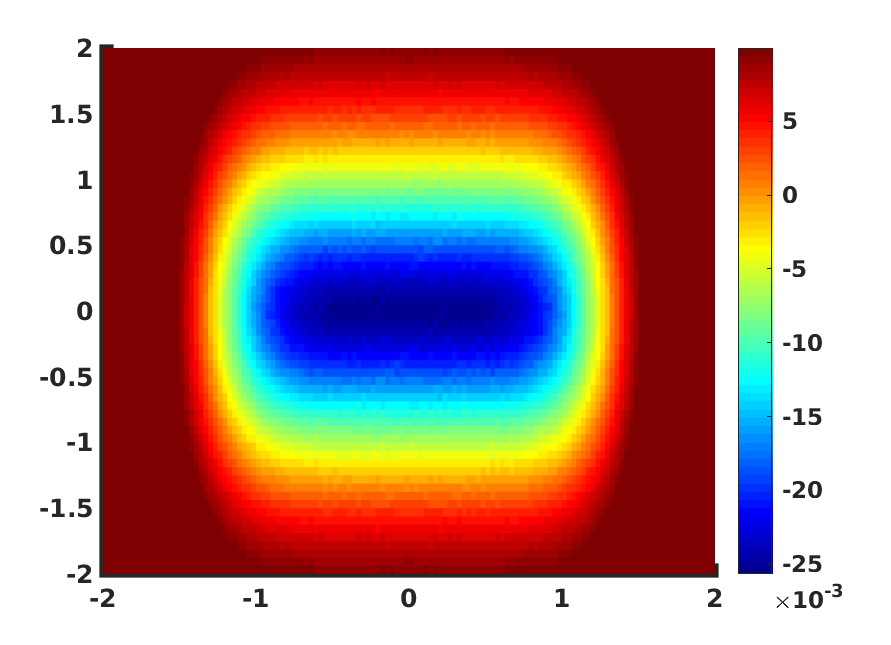} &
 \includegraphics[width=0.2\textwidth]{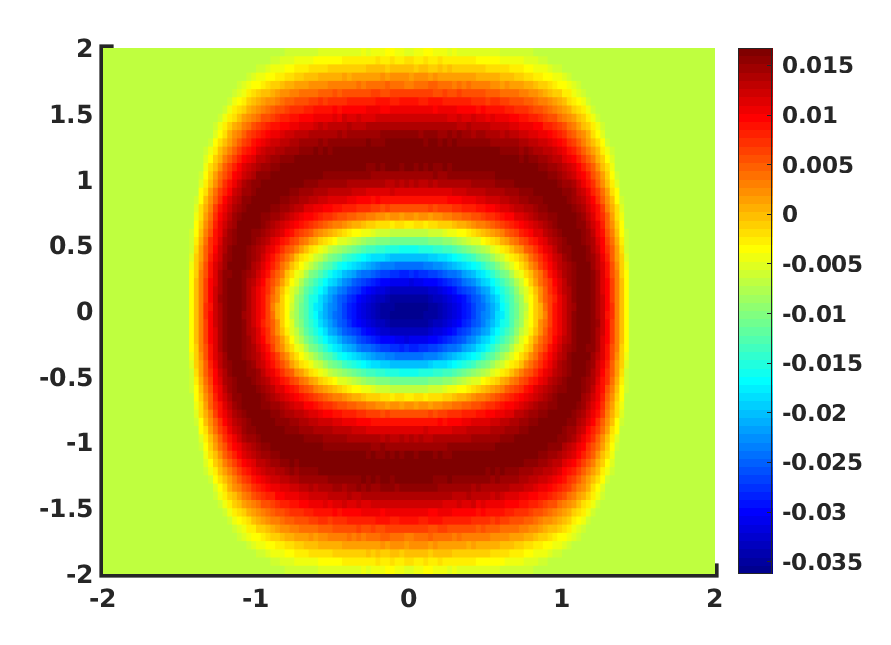} & 
  \includegraphics[width=0.2\textwidth]{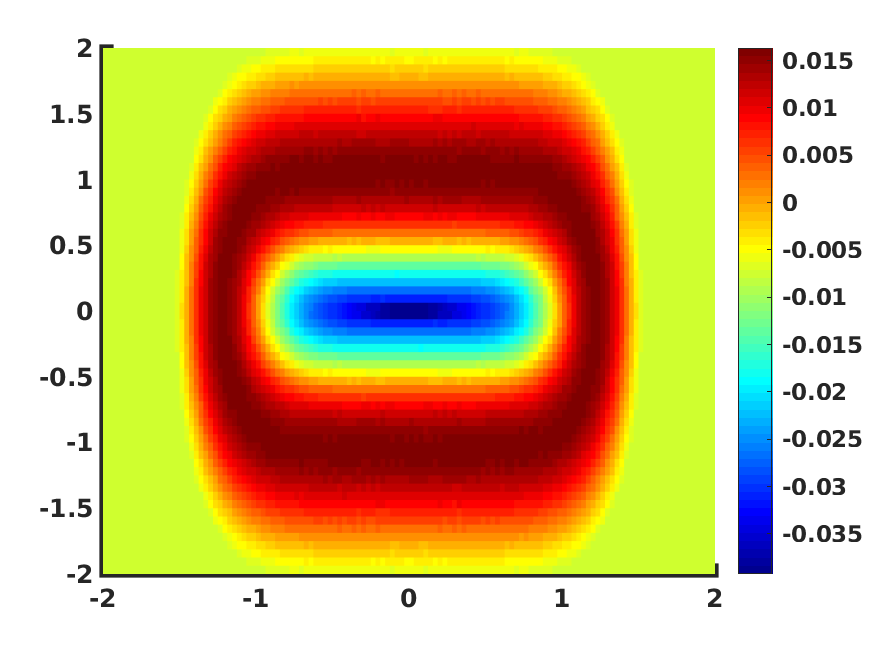}\\
  {\scriptsize $U_2(\epsilon, p<0)$} &   {\scriptsize $U_2(\epsilon, p=0)$} & {\scriptsize $U_3(\epsilon, p<0)$} &  {\scriptsize $U_3(\epsilon, p=0)$}\\
\end{tabular}
\end{center}
\caption{\label{fig:evatcriticalyone} $\,U_j(\epsilon,p\leq0),\,\,j=2,3$.}
\end{figure}

\begin{figure}[!htb]
  \centering
    \includegraphics[width=0.5\textwidth]{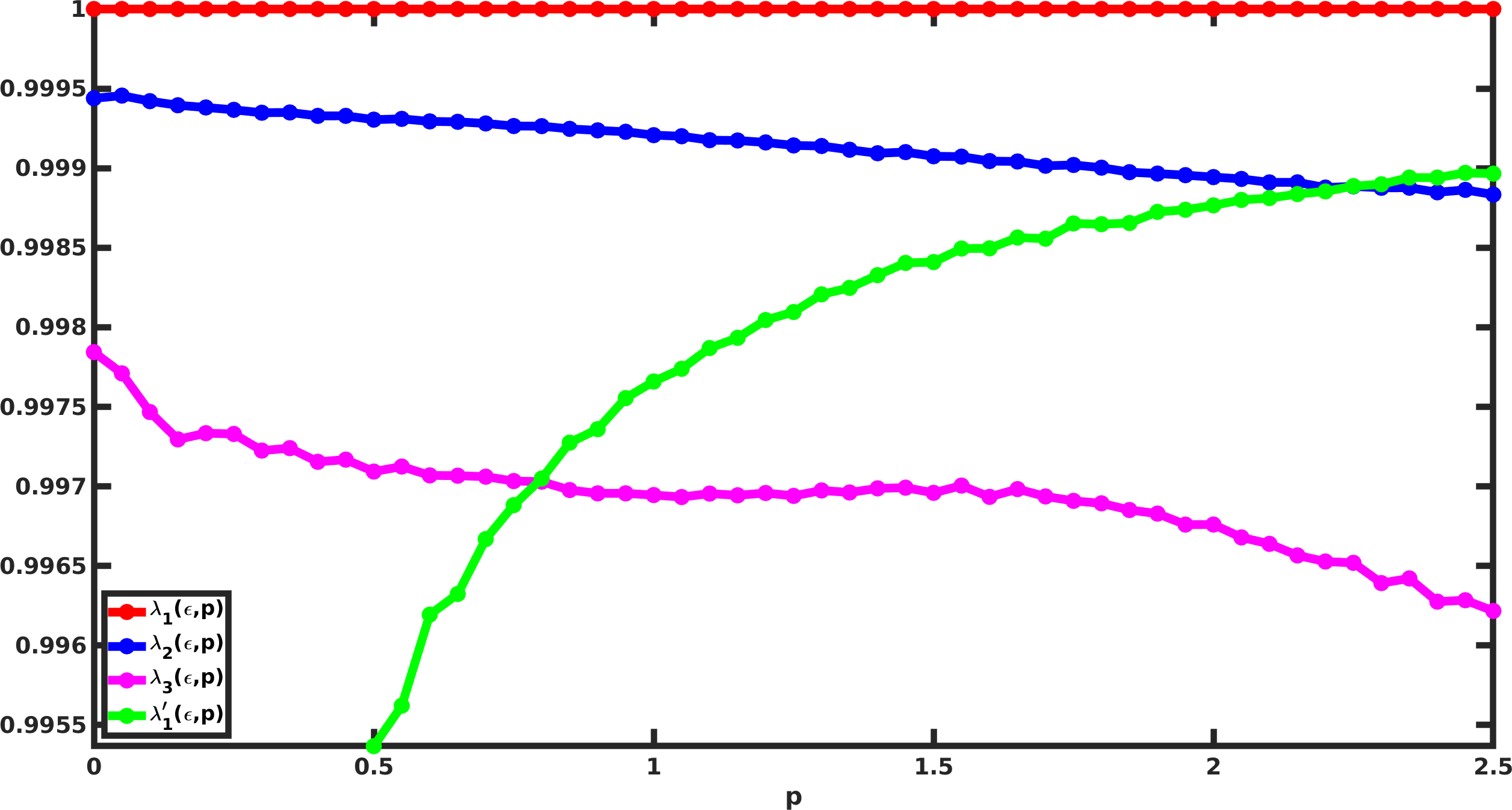}
\caption{\label{fig:alleigs} Spectral indicators of bifurcations of almost-invariant patterns in system  \eqref{psys1}. 
Three dominant eigenvalues $\lambda_1(\epsilon,p)=1$ (red),  $\lambda_3(\epsilon,p)<\lambda_2(\epsilon,p)<1$ (magenta, blue) and the rising eigenvalue $\lambda_1^\prime(\epsilon,p)$ (green) which intersects first the $\lambda_3(\epsilon,p)$-curve and then the $\lambda_2(\epsilon,p)$-curve.}
\end{figure}

Given the bifurcation diagram in figure \ref{fig:pteig} and the inequality constraint \eqref{degeneracy_eq}, we set the relation
\begin{equation}\label{pdegeneracy}
 \lambda_1^\prime(\epsilon,p) < \lambda_j(\epsilon,p) ,\,\, j=2,\ldots,k,
\end{equation}
which is true whenever $p\leq 0$. Moreover, for $p\leq 0$, the dominant patterns are stable in the sense that there is no qualitative 
change in their sign structure. In figure \ref{fig:evatcriticalyone}, we plot eigenvector patterns for $p\leq 0$. Note that the leading eigenvector is constant and therefore not shown. 
Even though there is a classical bifurcation in \eqref{psys1} 
at $p=0$, one notices that the dominant eigenvector patterns $U_2(\epsilon,p)$ and $U_3(\epsilon,p)$ in figure \ref{fig:evatcriticalyone} are qualitatively the same for $p<0$ and $p=0$. As a matter of fact, it suffices to investigate the bifurcation of almost-invariant patterns for positive values of $p$, as it can only happen in that parameter range. Indeed, when $p$ becomes positive, $\lambda_1^\prime(\epsilon,p)$ continues to increase monotonically to eventually become the second dominant 
eigenvalue after the eigenvalue $1$. In fact, $\lambda_1^\prime(\epsilon,p)$ will cross, in cascade, all the $k-1$ nontrivial dominant eigenvalues, as illustrated in 
 figure \ref{fig:alleigs}.

\begin{figure}[!htb]
\begin{center}
\begin{tabular}{cccc}
  \includegraphics[width=0.2\textwidth]{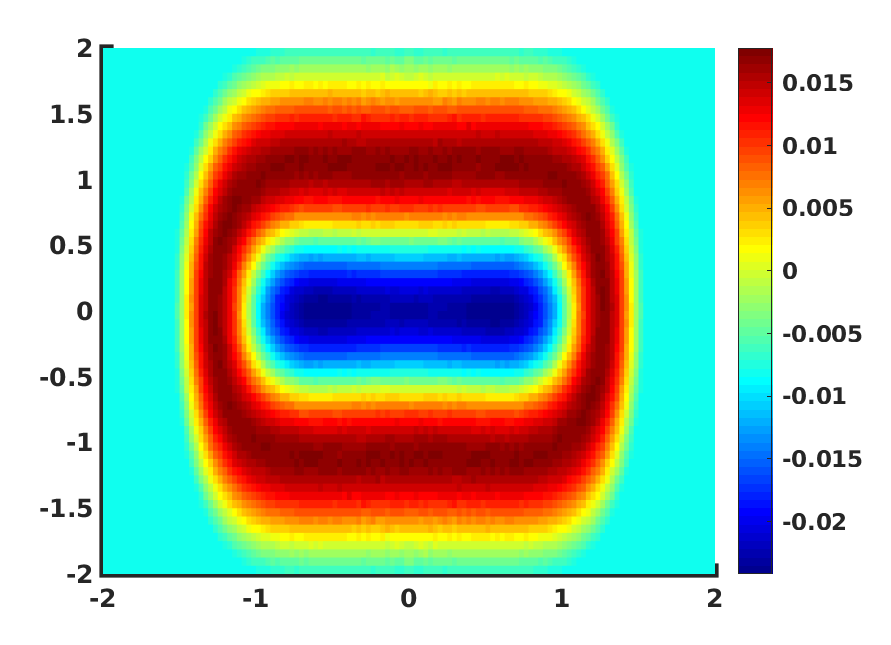} &
 \includegraphics[width=0.2\textwidth]{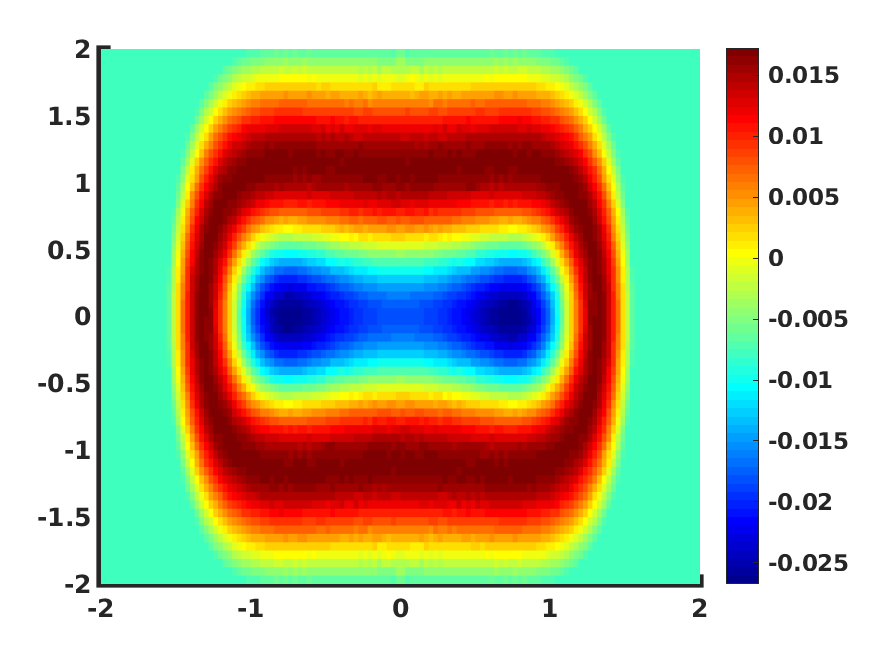} &
 \includegraphics[width=0.2\textwidth]{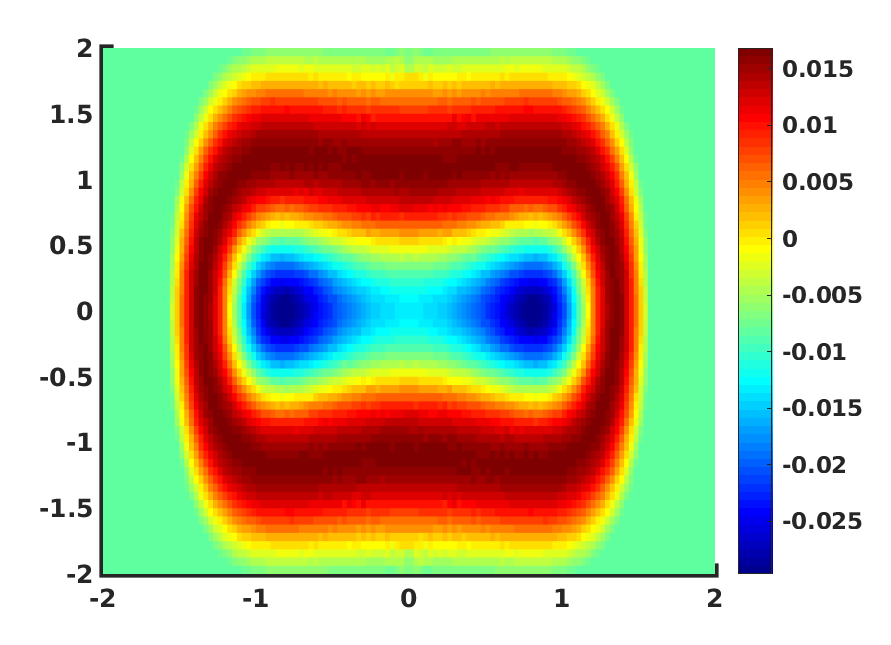} & 
  \includegraphics[width=0.2\textwidth]{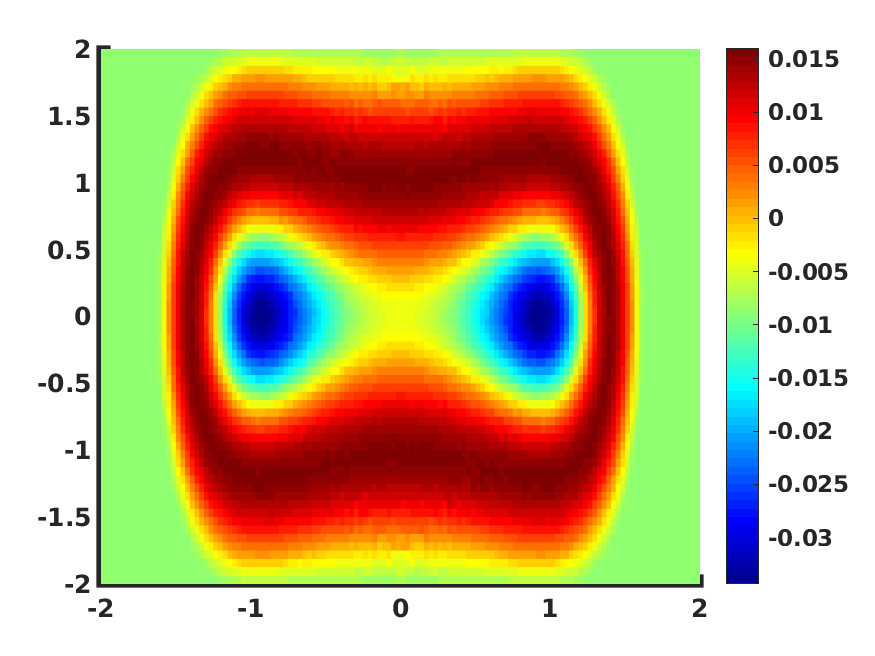}\\
  \end{tabular}
\end{center}
\caption{\label{fig:u3bifurc} Splitting process of the pattern $\mathcal{P}_3$ for increasing $p\geq 0$.}
\end{figure}

\begin{figure}[h!]
\begin{center}
 \includegraphics[width=0.48\linewidth]{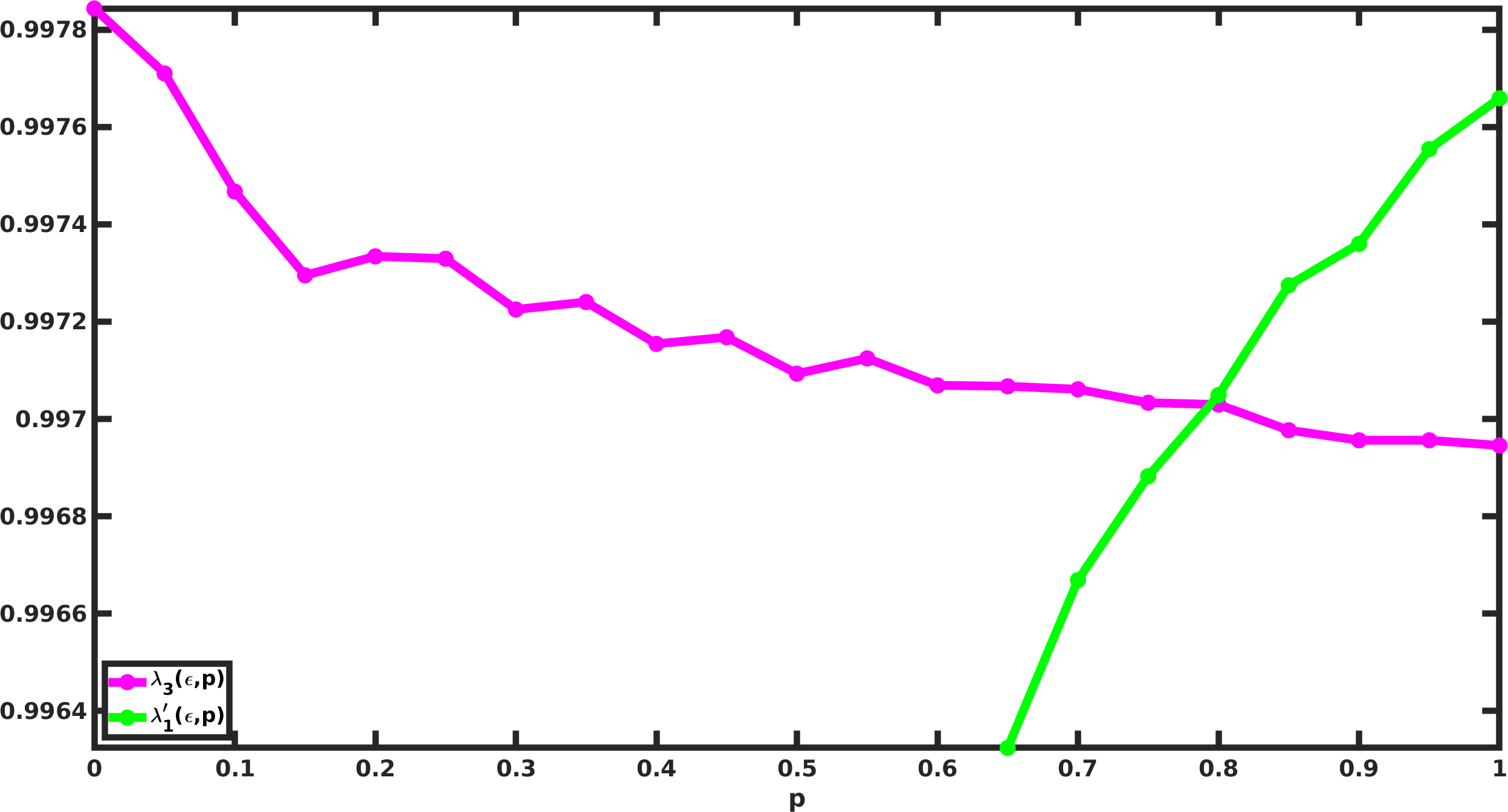}
 \includegraphics[width=0.48\linewidth]{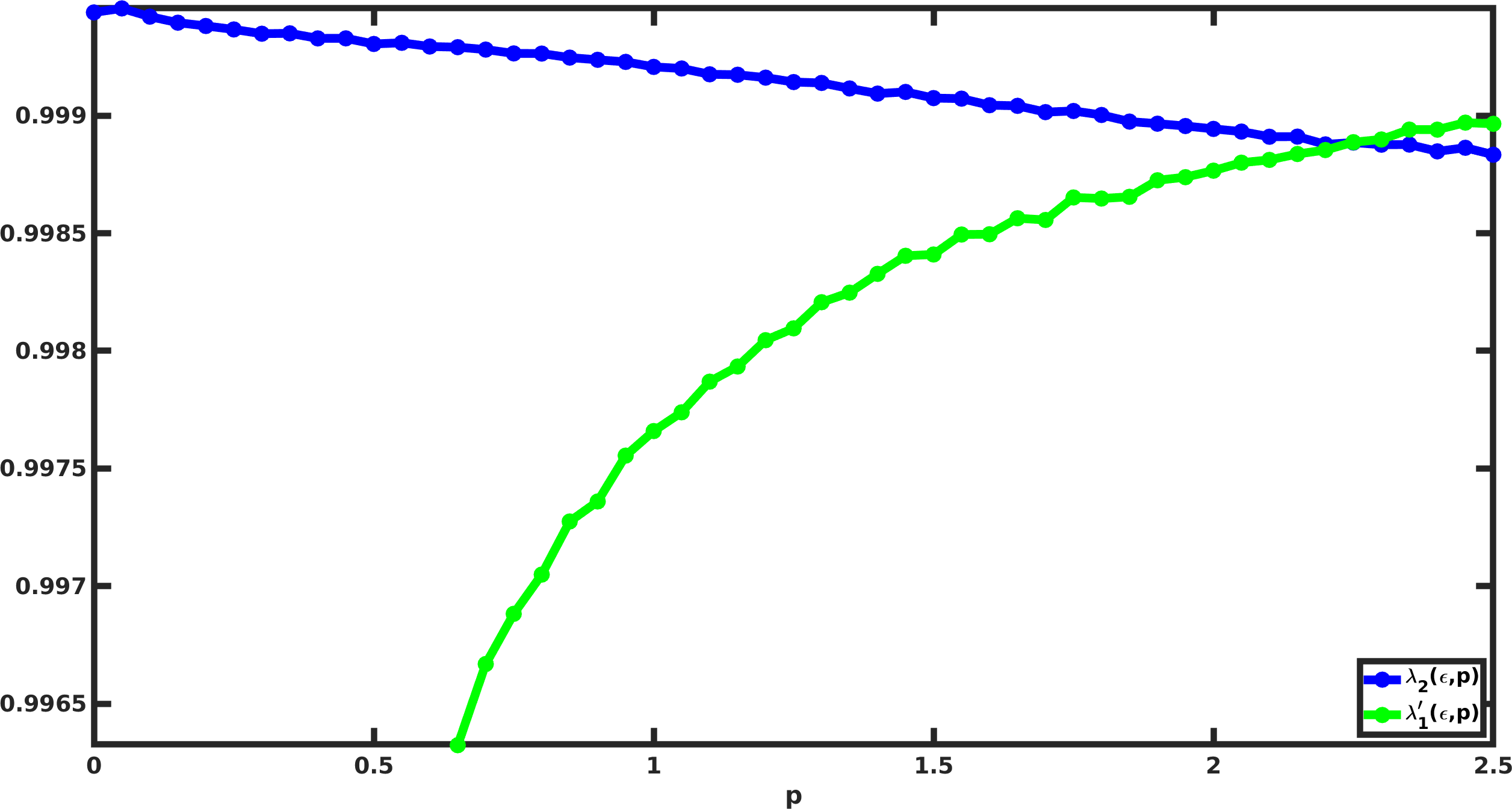}\\
 \end{center}
  \caption{Crossings of previously dominant eigenvalues when patterns $\mathcal{P}_3$ and $\mathcal{P}_2$ bifurcate under variation of $p$. Left: $\lambda_3(\epsilon,p)$ (magenta) vs. $\lambda_1^\prime(\epsilon,p)$ (green); right: $\lambda_2(\epsilon,p)$ (blue) vs.  $\lambda_1^\prime(\epsilon,p)$ (green). }\label{fig:onetonone_bif}
 \end{figure}

\begin{figure}[!htb]
\begin{center}
\begin{tabular}{cccc}
 \includegraphics[width=0.2\textwidth]{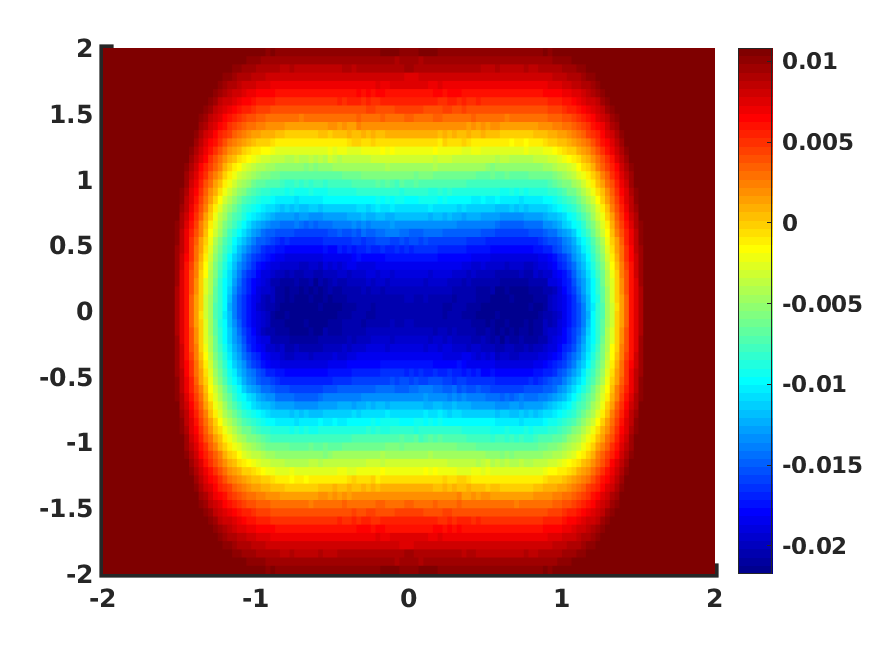} &
  \includegraphics[width=0.2\textwidth]{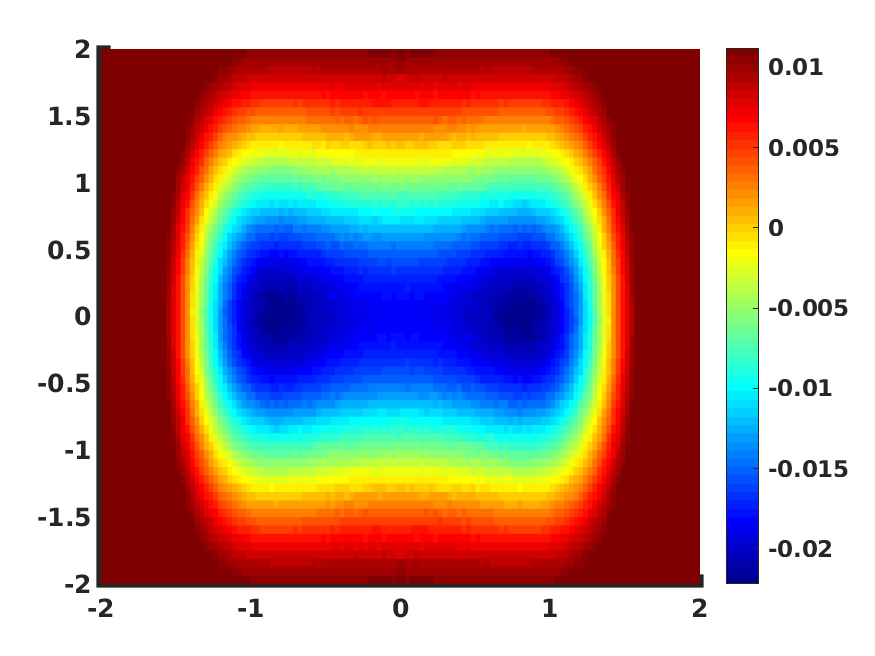} &
   \includegraphics[width=0.2\textwidth]{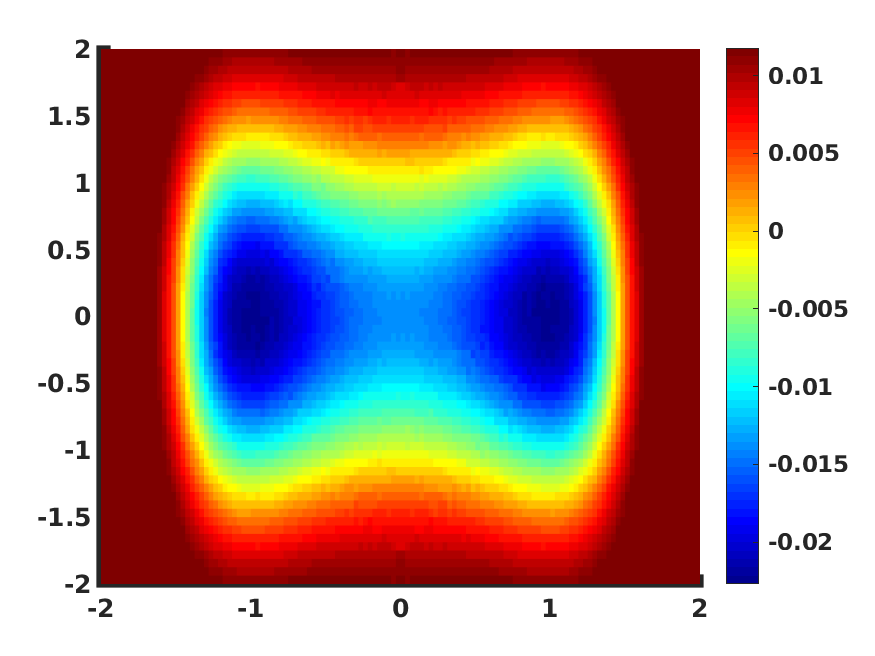} &
 \includegraphics[width=0.2\textwidth]{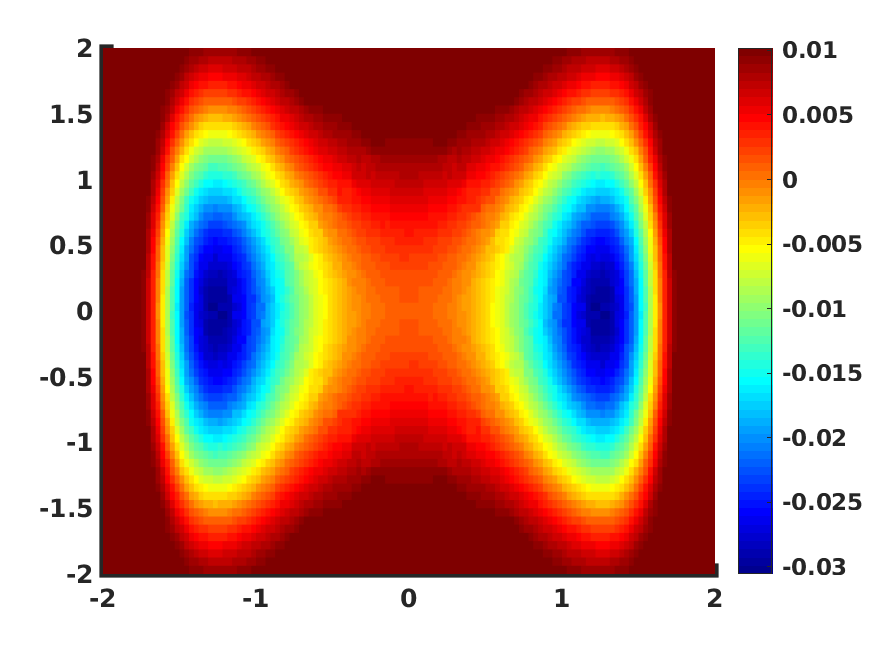}\\
 \end{tabular}
 \end{center}
\caption{\label{fig:u2bifurc} Splitting process of the pattern $\mathcal{P}_2$ for increasing $p\geq 0$.}
\end{figure}

For $p>0$, the global dynamics attracts along the $y$-axis 
and repels along the $x$-axis. As a consequence, the pattern generated by $U_1^\prime(\epsilon,p)$
expands in size, since it is supported in the region bounded by the two homoclinic orbits. On the other hand, the attractivity 
through the $y$-axis causes a shrinking process of the $k-1$ dominant eigenvectors patterns $U_j(\epsilon,p),\,j=2,\ldots,k$. 
Thus, as explained in the toy model experiments \ref{firstExample}-\ref{thirdExample}, we can define a set-oriented version of a degeneracy as  
\begin{equation}\label{bifurcation}
 \lambda_1^\prime(\epsilon,p) = \lambda_j(\epsilon,p) ,\,\, j=2,\ldots,k,\,\,\forall\,\,p>0.
\end{equation}
It follows that $\mathcal{P}_j,\,\,j=2,\ldots,k,$ bifurcates in the sense of a splitting, whenever equation \eqref{bifurcation} holds. Thus, according to figure \ref{fig:alleigs}, 
there is a cascade of two bifurcations. Every bifurcation occurs at a parameter $p>0$ where the support of $U_1^\prime(\epsilon,p)$ expands far enough to erupt out of 
$\mathcal{P}_j,\,\,j=2,3$. Indeed, at $p=0$ the support of $U_1^\prime(\epsilon,p)$ is, a priori, contained in the support of each $\,U_j(\epsilon,p),\,\,j=2,3$, specifically
in $\mathcal{P}_j,\,\,j=2,3$. 
This scenario changes radically the sign structures of  $U_j(\epsilon,p),\,\,j=2,3$ and, hence, the latter can no longer be 
expressed as in \eqref{peigvec}. The crossings occur in cascade from the smallest dominant eigenvalue to the biggest eigenvalue. In figure \ref{fig:onetonone_bif} (left), one can see 
that $\lambda_3(\epsilon,p)$ is crossed first. Later figure \ref{fig:onetonone_bif} (right) shows the last crossing scenario where 
$\lambda_1^\prime(\epsilon,p) = \lambda_2(\epsilon,p)$ after which $\lambda_1^\prime(\epsilon,p)$ becomes the dominant eigenvalue after the eigenvalue $1$. The 
eigenvector patterns from $U_2(\epsilon,p)$ and $U_3(\epsilon,p)$ undergo, respectively, a splitting process of $\mathcal{P}_3$ in figure \ref{fig:u3bifurc} and of
$\mathcal{P}_2$ in figure \ref{fig:u2bifurc}. Note that the splitting of the patterns $\mathcal{P}_j,\,\,j=2,3$ does not occur suddenly but gradually. In fact, one observes
the decreasing process of $\lambda_j(\epsilon,p),\,\,j=2,3$ before the crossing, which could be classified as an early warning signal.
 
\begin{figure}[!htb]
\begin{center}
\begin{tabular}{cccc}
\includegraphics[width=0.2\textwidth]{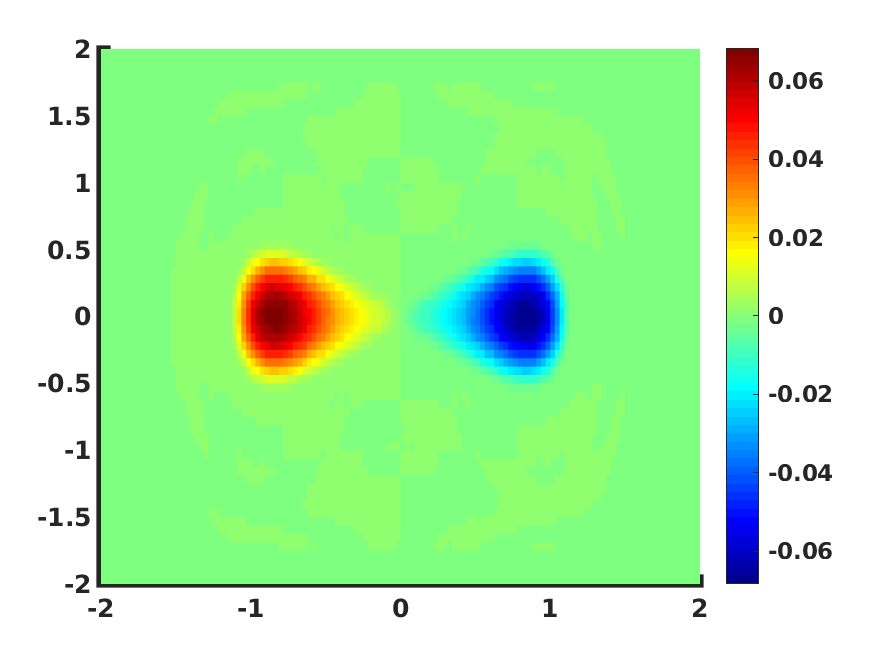} &
  \includegraphics[width=0.2\textwidth]{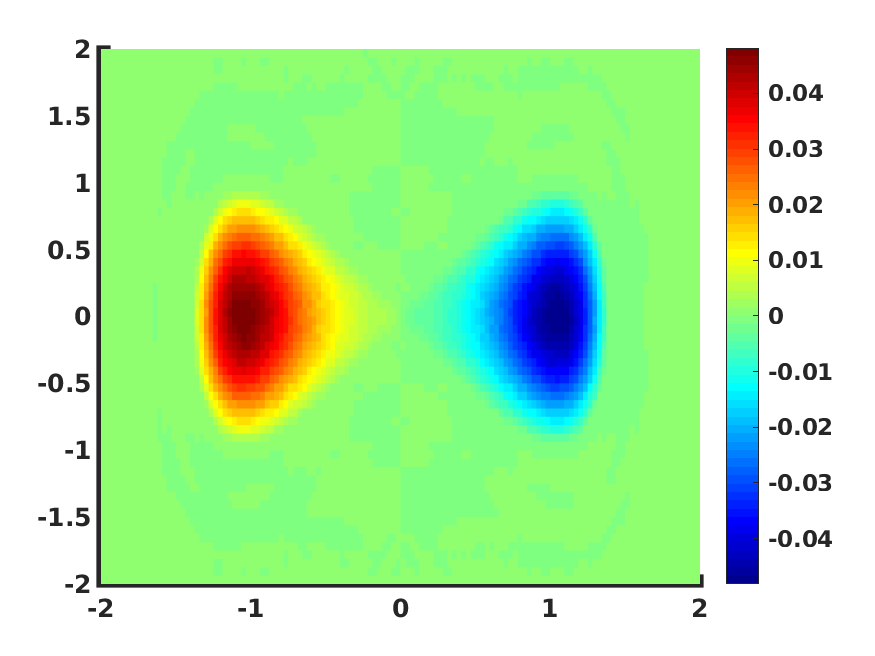} &
   \includegraphics[width=0.2\textwidth]{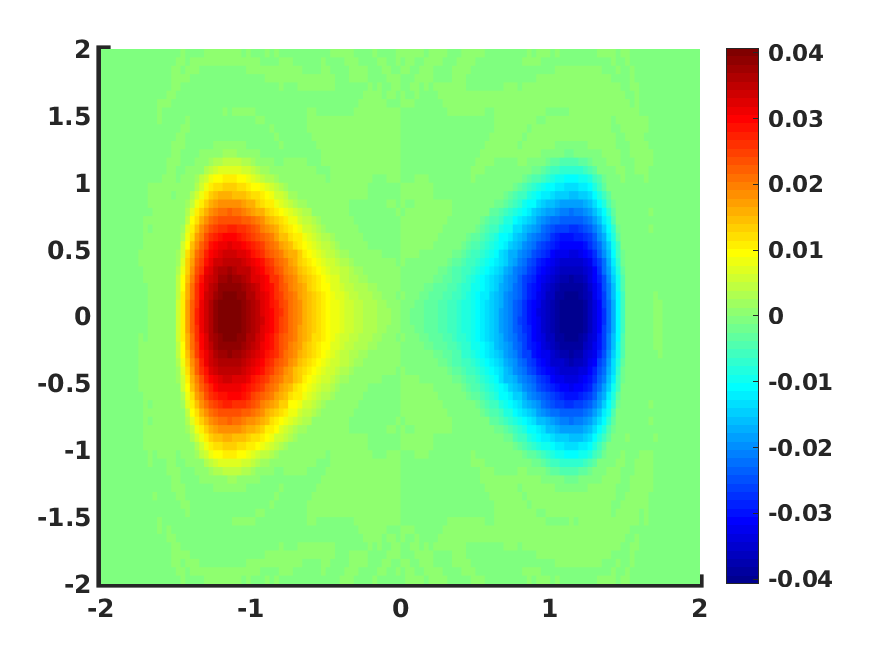} &
 \includegraphics[width=0.2\textwidth]{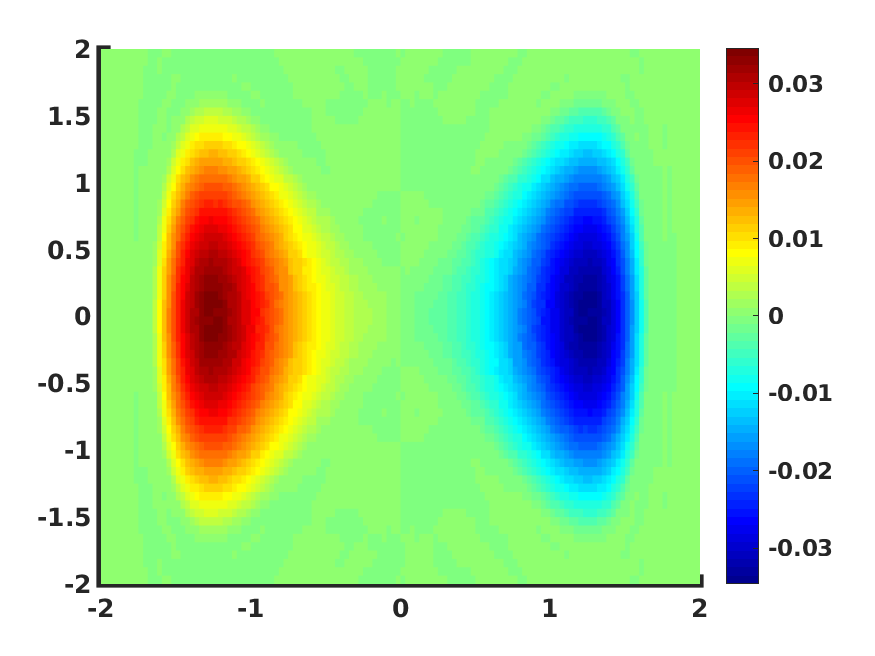}\\\end{tabular}
\end{center}
\caption{\label{fig:Dompatt} Dominant eigenvector pattern $U_1^\prime(\epsilon,p)$ post-bifurcation.}
\end{figure}

\begin{figure}[h!]
\begin{center}
 \includegraphics[width=0.3\linewidth]{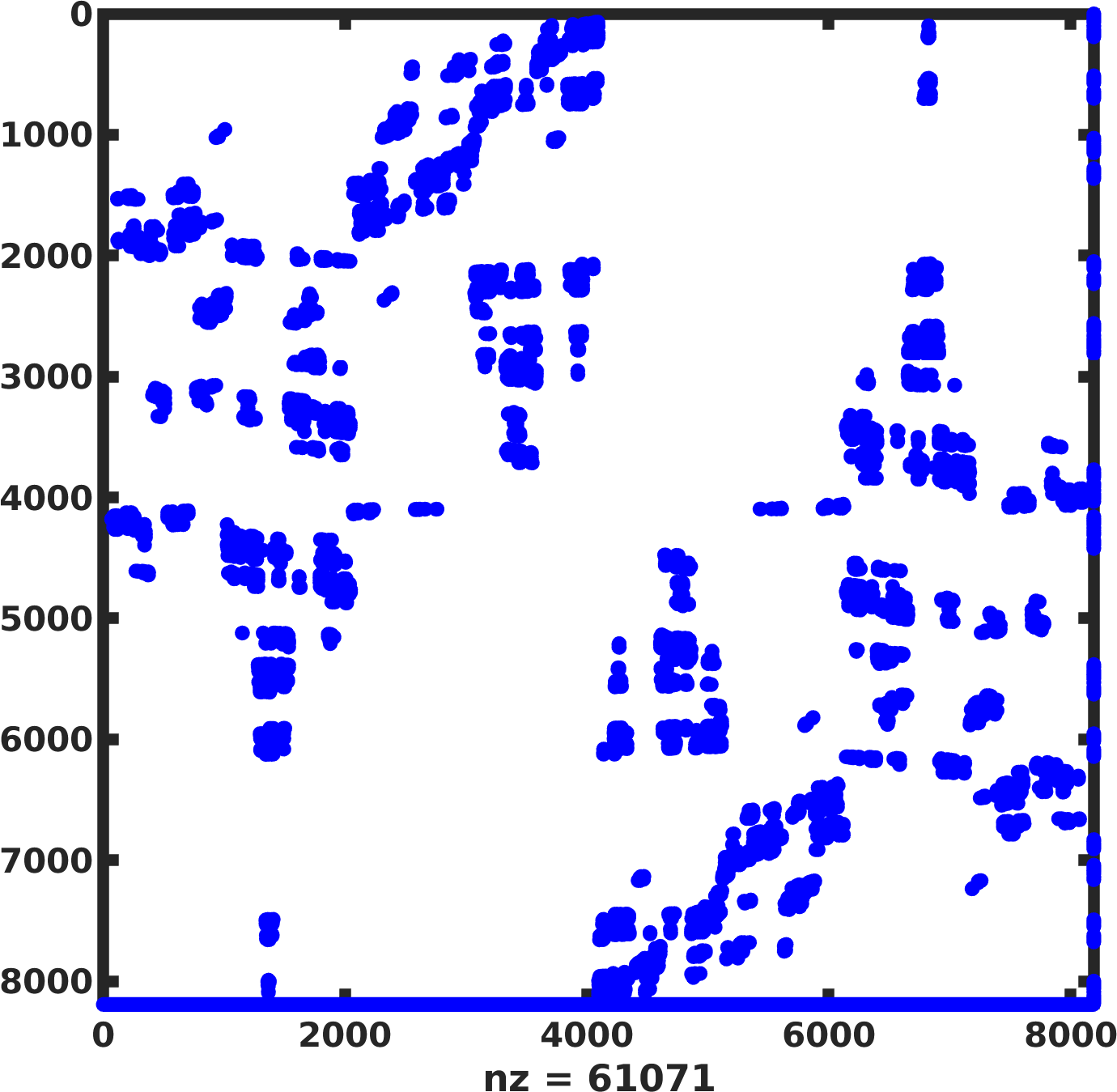}
 \end{center}
  \caption{Transition matrix post-bifurcation. }\label{fig:tmDF}
 \end{figure}

When the cascade of crossing eigenvalues (bifurcations) ends, $\lambda_1^\prime(\epsilon,p)$ becomes the second dominant eigenvalue after the eigenvalue $1$. 
Indeed, the global dynamics of 
\eqref{psys1} becomes nearly reducible with two coexisting symmetric vortices. The dominant eigenvector pattern for $p>0$ is shown in figure \ref{fig:Dompatt}.
The corresponding transition matrix is shown in figure \ref{fig:tmDF} and it is nearly reducible in accordance with the post-bifurcation global dynamics. In fact, 
one can see that the
global behavior of the system is now completely described by the support of $U_1^\prime(\epsilon,p)$.  

\subsection{Transition from one vortex to a double vortex dynamics}
\begin{figure}[!htb]
\begin{center}
\begin{tabular}{ccc}
\includegraphics[width=0.3\textwidth]{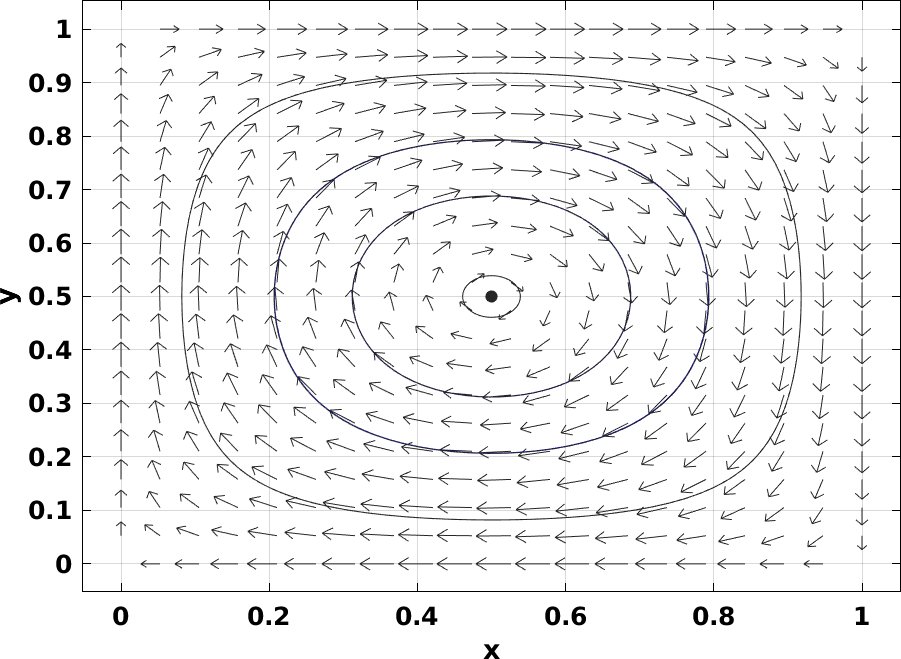} &
  \includegraphics[width=0.3\textwidth]{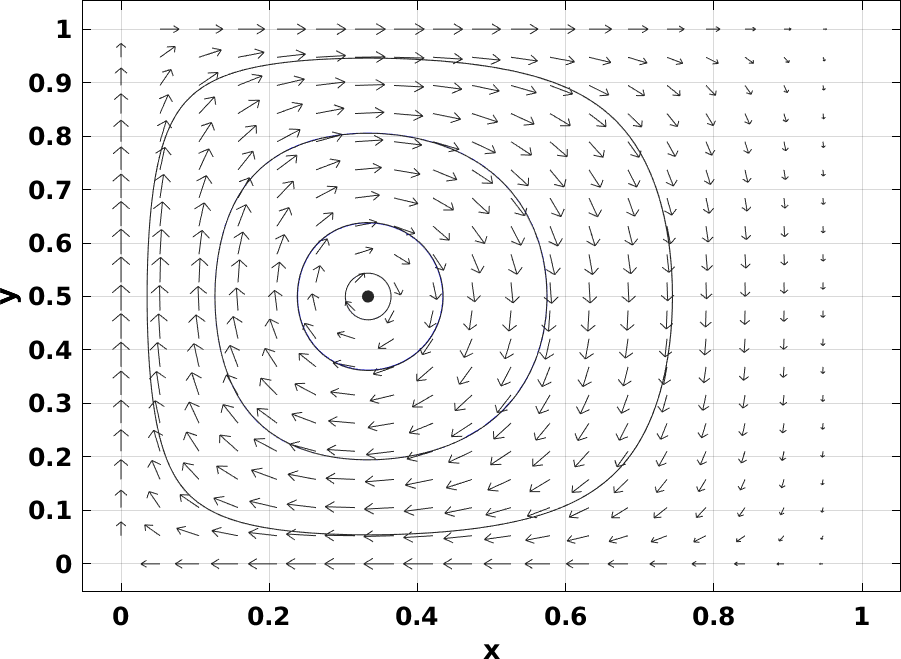} &
   \includegraphics[width=0.3\textwidth]{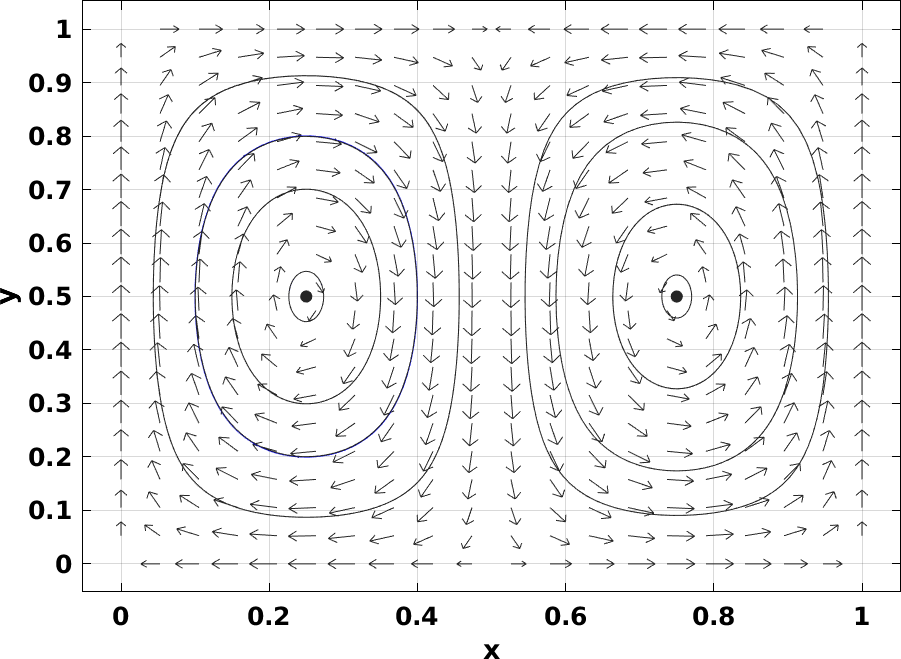} \\
   {\scriptsize $p=0$} &  {\scriptsize $p=1/3$}  &  {\scriptsize $p=1$}\\
\end{tabular}
\end{center}
\caption{\label{fig:DGvf} Changes of the velocity field in system \eqref{DG} for different $p$.}
\end{figure}

In this case study, we want to show an example of a transition of vortices that is not a bifurcation in the sense of a splitting.
The aim is to reinforce the results about the spectral indicators before a pattern splitting.  We study an incompressible two-dimensional vortex transition toy model known as the double gyre. Here, a single gyre pattern transitions to a double gyre pattern without any splitting process, which is in contrast to the setting that was studied in the previous paragraphs.
The velocity field for the system under consideration is given by 
$$V(x,t)=(-\frac{\partial\Psi}{\partial y},\frac{\partial\Psi}{\partial x})$$
with 
$$\Psi(x,y,p)=p\sin(2\pi x)\sin(\pi y)+ (1-p)\sin(\pi x) \sin(\pi y).$$
 being the parameter-dependent stream function with $p\in[0,\,1]$. We obtain the two-dimensional ordinary differential equation
\begin{equation}\label{DG}
 \begin{split}
  \dot{x}(t)   &= -(1-p)\pi\sin(\pi x)\cos(\pi y)-\pi p \sin(2\pi x)\cos(\pi x)  \\
    \dot{y}(t) &= (1-p) \pi\cos(\pi x)\sin(\pi y) + 2\pi p \cos(2 \pi x)\sin(\pi y)  \\
 \end{split}
\end{equation}
Note that the right hand side of \eqref{DG} is a convex combination of two velocity fields. For $p=0$, the dominant velocity field yields a single rotating vortex centred in 
the elliptic fixed point $(0,0)$, obtaining the system \eqref{sys2}
used earlier in section \ref{expModels} and figure \ref{fig:pEV}. For $p=1$, we have the coexistence of two counter-rotating vortices. These two values of $p$ correspond to the 
pre- and post-transition global dynamics of \eqref{DG}. The transition from a single rotating gyre to a rotating double gyre occurs for $p \in(0\,1)$, where the onset  of the emergence of the second gyre right is observed at $p=1/3$. What 
happens when $p\in (0,\,1/3)$ is that the single vortex only moves to the left side of the domain $M=[0,\,1]\times[0,\,1]$, see figure \ref{fig:DGvf}, where we illustrate the changes of the velocity field of \eqref{DG} 
with respect to $p$. Note that the motion of this single vortex to the left, 
before the transition, does not imply its expansion or shrinking. 

\begin{figure}[!htb]
  \begin{center}
  \begin{tabular}{cccc}
   \includegraphics[width=0.2\textwidth]{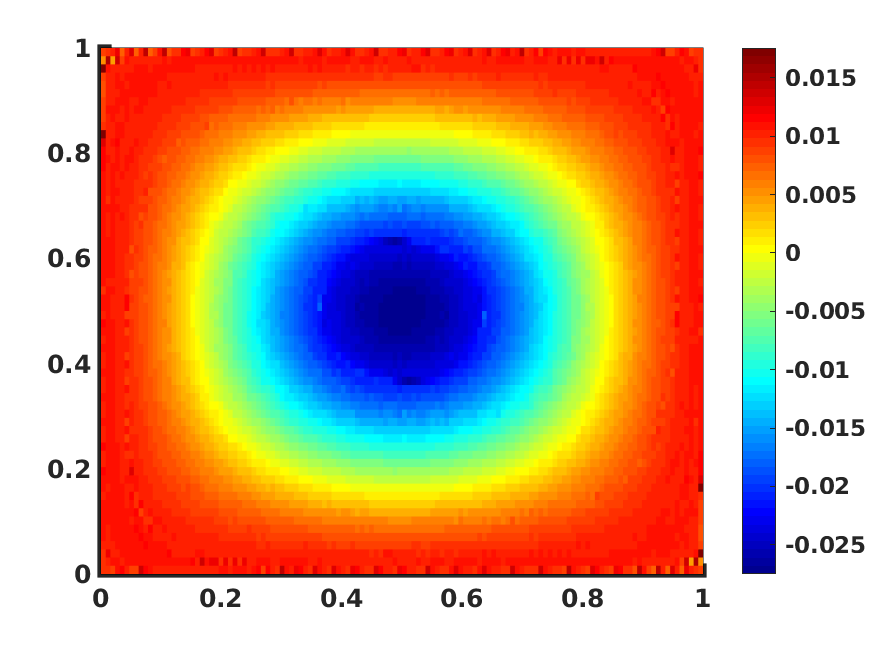}&
    \includegraphics[width=0.2\textwidth]{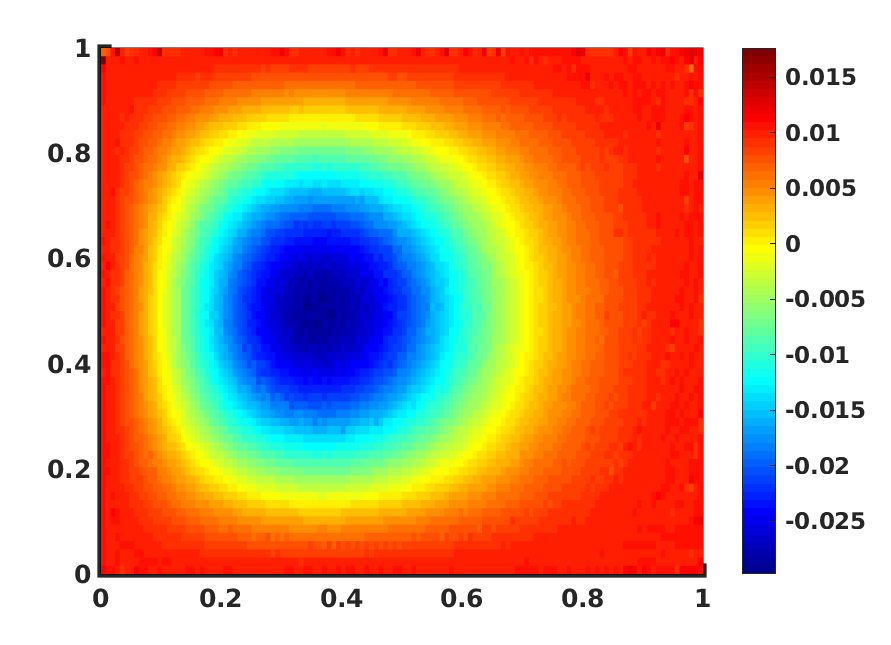}&
 \includegraphics[width=0.2\textwidth]{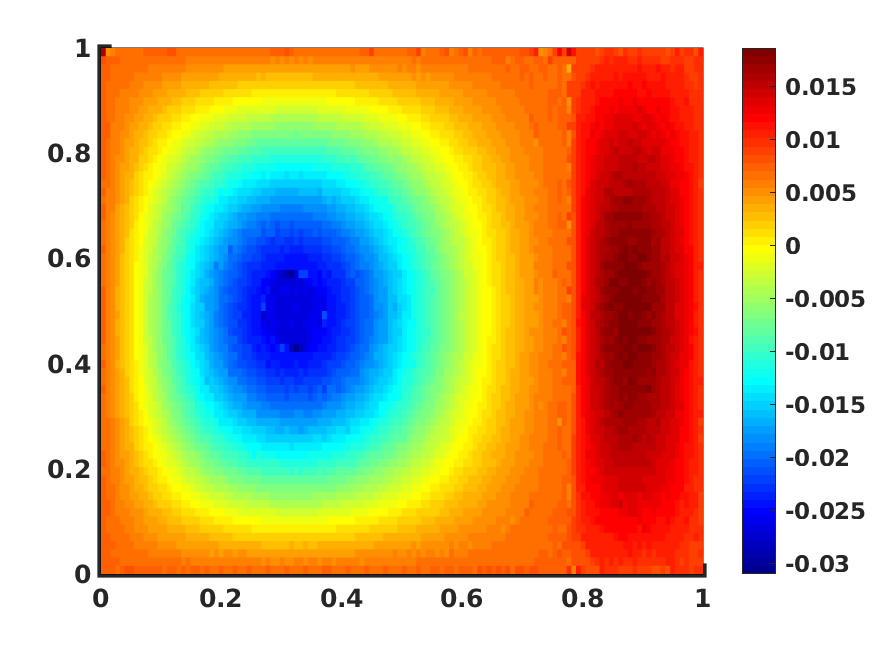}&
 \includegraphics[width=0.2\textwidth]{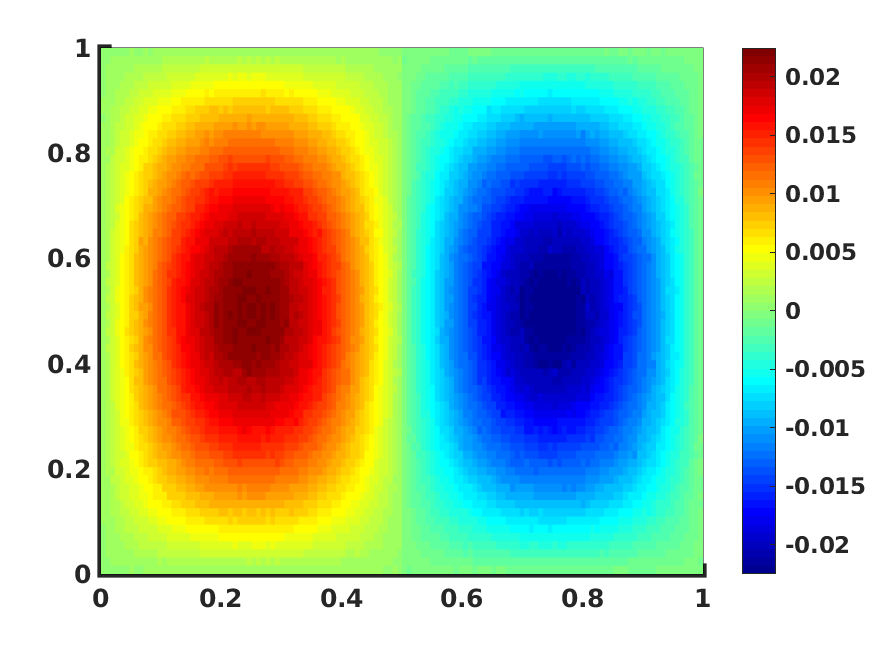}\\
  \includegraphics[width=0.2\textwidth]{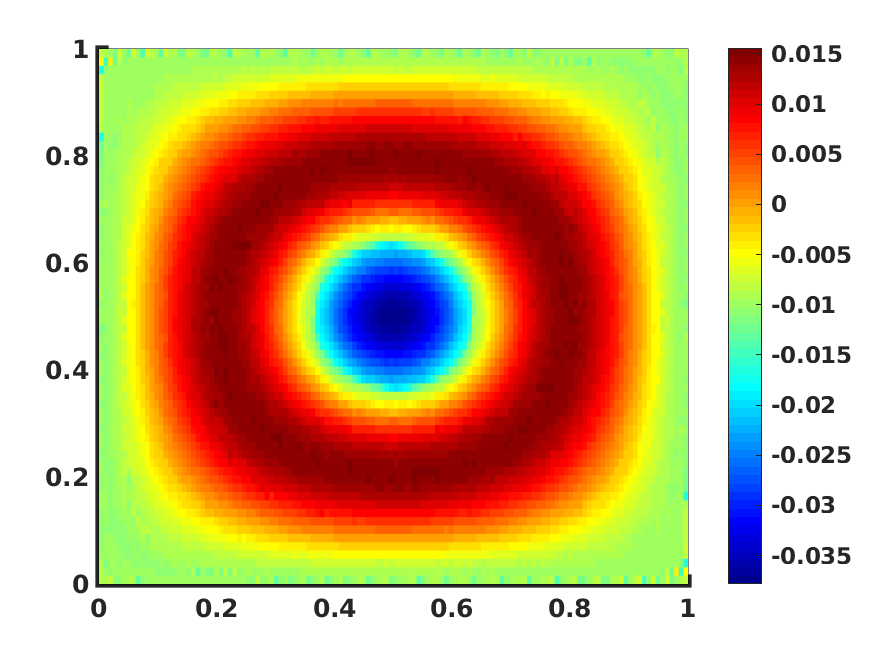}&
    \includegraphics[width=0.2\textwidth]{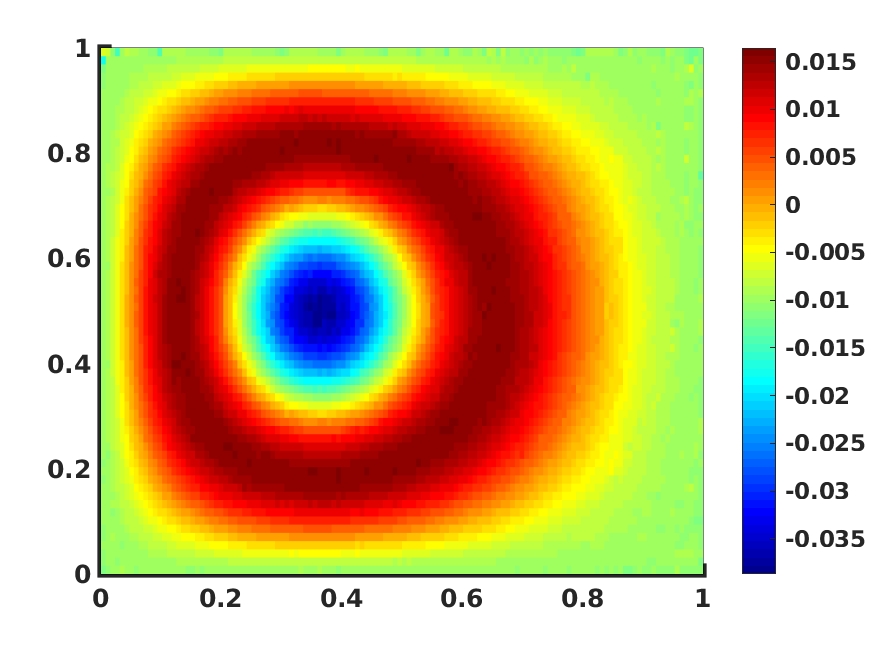}&
 \includegraphics[width=0.2\textwidth]{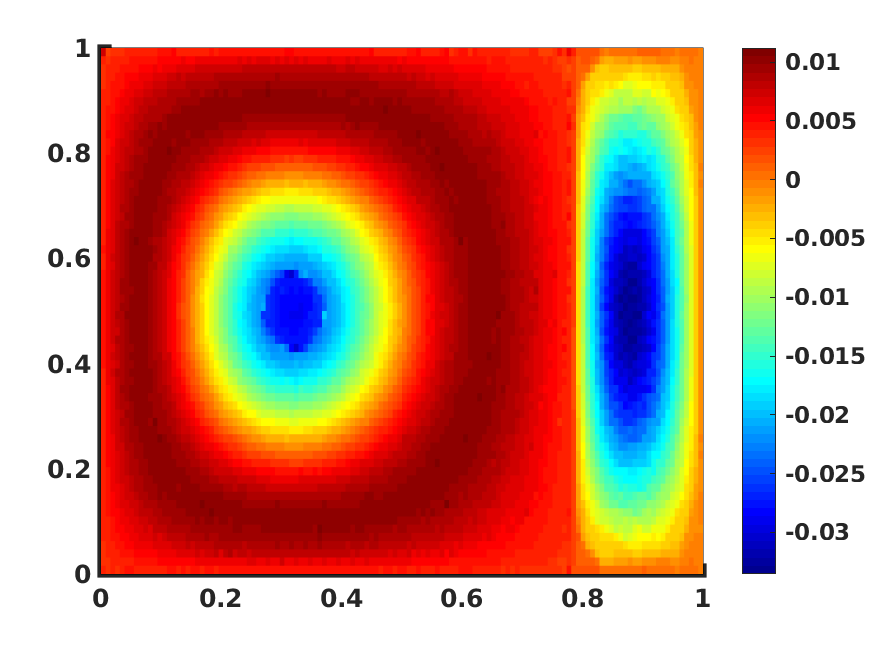}&
 \includegraphics[width=0.2\textwidth]{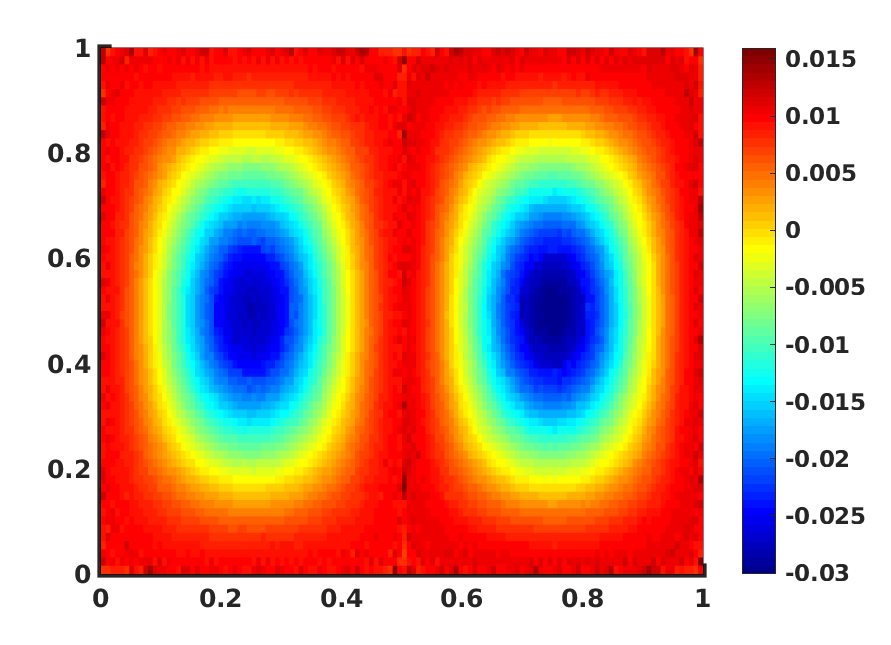}\\
  \end{tabular}
  \end{center}
\caption{\label{fig:DGu2} Changes of $U_2(\epsilon,p)$ (top) and $U_3(\epsilon,p)$ vs. $p$ .}
\end{figure}

Again, we use GAIO to numerically simulate the set-oriented dynamics of \eqref{DG} and find dominant patterns corresponding to the almost-invariant sets. 
For different values of $p$, the second and 
third dominant eigenvectors are shown in figure \ref{fig:DGu2}.

\begin{figure}[h!]
\begin{center}
 \includegraphics[width=0.6\linewidth]{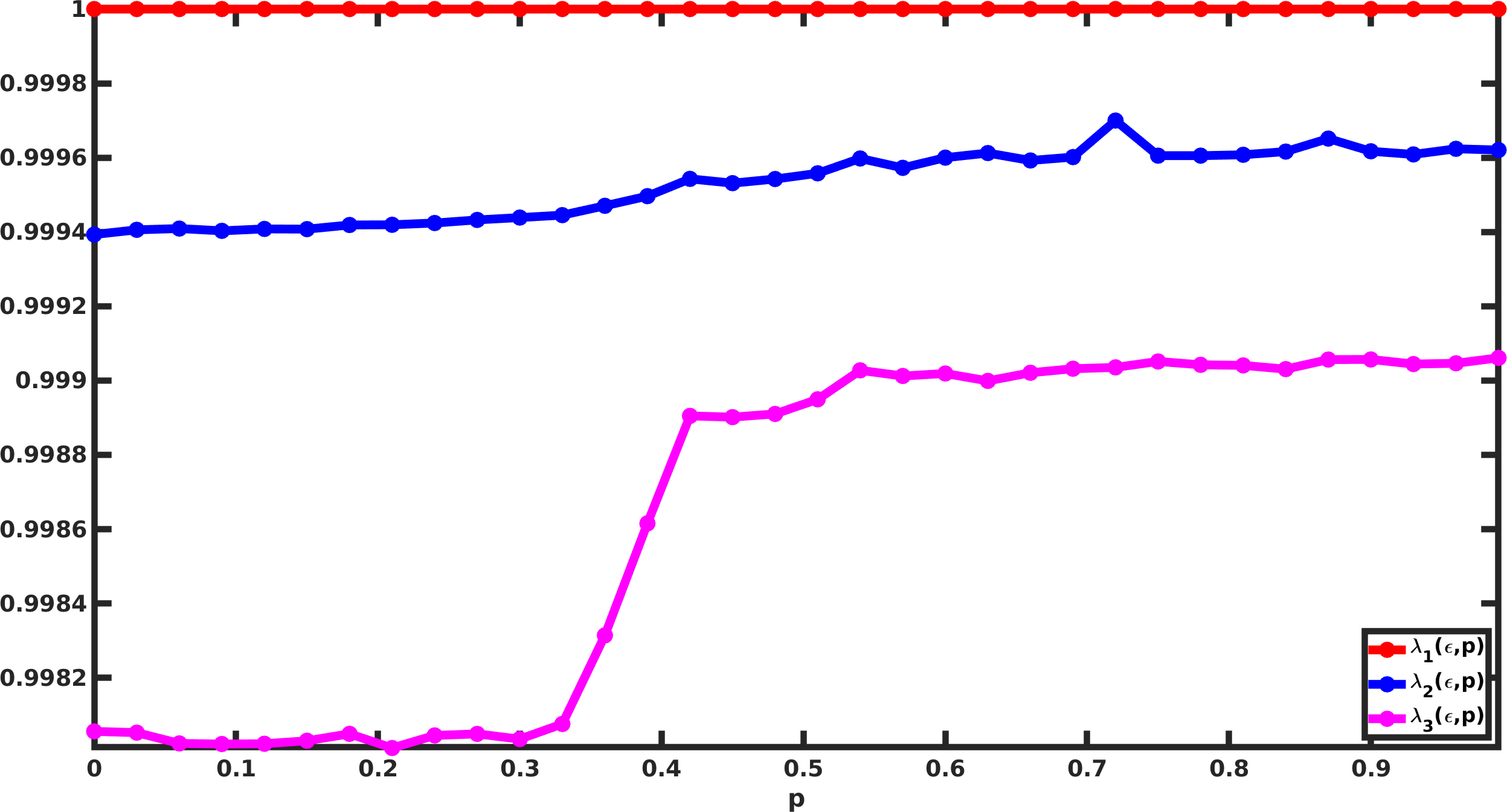}
 \end{center}
  \caption{\label{fig:DGspectrum} Spectral signature of the double gyre transition. }
 \end{figure}

The corresponding eigenvalues are shown in figure \ref{fig:DGspectrum}.
One clearly sees that the changes in the eigenvalues capture very well the behavior of the global dynamics. Indeed, for $p\in[0,\,1/3]$ 
$\lambda_2(\epsilon, p)$ and 
$\lambda_3(\epsilon,p)$ remain constant, since the single vortex only moves to the left, without shrinking or expanding. 
That explains why the eigenvalues have not decreased. For $p\in (1/3, \, 1]$, the birth of the second vortex separates the global dynamics into 
two distinct flow patterns. That is, the corresponding 
transition matrix becomes nearly reducible and it follows that the eigenvalues $\lambda_2(\epsilon,p)$ and 
$\lambda_3(\epsilon,p)$ increase to become closer to $1$. The rise of the eigenvalues $\lambda_2(\epsilon,p)$ and 
$\lambda_3(\epsilon,p)$ at $p>1/3$ can be compared to the trends of
$\lambda_1^\prime(\epsilon,p)$ in figure \ref{fig:alleigs}.

This simple transition in the dynamics of \eqref{DG} is clearly different from the critical transition caused 
by a splitting in the global dynamics of \eqref{psys1} as illustrated by figures \ref{fig:u3bifurc} and \ref{fig:u2bifurc}.
In particular, a decrease of dominant eigenvalues while another previously weak mode eigenvalue rises appears to be a spectral indicator of a splitting of almost-invariant patterns. 
\section{Conclusion}\label{conclusion}
From a set-oriented approach, we studied bifurcations of  particular almost-invariant patterns, which are supported in a neighborhood of an elliptic fixed point. These
almost-invariant sets result from invariant sets when the underlying stationary dynamical system is diffusively perturbed.  Near the splitting of patterns, generic indicators consist of a decrease of the dominant eigenvalues whose corresponding eigenvector patterns are in concern. In fact, the Duffing-type oscillator illustrates a cascade of splittings of the pattern supported in 
the neighborhood of the initially elliptic fixed point. The splitting occurs at the crossing between the dominant eigenvalues and a particular rising eigenvalue that initially belonged to the weak mode eigenvalues. It becomes the largest 
eigenvalue after the eigenvalue $1$ and its corresponding eigenvector is supported on the dominant phase space pattern post-bifurcation.

\begin{figure}[!htb]
\begin{center}
\begin{tabular}{ccc}
 \includegraphics[width=0.25\textwidth]{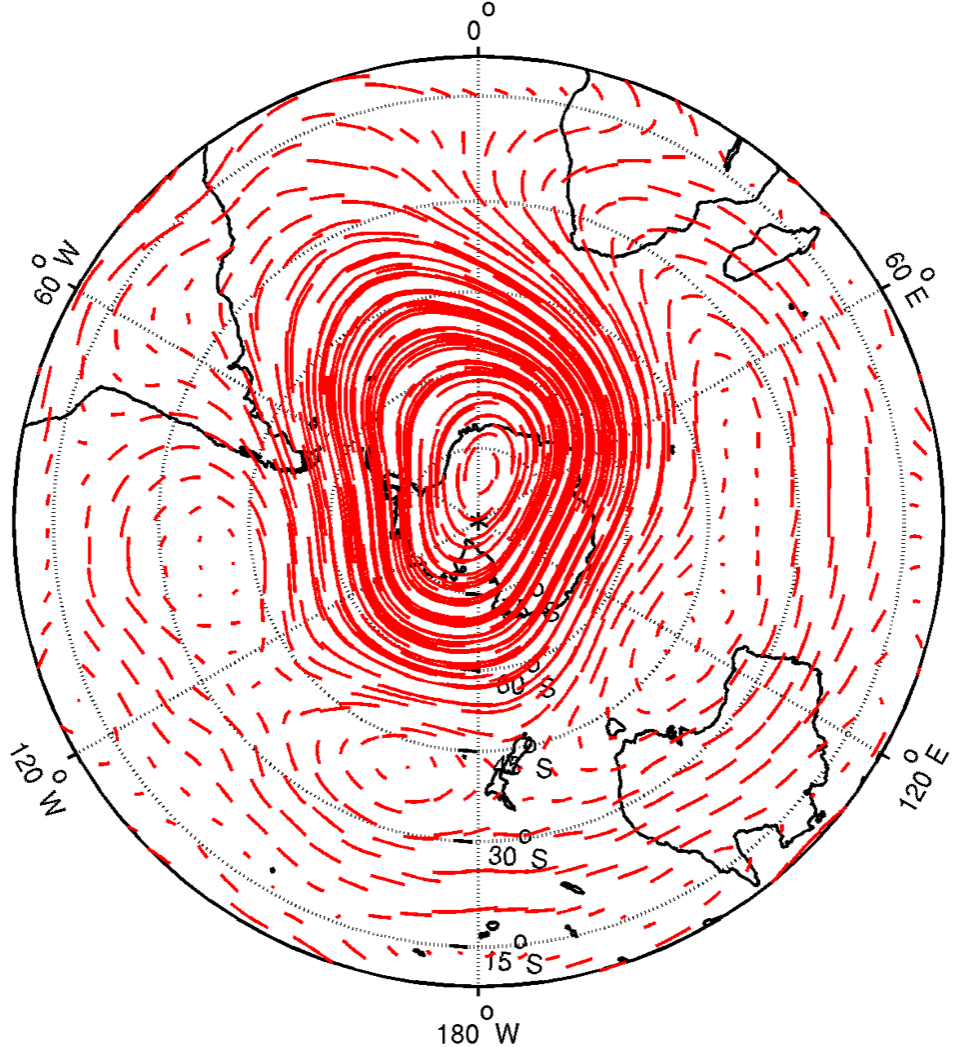} &
 \includegraphics[width=0.25\textwidth]{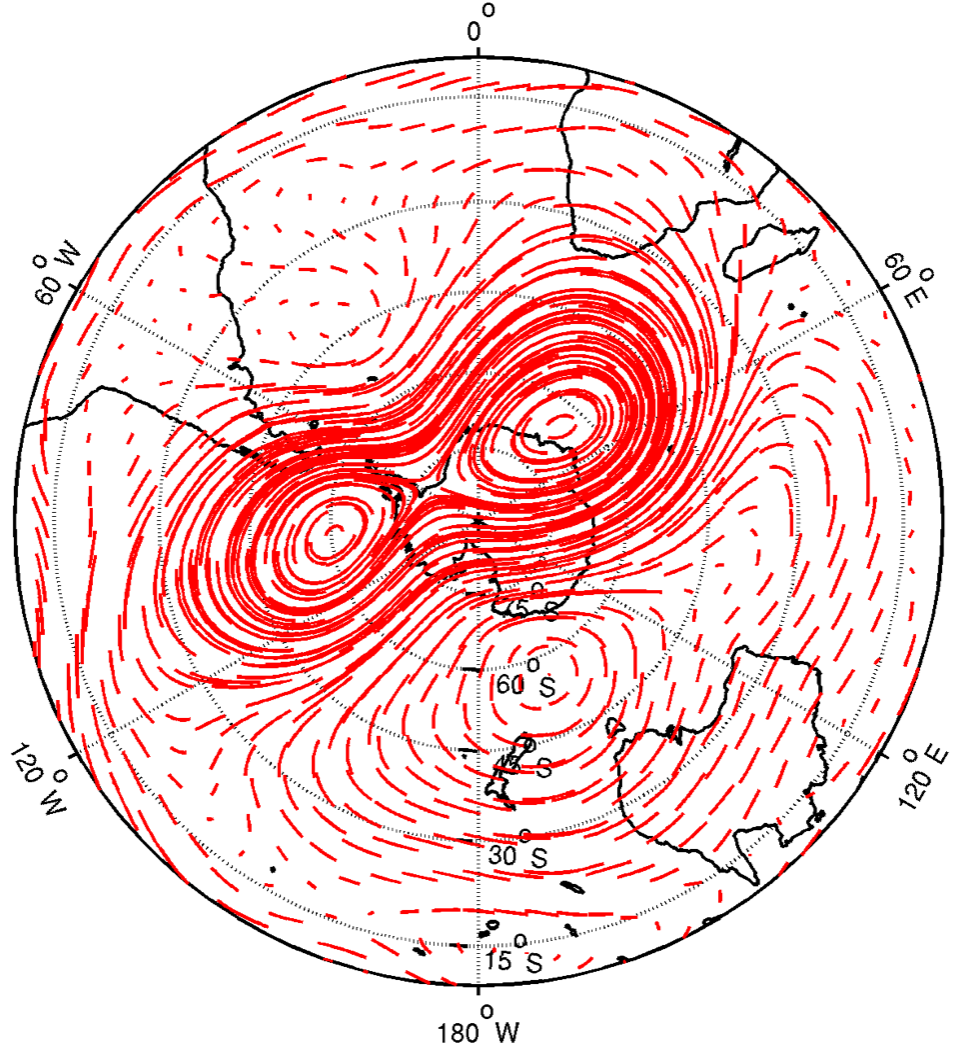} &
   \includegraphics[width=0.25\textwidth]{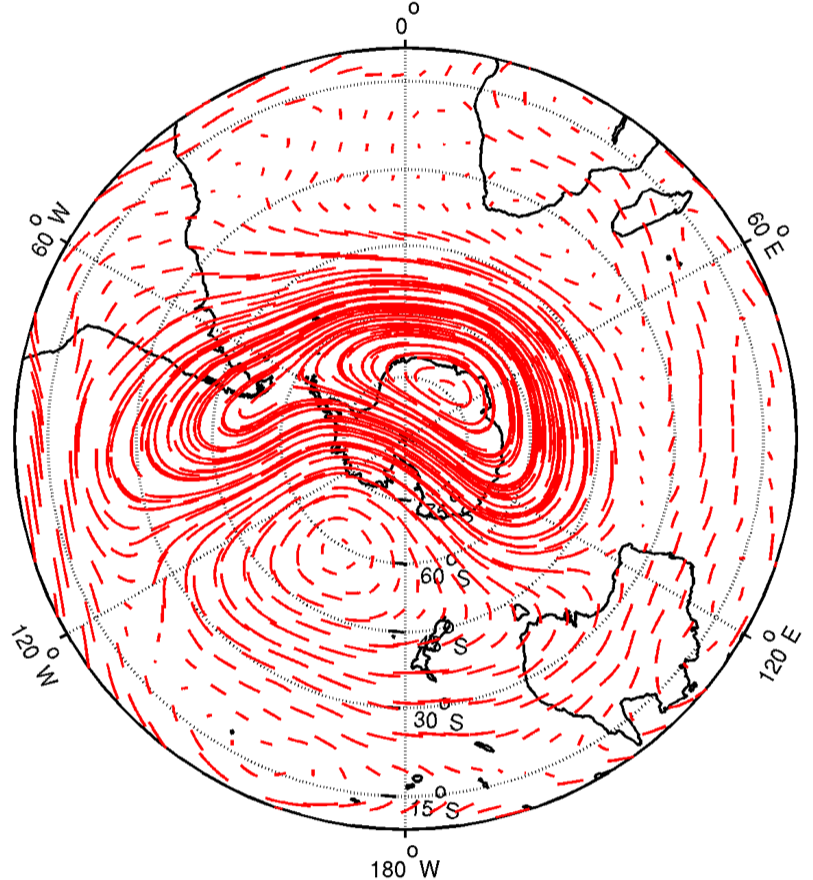}\\
   {\scriptsize September 20, 2002} &   {\scriptsize September 25, 2002} &   {\scriptsize September 30, 2002}\\
\end{tabular}
\end{center}
\caption{\label{fig:polvtex} Antarctic polar vortex splitting event in September 2002, visualized using two-dimensional velocity data from the ECMWF Interim data set (http://data.ecmwf.int/data/index.html).}
\end{figure}

Patterns emerging from complex dynamics of real- orld systems, such as the Antarctic polar vortex break up in late September 2002, suggest an 
analogous nonstationary framework. That is, in order to apply this set-oriented formulation of bifurcation analysis into real world applications, one may need to reconsider nonautonomous 
dynamical systems instead.
Thus, inspired by te present  study, future work will address a characterization of finite-time bifurcations of coherent sets, which emerge from a 
nonautonomous dynamical system. This will allow us to deduce finite-time generic early warning signals for sudden vortex splittings. These results will be used
to spectrally describe and characterize the Antarctic polar vortex splitting event from the recorded velocity data, see figure \ref{fig:polvtex}.


\end{document}